\documentclass[12pt]{article}

\usepackage{amsthm,amsmath,amsfonts,epsfig,amssymb,latexsym, bbm}
\usepackage{mathrsfs} 
\usepackage[T1]{fontenc}
\usepackage{graphicx}
\usepackage{enumerate}
\usepackage{xy,xypic}
\usepackage{graphics}
\usepackage{footnote}
\usepackage{scalefnt}
\usepackage{tikz}
\usepackage{color}
\usetikzlibrary{shapes,arrows,fit,calc,positioning}
\usetikzlibrary{matrix}
\xyoption{all}
\usepackage{imakeidx}
\makeindex

\title{Localization problems of Quillen} 
\author{
 Satya Mandal \\ University of Kansas,  Lawrence, Kansas 66045, USA
 }

\begin{document}

\pagenumbering{roman}
\setcounter{page}{0}

\renewcommand{\baselinestretch}{1.255}
\setlength{\parskip}{1ex plus0.5ex}
\date{8 October 2024}
\newcommand{\iso}{\stackrel{\sim}{\longrightarrow}}
\newcommand{\sur}{\twoheadrightarrow}

\newcommand{\eop}{\hfill \TCP{\rule{2mm}{2mm}}}
\newcommand{\pf}{\noindent{\bf Proof.~}}
\newcommand{\outlinePf}{\noindent{\bf Outline of the Proof.~}}

\newcommand{\PD}{\dim_{{\SA}}}
\newcommand{\PDV}{\dim_{{\SV(X)}}}

\newcommand{\ra}{\rightarrow}
\newcommand{\lra}{\longrightarrow}
\newcommand{\hra}{\hookrightarrow} 
\newcommand{\Lra}{\Longrightarrow}
\newcommand{\Lla}{\Longleftarrow}
\newcommand{\Llra}{\Longleftrightarrow}
\newcommand{\pic}{The proof is complete.~}
\newcommand{\Sp}{\mathrm{Sp}}
\newcommand{\BiSp}{\mathrm{BiSp}}
\newcommand{\colim}{\mathrm{colim}}
\newcommand{\codim}{\mathrm{co}\dim}

\newcommand{\Dia}{\diagram}

\newcommand{\bE}{\begin{enumerate}}
\newcommand{\eE}{\end{enumerate}}

\newtheorem{theorem}{Theorem}[section]
\newtheorem{proposition}[theorem]{Proposition}
\newtheorem{lemma}[theorem]{Lemma}
\newtheorem{definition}[theorem]{Definition}
\newtheorem{corollary}[theorem]{Corollary}
\newtheorem{construction}[theorem]{Construction}
\newtheorem{notation}[theorem]{Notation}
\newtheorem{notations}[theorem]{Notations}
\newtheorem{remark}[theorem]{Remark}
\newtheorem{question}[theorem]{Question}
\newtheorem{example}[theorem]{Example} 
\newtheorem{examples}[theorem]{Examples} 
\newtheorem{exercise}[theorem]{Exercise} 
\newtheorem{setup}[theorem]{Setup} 
\newtheorem{openProblem}[theorem]{Open Problem}

\newtheorem{problem}[theorem]{Problem} 
\newtheorem{conjecture}[theorem]{Conjecture} 

\newcommand{\bD}{\begin{definition}}
\newcommand{\eD}{\end{definition}}
\newcommand{\bP}{\begin{proposition}}
\newcommand{\eP}{\end{proposition}}
\newcommand{\bL}{\begin{lemma}}
\newcommand{\eL}{\end{lemma}}
\newcommand{\bT}{\begin{theorem}}
\newcommand{\eT}{\end{theorem}}
\newcommand{\bC}{\begin{corollary}}
\newcommand{\eC}{\end{corollary}} 
\newcommand{\TCP}{\textcolor{purple}}
\newcommand{\TCM}{\textcolor{magenta}}
\newcommand{\TCR}{\textcolor{red}}
\newcommand{\TCB}{\textcolor{blue}}
\newcommand{\TCG}{\textcolor{green}}
\def\spec#1{\mathrm{Spec}\left(#1\right)}
\def\proj#1{\mathrm{Proj}\left(#1\right)}
\def\supp#1{\mathrm{Supp}\left(#1\right)}
\def\Sym#1{{\CS}\mathrm{ym}\left(#1\right)} 
\def\norm#1{\parallel #1\parallel}
\def\LRf#1{\left( #1\right)}
\def\LRs#1{\left\{ #1\right\}}
\def\LRt#1{\left[ #1\right]}
\def\Cat{\underline{\bf Cat}}
\def\eCat{\underline{\bf eCat}}
\def\dgCat{\underline{\bf dgCat}}

\def\a{\mathfrak {a}}
\def\b{\mathfrak {b}}
\def\c{\mathfrak {c}}

\def\d{\mathfrak {d}}
\def\e{\mathfrak {e}}
\def\f{\mathfrak {f}}
\def\g{\mathfrak {g}}
\def\i{\mathfrak {i}}
\def\j{\mathfrak {j}}
\def\k{\mathfrak {k}}
\def\l{\mathfrak {l}}
\def\m{\mathfrak {m}}
\def\n{\mathfrak {n}}
\def\p{\mathfrak {p}}
\def\q{\mathfrak {q}}
\def\r{\mathfrak {r}}
\def\s{\mathfrak {s}}
\def\t{\mathfrak {t}}
\def\u{\mathfrak {v}}
\def\w{\mathfrak {w}}

\def\A{\mathfrak {A}}
\def\B{\mathfrak {B}}
\def\C{\mathfrak {C}}
\def\D{\mathfrak {D}}
\def\E{\mathfrak {E}}

\def\F{\mathfrak {F}}

\def\G{\mathfrak {G}}
\def\H{\mathfrak {H}}
\def\I{\mathfrak {I}}
\def\J{\mathfrak {J}}
\def\K{\mathfrak {K}}
\def\L{\mathfrak {L}}
\def\M{\mathfrak {M}}
\def\N{\mathfrak {N}}
\def\O{\mathfrak {O}}
\def\P{\mathfrak {P}}

\def\Q{\mathfrak {Q}}

\def\R{\mathfrak {R}}

\def\Sf{\mathfrak {S}} 
\def\T{\mathfrak {T}}
\def\U{\mathfrak {U}}
\def\V{\mathfrak {V}}
\def\W{\mathfrak {W}}

\def\CA{\mathcal {A}}
\def\CB{\mathcal {B}}
\def\CP{\mathcal {P}}
\def\CC{\mathcal {C}}
\def\CD{\mathcal {D}}
\def\CE{\mathcal {E}}
\def\CF{\mathcal {F}}
\def\CG{\mathcal {G}}
\def\CH{\mathcal{H}}
\def\CI{\mathcal{I}}
\def\CJ{\mathcal{J}}
\def\CK{\mathcal{K}}
\def\CL{\mathcal{L}}
\def\CM{\mathcal{M}}
\def\CN{\mathcal{N}}
\def\CO{\mathcal{O}}
\def\CP{\mathcal{P}}
\def\CQ{\mathcal{Q}}

\def\CR{\mathcal{R}}

\def\CS{\mathcal{S}}
\def\CT{\mathcal{T}}
\def\CU{\mathcal{U}}
\def\CV{\mathcal{V}}
\def\CW{\mathcal{W}} 
\def\CX{\mathcal{X}}
\def\CY{\mathcal{Y}}
\def\CZ{\mathcal{Z}}

\newcommand{\smallcirc}[1]{\scalebox{#1}{$\circ$}}

\def\BA{\mathbb{A}}
\def\BB{\mathbb{B}}
\def\BC{\mathbb{C}}
\def\BD{\mathbb{D}}
\def\BE{\mathbb{E}}
\def\BF{\mathbb{F}}
\def\BG{\mathbb{G}}
\def\BH{\mathbb{H}}
\def\BI{\mathbb{I}}
\def\BJ{\mathbb{J}}
\def\BK{\mathbb{K}}
\def\BL{\mathbb{L}}
\def\BM{\mathbb{M}}
\def\BN{\mathbb{N}} 
\def\BO{\mathbb{O}}
\def\BP{\mathbb{P}} 
\def\BQ{\mathbb{Q}}
\def\BR{\mathbb{R}}
\def\BS{\mathbb{S}}
\def\BT{\mathbb{T}}
\def\BU{\mathbb{U}}
\def\BV{\mathbb{V}}
\def\BW{\mathbb{W}}
\def\BX{\mathbb{X}}
\def\BY{\mathbb{Y}}
\def\BZ{\mathbb {Z}}

\def\SA{\mathscr {A}} 
\def\SB{\mathscr {B}}
\def\SC{\mathscr {C}}
\def\SD{\mathscr {D}}
\def\SE{\mathscr {E}}
\def\SF{\mathscr {F}}
\def\SG{\mathscr {G}}
\def\SH{\mathscr {H}}
\def\SI{\mathscr {I}}
\def\SK{\mathscr {K}} 
\def\SL{\mathscr {L}} 
\def\SM{\mathscr {M}}
\def\SN{\mathscr {N}}
\def\SO{\mathscr {O}}
\def\SP{\mathscr {P}}
\def\SQ{\mathscr {Q}}
\def\SR{\mathscr {R}}
\def\SS{\mathscr {S}}
\def\ST{\mathscr {T}}
\def\SU{\mathscr {U}}
\def\SV{\mathscr {V}}
\def\SW{\mathscr {W}}
\def\SX{\mathscr {X}}
\def\SY{\mathscr {Y}}
\def\SZ{\mathscr {Z}} 
 
\def\bfA{\bf A}
\def\bfB{\bf B}
\def\bfC{\bf C}
\def\bfD{\bf D}
\def\bfE{\bf E}
\def\bfF{\bf F}
\def\bfG{\bf G}
\def\bfH{\bf H}
\def\bfI{\bf I}
\def\bfJ{\bf J}
\def\bfK{\bf K}
\def\bfL{\bf L}
\def\bfM{\bf M}
\def\bfN{\bf N} 
\def\bfO{\bf O}
\def\bfP{\bf P} 
\def\bfQ{\bf Q}
\def\bfR{\bf R}
\def\bfS{\bf S}
\def\bfT{\bf T}
\def\bfU{\bf U}
\def\bfV{\bf V}
\def\bfW{\bf W}
\def\bfX{\bf X}
\def\bfY{\bf Y}
\def\bfZ{\bf Z}

\maketitle









\pagenumbering{arabic}
 \noindent{\bf Abstract:} Let $X$ be a quasi projective scheme over a noetherian affine scheme $Spec(A)$, $U\subseteq X$ be an open subset, and $Z=X-U$.Assume that $Z$ is complete intersection, with $k=\codim Z$. Consider the map   
  $$
  \q:{\mathbb K}\left({\mathscr V}(X)\right) \rightarrow {\mathbb K}\left({\mathscr V}(U)\right)
  $$
 of  the  ${\mathbb K}$-theory spectra. We give a description of the homotopy fiber of $\q$. 
 Let  $C{\mathbb M}^Z\left(X\right)$ denote the full subcategory of perfect modules ${\mathscr F} \in Coh(X)$ such that(1) ${\mathscr F} _{|U}=0$, (2) $grade({\mathscr F} )=\dim_{{\mathscr V}(X)}{\mathscr F}=k $.
 It turns out that the homotopy fiber of $\q$ is the ${\mathbb K}$-theory spectra 
 ${\mathbb K}\left(C{\mathbb M}^Z\left(X\right)\right)$. 
 Likewise, we compute the homotopy fiber of the   pullback map
 $$
 \g:  {\mathbb G}W\left({\mathscr V}(X)\right) \rightarrow {\mathbb G}W\left({\mathscr V}(U)\right) 
 $$
 of Karoubi Grothendieck-Witt bispectra. Consequently, we obtain long exact sequences of ${\mathbb K}$-groups and of ${\mathbb G}W$-groups. These results settle some of the long standing open problems. 
\section{Introduction} 

In general, refer to section \ref{SecVorAus} for notations. Throughout, $GW$ abbreviates 
"Grothedieck-Witt", with possible variations of  fonts. For notations and background on $GW$ theory, refer to section \ref{BISP4Exact}.  Whenever we discuss $GW$ theory (also known as hermitian theory), we assume $2$ is invertible in the respective category. 
Section \ref{SecGrade} establishes some background on grade and annihilators of coherent modules over quasi projective schemes. 
All schemes considered in this article are assumed to be noetherian.

In this article we discuss the homotopy fiber of the  pullback 
maps
\begin{equation}\label{372First}
\left\{\begin{array}{ll}
\diagram 
{\BK}\left({\SV}(X)\right) \ar[r]^{\q} & {\BK}\left({\SV}(U)\right)\\
\enddiagram & in ~\Sp\\
\diagram 
{\BG}W\left({\SV}(X)\right) \ar[r]_{\p} & {\BG}W\left({\SV}(U)\right)\\
\enddiagram & in~\BiSp\\
\end{array}\right.
\end{equation}
where $X$ is a quasi projective scheme, over a noetherian affine scheme $\spec{A}$, and $U\subseteq X$ is an open subset. Here ${\BK}\left(-\right)\in {\Sp}$ denotes the nonconnective ${\BK}$-theory spectra and ${\BG}W(-)\in \BiSp$ denotes the Karoubi ${\BG}W$ bispectra (see section \ref{BISP4Exact}). 
%
We will write $Z=X-U$. Recall the classical fundamental result that, any map $\f:{\CY} \lra {\CZ}$ of pointed spaces fits into a homotopy fibration $\diagram {\SF}(\f) \ar[r] &{\CY} \ar[r]^{\f} & {\CZ} \enddiagram$ leading up to a long exact sequence of homotopy groups, ending at the degree zero term (see \cite{W78}). Therefore, the question of existence of the homotopy fibers of $\q$, $\p$ was never in contention. Right from the inception 
\cite{Q73}, it has been an open 
question, how best to describe the homotopy fiber of the maps $\q$ ({\it originally, of the maps of $K$-theory spaces,  then of spectra and so on}), in terms of the closed complement $Z$. The questions on homotopy fiber of $\p$ followed, as the $GW$ theory evolved. For the  ${\bfG}$-theory spaces ({\it meaning, $K$-theory spaces of $Coh(X)$}), Quillen gave fully satisfactory description of the homotopy fibers. For regular schemes,  $K$-theory coincides with ${\bfG}$-theory. Therefore, the issue of such homotopy fibers remains to be interesting only when $X$ is non regular.

This project was triggered when the  result (\ref{GrayQuillen}) of Quillen came to my attention. 

 \bT[Quillen  \cite{G74, G76}]\label{GrayQuillen}{\rm 
Suppose $X$ is a quasi compact scheme,  $U=\spec{A}\subseteq X$ is an affine open subschemes of $X$, and $Z=X-U$. Assume that $Z$ has a subscheme structure,  defined by an invertible ideal sheaf ${\SI}\subseteq {\CO}_X$.
Let 
$$
{\SH}=\left\{{\SF}\in QCoh(X): {\SF}_{|U}=0, \PDV({\SF})=1 \right\} \qquad {\rm be~the~full~subcategory}. 
$$
Then there is an exact sequence 
$$
\diagram 
\cdots \ar[r] &
{K}_n\left({\SH}\right) \ar[r] &{K}_n\left({\SV}(X)\right) \ar[r] & {K}_n\left({\SV}(U)\right)
\ar[r] &{K}_{n-1}\left({\SH} \right)\ar[r]& \cdots\\
\enddiagram 
$$ 
ending at degree zero.
}
\eT

It was instantly clear to me, while $\codim(X, Z)=1$,  the category ${\SH}=C{\BM}^Z(X)$ is the category of perfect modules with support on $Z$. The $K$-theory of various categories of perfect modules were previously studied in \cite{M15, M17, M20, M23} and others. 

After  theorem \ref{GrayQuillen}, it was clear to the experts what to expect, in higher codimension.
Following formulation was communicated to me by Daniel Grayson,  which we write as an open problem.
\begin{openProblem}\label{opnCMZ}{\rm
Let $X$ be a quasi projective scheme  over a noetherian affine scheme $\spec{A}$ and $Z\subseteq X$ be a closed complete intersection subset, as in
(\ref{CIsetUp}) and $U=X-Z$.  ({\it One says, $Z$ is cutout by a regular sequence.}) Then the sequence
$$
\diagram
{\bfK}\LRf{{\BM}^Z(X)} \ar[r] & {\bfK}\LRf{{\SV}(X)} \ar[r]^{\q} & {\bfK}\LRf{{\SV}(U)} \\
\enddiagram
$$
of $K$-theory spaces, is a homotopy fibration. By now, we think about this problem in terms of nonconnective ${\BK}$-theory spectrum. 
}
\end{openProblem}
In the case, $\codim(Z)=1$ the problem was stated as conjectures, in the paper of Gersten \cite[pp 343]{G74}
and in the paper of Weibel \cite[pp 199]{W89}.

Using the methods developed on various categories of perfect modules ({\it loc. cit}), we 
settle this long standing open problem (\ref{opnCMZ}), in complete generality. 
We prove when $Z$ is complete intersection, as in (\ref{CIsetUp}),  the homotopy fiber of the map $\q$ is given by the ${\BK}$-theory spectra 
${\BK}\left(C{\BM}^Z(X)\right)\in {\Sp}$ (see (\ref{defPerModd}, \ref{nota},) for notations). 
Analogously, the homotopy fiber 
of $\p$ is given by ${\BG}W\left(C{\BM}^Z(X)\right)\in {\BiSp}$. In other words (\ref{Equi1909ALL}, \ref{BiSpectrEHmoFib}),
there are homotopy fibrations 
\begin{equation}\label{372Firstresulatr}
\left\{\begin{array}{ll}
\diagram 
{\BK}\left(C{\BM}^Z(X)\right) \ar[r] &
{\BK}\left({\SV}(X)\right) \ar[r]^{\q} & {\BK}\left({\SV}(U)\right)\\
\enddiagram & in ~\Sp\\
\diagram 
{\BG}W^{[-k]}\left(C{\BM}^Z(X)\right) \ar[r] &
{\BG}W\left({\SV}(X)\right) \ar[r]_{\p} & {\BG}W\left({\SV}(U)\right)\\
\enddiagram & in~\BiSp\\
\end{array}\right.
\end{equation}
where  $k=\codim(Z)\geq 1$ and $Z$ is complete intersection. Note, we have
$$
{\BG}W^{[-k]}\left(C{\BM}^Z(X)\right) =
\left\{\begin{array}{ll}
{\BG}W\left(C{\BM}^Z(X)\right) &if~k=0~mod~4\\
{\BG}W\left(C{\BM}^Z(X)\right) ^-&if~k=2~mod~4\\
\end{array}\right.
$$
Further, 
we  prove  
\begin{equation}\label{372TWOesulatr}
{\BK}\left(C{\BM}^Z(X)\right) \iso {\BK}\left({\BM}^Z(X)\right) \quad {\rm is~a~homotopy~equivalence~in}~{\Sp}.
\end{equation} 

  Note, 
  if  $X$ is noetherian and regular then ${\BM}(X)=Coh(X)$, ${\BM}^Z(X)=Coh^Z(X)$. 
  Further note that,  statements analogous to (\ref{372TWOesulatr}) in $GW$-theory would not even make sense, because ${\BM}^Z(X)$ does not have a natural duality. 
 However, $C{\BM}^Z(X)$ has  natural dualities $M \mapsto {\SE}xt^k\LRf{M, {\SL}}$, when $k=grade(Z, X)$ and ${\SL}$ is any line bundle on $X$. Because of this, our results on $GW$-theory are sharper. 

As a consequence to (\ref{372Firstresulatr}), if $Z$ is complete intersection,  $\forall n\in {\BZ}$, there are exact sequences of corresponding groups
\begin{equation}\label{MainOur381}
\left\{\begin{array}{l}
\diagram 
{\BK}_n\left(C{\BM}^Z(X)\right) \ar[r] &{\BK}_n\left({\SV}(X)\right) \ar[r] & {\BK}_n\left({\SV}(U)\right)
\ar[r] &{\BK}_{n-1}\left(C{\BM}^Z(X)\right)\\
\enddiagram \\
\diagram 
{\BG}W_n\left(C{\BM}^Z(X)\right)^{\pm} \ar[r] &{\BG}W_n\left({\SV}(X)\right) \ar[r] & 
{\BG}W_n\left({\SV}(U)\right)
\ar[r] &{\BG}W_{n-1}\left(C{\BM}^Z(X)\right)^{\pm}\\
\enddiagram \\
\end{array}\right. 
\end{equation}
while we assume $\codim Z$ is even for the latter exact sequence to make sense. 
%

%
%
  Thomason developed a complete throry \cite{TT90}, in terms of $K$-theory of perfect complexes. 
After the
 work of Waldhausen \cite{W83}, followed by  \cite{TT90}, $K$-theory spaces and spectra became  invariants of the corresponding chain complex categories. Eventually, it percolated to the $K$-theory of perfect complexes \cite{TT90}. The homotopy fiber of $\q$ is described as $K$-theory of
${\bf P}erf_Z(X)$ in  \cite{TT90} ({\it see \S ~\ref{SecVorAus} for notation}).
This description of the homotopy fiber of $\q$ in  \cite{TT90} was label as a "revolutionary advance"  \cite[pp 337]{TT90} compared to Theorem \ref{GrayQuillen}. 
However,
%
I believe that $K$-theory should be a story of exact categories, which is what we do in this article. 

We insert the following example, for perspective. 
\begin{example}[Deligne \cite{G74}]\label{DelineExample}{\rm 
Let $(A, \m)$ be a local noetherian commutative ring, which is not Cohen Macaulay. 
Let $X=\spec{A}$,  $Z=V(\m)$ and $U=X-Z$. In this case,  $C{\BM}^Z(X)=\{0\}$. 
There is an example of a non Cohen Macaualay ring $A$, suggested by Deligne \cite{G74},  such that ${\bf K}\LRf{{\SV}(X) }\neq {\bf K}\LRf{{\SV}(U) }$.  Therefore, 
the sequence 
\begin{equation}\label{DelinemArla}
\diagram 
{\BK}\LRf{C{\BM}^Z(X)} \ar[r] &{\BK}\LRf{{\SV}(X)} \ar[r]^{\q} & {\BK}\LRf{{\SV}(U)} \\
\enddiagram
\end{equation}
is not a homotopy fibration. However, Thomason computes \cite{TT90} the homotopy fiber of $\q$, as ${\BK}\LRf{{\bf P}erf_Z(X)}$, which is non trivial because $\q$ is not an equivalence. It would be hard to imagine a subcategory ${\SH}\subseteq Coh(X)$, such that ${\BK}\LRf{{\bf P}erf_Z(X)}\cong {\BK}\LRf{{\SH}}$. Note, while $C{\BM}^Z(X)=0$, there is no natural functor
 $Coh^Z(X)\lra {\bf P}erf_Z(X)$. 
}
\end{example}

Differences between our results and that in \cite{TT90} is clarified in (\ref{DelineExample}). While the spectrum ${\BK}\LRf{{\bfP}erf_Z(X)}$ is a description of homotopy fiber of $\q$, it does not answer the question \ref{opnCMZ}. To respond to Problem \ref{opnCMZ}, we need a subcategory ${\SH}\subseteq Coh(X)$ and an equivalence ${\BK}\LRf{{\bf P}erf_Z(X)}\cong {\BK}\LRf{{\SH}}$, which is  essentially impossible in the case of 
(\ref{DelineExample}). In \cite{TT90, S17} there is no natural candidate to take the place of ${\SH}$. 
When $Z$ is complete intersection, we establish that ${\SH}=C{\BM}^Z(X)$ will do the job, and we establish in the Agreement theorems (\ref{AgrBkt}, \ref{AgrBktBGW}) that
${\BK}\LRf{C{\BM}^Z(X)} \cong {\BK}\LRf{{\bfP}erf_Z(X)}$ and ${\BG}W\LRf{C{\BM}^Z(X)} \cong {\BG}W\LRf{{\bfP}erf_Z(X)}$, as two homotopy fibers must agree. 
In particular,  in this article, we settle the conjecture of 
Gersten \cite[pp 343]{G74}
and of Weibel \cite[pp 199]{W89}, which is the   case $\codim(Z)=1$ of  (\ref{opnCMZ});  along with the problem \ref{opnCMZ} in all codimensions $\codim(Z)\geq 1$.

As usual, work on $GW$-theory follows the advancements in $K$-theory \cite{W83, TT90}. We prove the $GW$ analogues of the ${\BK}$-theory conjecture. In \cite{S17}, 
${\BG}W$-theory was developed in the framework is dg categories. As far as this article is concerned, all we need to know is, for an exact category ${\SE}$,  the category $Ch^b\LRf{\SE}$ of bounded complexes in ${\SE}$ has a dg category structure, which we denote by ${\bf dg}{\SE}$. If ${\SE}:=\LRf{{\SE}, ^{\vee}, \varphi}$ has a duality, then ${\bf dg}{\SE}$ is a dg category with duality and weak equivalences (the  quasi isomorphisms).  Due to the shift funtor $T$ on ${\bf dg}{\SE}$, we obtain one shifted ${\BG}W^{[n]}\LRf{{\bf dg}{\SE}}\in \BiSp$ Karoubi ${\BG}W$ bispectrum, for each $n\in {\BZ}$. However, ${\BG}W^{[n]}\LRf{{\bf dg}{\SE}}$ has four periodicity. 
Further, ${\SE}$ itself can be considered as a 
dg category, concentrated at degree zero, which has its  Karoubi ${\BG}W$ bispectrum.
In summary, we have 
$$
\left\{\begin{array}{l}
{\BG}W\LRf{{\SE}}:={\BG}W\LRf{{\SE}, ^{\vee}, \varpi}={\BG}W^{[0]}\LRf{{\bf dg}{\SE}}\\
{\BG}W\LRf{{\SE}}^-:={\BG}W\LRf{{\SE}, ^{\vee}, -\varpi}={\BG}W^{[2]}\LRf{{\bf dg}{\SE}}\\
\end{array}\right.
$$
In section \ref{BISP4Exact} we give an overview of all these on Karoubi ${\BG}W$ bispectrum. In this introduction, we only stated our results on ${\BG}W\LRf{-}^{\pm}$ above (\ref{372Firstresulatr}). We refer to the main body of this article for results on 
${\BG}W^{[n]}\LRf{{\bf dg}-}$.

In section \ref{SecAPPP} we  give a number of applications of our main 
results (\ref{372Firstresulatr}), both on ${\BK}$ theory and ${\BG}W$ theory. These include (1) Cohen Macaulay local rings $(A, \m)$, (2) The case, when $Z=\LRs{\m_1, \m_2, \ldots, m_{\ell}}$ is a set closed points, which is  complete intersection. (3) The case when $U$ is regular, (4) Agreement with Thomason's description on ${\BK}$ theory fiber of $\q$, and Agreement with  
Schlichting's description on Karoubi  ${\BG}W$ theory fiber of $\p$, (5) Excision, (6) Mayer-Vietories.

For completeness of the story we pose the following conjecture.
\begin{conjecture}{\rm 
Suppose $X$ is a quasi projective scheme over a noetherian affine scheme $\spec{A}$, and $d=\dim X$. Let $U$ be an open subset and $Z=X-U$. Assume that $Z$ has a subscheme structure, such that $Z\hra X$ is a regular embedding, with $k=\codim Z=grade(Z, X)$. Then the sequences
\begin{equation}\label{QuiConjecture}
\left\{\begin{array}{ll}
\diagram 
{\BK}\left(C{\BM}^Z(X)\right) \ar[r] &
{\BK}\left({\SV}(X)\right) \ar[r]^{\q} & {\BK}\left({\SV}(U)\right)\\
\enddiagram & in ~\Sp\\
\diagram 
{\BG}W\left(C{\BM}^Z(X)\right)^{\pm} \ar[r] &
{\BG}W\left({\SV}(X)\right) \ar[r]_{\p} & {\BG}W\left({\SV}(U)\right)\\
\enddiagram & in~\BiSp,~k\in 2{\BZ}\\
\diagram 
{\BG}W^{[-k]}\left(C{\BM}^Z(X)\right) \ar[r] &
{\BG}W\left({\SV}(X)\right) \ar[r]_{\p} & {\BG}W\left({\SV}(U)\right)\\
\enddiagram & in~\BiSp\\
\end{array}\right.
\end{equation} 
are homotopy fibrations, in the respective categories.
}
\end{conjecture}
Question of localization has always been an open problem among the expert, and its absence was always missed.  The description in \cite{TT90} did not level it up with the formulation in (\ref{GrayQuillen}). 
 It would be hard to imagine that Quillen himself only had bank sheet,
 beyond the case of codimension one, in spite of Theorem 
\ref{GrayQuillen}. Again, in the light of (\ref{GrayQuillen}), it is reasonable to speculate the Quillen's formulation would been the same as (\ref{QuiConjecture}).
For this reason I would refer to this conjecture as the {\bf localization conjecture of Quillen}. 
To the best of our search there is no documented statement on
 the localization problem. %
The results in this article is an evidence in support of this conjecture. Secondly, regular embedding has always been a gold standard in many facets of mathematics. For example, intersection products behave fairly well for regular embeddings $Z\hra X$ \cite[Chap 6]{Fu}.


Regarding the organization of this article, most notations are summarized in section \ref{SecVorAus}. Notations and overview on hermitian theory are in section \ref{BISP4Exact}. Some background on scheme theory (grades, annihilators, perfect modules) is in section 
\ref{SecGrade}. In section \ref{secFoxbyMor}, we construct certain maps from some Koszul complexes, to certain 
 complexes of modules, in order to approximate homologies, by perfect modules. 
 This section \ref{secFoxbyMor} is the heart of the technical machineries,  that are used in this article. 
In section \ref{equivSec} we establish  equivalences of certain derived categories, which are
used to prove our main results. Our main results on ${\BK}$ theory fiber are in section
\ref{SECBKTHEY}, and the main results on Karoubi ${\BG}W$ theory fiber are in section
\ref{SEcBGWSpSec}. Application are dealt with in section \ref{SecAPPP}.

\vspace{3mm}
\noindent{\bf Acknowledgement:}{\it  I would like to thank Daniel Grayson for many helpful communication and in particular, for his formulation of Problem \ref{opnCMZ}. I also thank Amalendu Krishna for helpful discussions and encouragement.  I am also thankful to  Anand Swant for his hospitality and discussions during my visit to TIFR in summer 2024.}

%

\section{Koszul construction of Foxby}\label{secFoxbyMor}
In this section we prove out main technical machinery. 
Throughout this article we will be working with the following setup.
\begin{setup}\label{CIsetUp}{\rm 
As in section \ref{SecVorAus}, let
 $X$  be a quasi projective scheme over a noetherian affine scheme $\spec{A}$. In particular, $S=\bigoplus_{n\geq 0}S_n$ is a graded noetherian ring and $Y=\proj{S}$,
and $X\subseteq Y$ is open (and dense). In general $Z \subseteq X$ will denote a closed subset,
and $U=X-Z$. 

 Our focus will be when $Z=V(\f_1, \f_2, \ldots, \f_k)\cap X$, where  $\f_1, \f_2, \ldots,  \f_k\in S$ are homogenous. 
Further, $\f_{i_1}, \f_{i_2}, \ldots, \f_{i_j}$ induce regular   sequences on $X$,
  $\forall~1\leq i_1<i_2<\cdots<i_j\leq k$. ({\it This means that the corresponding Koszul complexes on $X$ is exact at all degrees higher than zero.}) We will say that $Z$ is {\bf complete intersection} (with $grade(Z, X)=k$ or $\codim(Z)=k$).
}
\end{setup}

\vspace{3mm}
We recall the following Koszul complex construction.
\begin{construction}\label{consKosz}{\rm 
Let $S=\bigoplus_{n\geq 0}S_n$ be a graded noetherian ring and $Y=\proj{S}$.
Let $f_1, f_2, \ldots, f_k\in S$ be homogeneous elements, with $\deg(f_i)=r_i$. 
The Koszul complex $K_{\bullet}(f_1, f_2, \ldots, f_k)$ is constructed, as usual, while we view at 
it as a complex of graded free (shifted) modules, and of graded $S$-module maps of degree zero. More explicitly, consider the graded free $S$-module 
$$
P=\bigoplus_{i=1}^kS(-r_i)e_i\quad{\rm and ~the ~map}\quad\varphi: P\lra S, \quad {\rm defined ~by}\quad  \varphi(e_i)=f_i\in S.
$$
So,  $\varphi \in Hom(P, S)$ is a map of degree zero. Therefore, 
$K_{\bullet}\LRf{\varphi}=K_{\bullet}\LRf{f_1, f_2, \ldots, f_k}$ 
is the Koszul complex of $\varphi$. In particular, the degree $p$-term of $K_{\bullet}(f_1, f_2, \ldots, f_k)$ is given by
$$
\Lambda^pP=\bigoplus_{1\leq i_1<i_2<\cdots<i_p\leq k} S\left(-(r_{i_1}+r_{i_2}+\cdots+r_{i_p})\right)e_{i_1}\wedge e_{i_2}\wedge \cdots \wedge e_{i_p}
$$
Sheafifying, we obtain the Koszul complex 
${\SK}_{\bullet}\LRf{\widetilde{\varphi}}={\SK}_{\bullet}\LRf{f_1, f_2, \ldots, f_k}$ of locally free sheaves, on $Y$. Again, the degree $p$-term is given by 
$$
\Lambda^p\widetilde{P}=\bigoplus_{1\leq i_1<i_2<\cdots<i_p\leq k} {\CO}_Y\left(-(r_{i_1}+r_{i_2}+\cdot+r_{i_p})\right)e_{i_1}\wedge e_{i_2}\wedge \cdots \wedge e_{i_p}
$$
The complex is exact, only at the points $\wp\in Y$, such that $f_1, f_2, \ldots, f_k$ induce $S_{(\wp)}$-regular sequences. 
}
\end{construction}

%
The following is an improvisation   of a construction of Foxby \cite{F}. 
\bL\label{AppFox1}
{\rm 
Let $A$ be a commutative noetherian ring  and 
$$
\diagram 
\cdots \ar[r] & M_r \ar[r] & M_{r-1} \ar[r] & \cdots\ar[r] & M_1 \ar[r]^{\partial_1} & M_{0} \ar[r]^{\partial_0} &M_{-1} \ar[r] &\cdots\\
\enddiagram 
$$
be a complex of $A$-modules. Suppose $f_1, f_2, \ldots, f_r\in Ann\LRf{\bigoplus H_i(M_{\bullet})}$.
As usual, subsequently we write $Z_0:=Z_0\left(M_{\bullet}\right):=\ker\LRf{\partial_0}$,
$B_0:=B_0\left(M_{\bullet}\right):=Image\LRf{\partial_1}$, and likewise.
Suppose $z\in Z_0\LRf{M_{\bullet}}$ and consider the map $\nu_0:A \lra Z_0$ such that 
$\nu_0(1)=z$. 
Then $\forall m\geq r$ there is a  map of  complexes 
$$
\nu_{\bullet}^m:K_{\bullet}\left(f_1^m, f_2^m, \ldots, f_r^m\right) \lra M_{\bullet} \qquad 
\ni \qquad \nu_{0}^m=\nu_0
$$
Further, with diagonal map $\Delta:A^r \lra A^r$, defined by $(f_1, f_2, \ldots, f_r)$, 
the maps $\nu_{\bullet}^m$ can be chosen consistently, as in the commutative diagram %
\begin{equation}\label{consi491nubull}
\diagram
K_{\bullet}(f_1^m, f_2^m, \ldots, f_r^m) \ar[r]^{\qquad \nu_{\bullet}^m}\ar[d]_{\Lambda\Delta} & M_{\bullet}\ar@{=}[d]\\
K_{\bullet}(f_1^{m+1}, f_2^{m+1}, \ldots, f_r^{m+1})\ar[r]_{\qquad \qquad \nu_{\bullet}^{m+1}} &M_{\bullet} \\
\enddiagram 
\end{equation}
 }
\eL
\pf Let $r=1$. Then we have $f\in Ann\left(\bigoplus H_i\left(M_{\bullet}\right)\right)$. 
We note that $fz\in {B}_0(M_{\bullet})$ is a boundary. So, there is a choice $m_1\in M_1$ such that $d_1(m_1)
=fz$. Therefore, 
we have the diagram
$$
\diagram
 &  &   & 0\ar[r] & A \ar[r]^{f^2} \ar@{-->}[d]_{f} & A \ar[r]\ar[d]^1  &0 &\\
 &  &   & 0\ar[r] & A \ar[r]^f \ar@{-->}[d]_{\nu_1} & A \ar[r]\ar[d]^{\nu_0}  &0 &\\
\cdots \ar[r] & M_r \ar[r] & M_{r-1} \ar[r] & \cdots\ar[r] & M_1  \ar[r]_{d_1} & M_{0} \ar[r] &M_{-1} \ar[r] &\cdots\\
\enddiagram
$$
This settles the proof when $r=1$. 
Now suppose the proof is done for any sequence of length $r-1$. Let 
$$
\left\{\begin{array}{l}
F=\oplus_{j=1}^r Ae_j\quad {\rm denote~the~free~module}\\
\hat{F}_i=\oplus_{j=1, i\neq j}^r Ae_i\\
e_{i_1i_2\cdots i_p}=e_{i_1}\wedge e_{i_2} \wedge \cdots \wedge e_p\quad \forall 1\leq i_1<i_2<\cdots <i_p\leq r\\
\hat{e}_i=e_{1}\wedge \cdots e_{i-1} \wedge e_{i+1}\wedge \cdots\wedge e_{r-1} \wedge e_r\qquad {\rm skipping}\qquad e_i\\
\end{array}\right. 
$$
 By induction we assume that  for each basis element $e_{i_1i_2\ldots i_p}$, with $p\leq r-1$, we made one choice 
$m_{i_1i_2\ldots i_p}\in M_p$, so that the following maps $\nu^p_p$ are defined, sending 
$e_{i_1i_2\ldots i_p}\mapsto m_{i_1i_2\ldots i_p}$. 
In particular, for all $i$, we made choices $\hat{m}_i\in M_{r-1}$ such that, we have maps
$$
\diagram
 & 0\ar[r] & A\hat{e_i}\ar[d]_{\Lambda^{r-1}{\Delta_i}} \ar[r]&\Lambda^{r-2}\hat{F}_i \ar[d] \ar[r]& \cdots\ar[r] & \hat{F}_i \ar[r]^{\hat{f}_i^{r}} \ar@{-->}[d]_{\Delta_i} & A \ar[r]\ar[d]^{1}  &0 &\\
 & 0\ar[r] & A\hat{e}_i\ar[d]_{\hat{e}_i\mapsto \hat{m}_i} \ar[r]&\Lambda^{r-2}\hat{F}_i \ar[d] \ar[r]& \cdots\ar[r] & \hat{F}_i \ar[r]^{\hat{f}_i^{r-1}} \ar@{-->}[d]_{\nu_1} & A \ar[r]\ar[d]^{\nu_0}  &0 &\\
\cdots \ar[r] & M_r \ar[r] & M_{r-1} \ar[r] & M_{r-2}\ar[r] & \cdots  \ar[r] & M_{1} \ar[r] &M_{0} \ar[r] &\cdots\\
\enddiagram
$$
Clearly, these vertical maps agree on the canonical bases of $\Lambda^pF$ for all $p\leq r-1$. Combining these, we get the maps for $f_1^{r-1}, f_2^{r-1}, \ldots, f_r^{r-1}$. We obtain 
$$
\diagram
0\ar[r] & \Lambda^rF\ar[r]^{\kappa_r^{r}}\ar[d]_{\det(\Delta)}\ar@/_/@{-->}[dddl] & \Lambda^{r-1}F\ar[d]   \ar[r]&\Lambda^{r-2}F \ar[d] \ar[r]& \cdots\ar[r] & F \ar[r]^{\hat{f}_i^{r}} \ar@{-->}[d]_{\Delta} & A \ar[r]\ar[d]^{1}  &0 &\\
 0\ar[r] & \Lambda^rF\ar[r]^{\kappa_r^{r-1}}\ar@{-->}[ddr] &
  \Lambda^{r-1}F\ar[d]_{\hat{e_i}\mapsto \hat{m}_i}^{\nu_{r-1}^{r-1}} \ar[r]&\Lambda^{r-2}F \ar[d] \ar[r]& \cdots\ar[r] & F \ar[r]^{\hat{f}_i^{r-1}} \ar@{-->}[d]_{\nu_1} & A \ar[r]\ar[d]^{\nu_0}  &0 &\\
\cdots \ar[r] & M_r \ar[r]_{d_r} \ar@{-->>}[dl]_{d_r} & M_{r-1} \ar[r]_{d_{r-1}} & M_{r-2}\ar[r] & \cdots  \ar[r] & M_{1} \ar[r] &M_{0} \ar[r] &\cdots\\
 B_{r-1} \ar@{^(-->}[rr] && Z_{r-1} \ar@{^(-->}[u]  &  &   &   &  & \\
\enddiagram
$$
Write $e=e_1\wedge e_2\wedge \cdots \wedge e_r$. We have
$$
\nu_{r-1}^{r-1}\kappa_r^{r-1}(e)= \sum \pm f_i^{r-1}\nu_{r-1}^{r-1}(\hat{e_i}) =0
\quad
{\rm So, }\quad \nu_{r-1}^{r-1}\kappa_r^{r-1}(e)\in Z_{r-1}\left(M_{\bullet}\right)\quad {\rm is~a~cycle.}
$$
Therefore
$$
\nu_r^r\kappa_r^r(e)=\det(\Delta)\left(\nu_{r-1}^{r-1}\kappa(e)\right)\in B_{r-1}\left(M_{\bullet}\right)\quad {\rm is~a~boundary.}
$$
So, we can pick 
$$
m_{12\cdots r}\in M_r \qquad \ni \qquad d_r\left(m_{12\cdots r}\right)= \nu_r^r\kappa_r^r(e)
$$
Define 
$$
\nu_r^r(e)=m_{12\cdots r}
$$
This settles the inductive step, for $m=r$. For $m\geq r+1$, $\nu_{\bullet}^m$ is defined
 by the consistency diagram 
(\ref{consi491nubull}). 

This completes the inductive case.
\pic $\eop$ 

\vspace{3mm} 
Now we do the same (\ref{AppFox1}) for  projective schemes $Y$.

 \bP\label{projVerFox}{\rm 
 Let $Y=\proj{S}$ be projective scheme over a noetherian affine scheme $\spec{A}$. 
 Let 
 \begin{equation}\label{YmeKaroCx}
 \diagram 
\cdots \ar[r]& {\SF}_k\ar[r]&  {\SF}_{k-1} \ar[r]& \cdots \ar[r]  & {\SF}_{1} \ar[r]& {\SF}_{0} \ar[r]^{\partial_0}&{\SF}_{-1}\\
\enddiagram
\end{equation} 
be a bounded complex of coherent ${\CO}_Y$-module. 
Write $M=\bigoplus {\CH}_i\LRf{{\SF}_{\bullet}}$ and $I\subseteq S$ be homogeneous ideal such that $\tilde{I}=ann(M)$. Let $\f_1, \f_2, \ldots, \f_k\in I$ be homogenous elements. Then there is a map of complexes 
$\nu_{\bullet}:
{\SE}_{\bullet}\lra {\SF}_{\bullet}$ such that 
\bE
\item\label{OneFoxE0} First,  ${\SE}_i\in {\SV}(Y)~\forall ~i$, and ${\SE}_i=0$ unless $0\leq i\leq k$. 
\item\label{TWoFoxBabu} The map $\nu_0:{\SE}_0\sur {\CZ}_0\LRf{{\SG}_{\bullet}}$ is surjective,  and hence
${\CH}_0(\nu_{\bullet}): {\CH}_0\LRf{{\SE}_{\bullet}} \sur {\CH}_0\LRf{{\SG}_{\bullet}}$ is surjective. 
\item\label{CH0Kcx} In fact, ${\SE}_{\bullet}$ will be the direct sum of certain Koszul complexes, similar to the map 
(\ref{consi491nubull}). And 
\begin{equation}\label{H0EbullIsDirsom}
{\CH}_0\LRf{{\SE}_{\bullet}}\otimes {\CL}= \LRf{\frac{{\CO}_Y}{\widetilde{\LRf{\f_1^m, \f_2^m, \ldots, \f_k^m}}}}^t\quad m\gg 0, ~for~some~t
\end{equation} 
and for some line bundle ${\CL}$ (i.e. {\it up to a twist}). 
\item\label{CMZFoxReg} 
Assume $\f_1, \f_2, \ldots, \f_k$  induce an regular sequence on an open subset $X\subseteq Y$. 
Then 
$$
\LRf{{\CH}_0\LRf{{\SE}_{\bullet}}}_{|X}\in C{\BM}^Z(X)\qquad with\quad Z=V\LRf{\f_1, \f_2, \ldots, \f_k}\cap X.
$$
\eE 
 }
 \eP
 \pf As the diagram (\ref{H0EbullIsDirsom}) suggets, for fixed $k$ and $m\geq k$ we will construct one map $\nu_{\bullet}^m$ and $\nu_{\bullet}:=\nu_{\bullet}^m$ will be a map, as requited. 
 By twisting (\ref{YmeKaroCx}) high enough, we assume ${\SF}_i$, the boundaries
 ${\CB}_i\left({\SF}_{\bullet}\right)$, the cycles ${\CZ}_i\left({\SF}_{\bullet}\right)$, and the homologies ${\CH}_i\LRf{{\SF}_{\bullet}}$ are globally generated and 
 $$
 \left\{\begin{array}{ll}
 {H}^i\left(Y, {\SF}_{i}\right) =0 &\forall i\geq 1\\
  {H}^i\left(Y, {\CB}_j\left({\SF}_{\bullet}\right)\right) =0 &\forall i\geq 1, \forall j\\
    {H}^i\left(Y, {\CZ}_j\left({\SF}_{\bullet}\right)\right) =0 &\forall i\geq 1, \forall j\\
        {H}^i\left(Y, {\CH}_j\left({\SF}_{\bullet}\right)\right) =0 &\forall i\geq 1, \forall j\\
 \end{array}\right.
 \qquad \qquad {\rm \cite[III.5.2, III.5.3]{H}}
 $$
Let $z\in H^0\left(Y, {\CZ}_0\left({\SF}_{\bullet}\right)\right)$ and let 
$\mu_z:{\CO}_Y \lra {\CZ}_0\left({\SF}_{\bullet}\right)\subseteq {\SF}_0$ be the map induced by $1\mapsto z$.
By induction on $k$, we prove that for any selection of $1\leq i_1<i_2< \cdots < i_r\leq k$, there is a 
 map of complexes
\begin{equation}\label{787IndHyp}
\nu_{\bullet}^r:{\CK}_{\bullet}\left(\f_{i_1}^m, \f_{i_2}^m, \ldots, \f_{i_r}^m \right)\lra {\SF}_{\bullet}
\qquad \forall m\geq r
\end{equation}
\bE
\item This Koszul complex 
${\CK}_{\bullet}\left(\f_{i_1}^m, \f_{i_2}^m, \ldots, \f_{i_r}^m \right)={\CK}_{\bullet}\LRf{\widetilde{\varphi^r_{i_1\ldots i_r}}}$ is built 
from the map of graded modules  $\varphi^r_{i_1\ldots i_r}:F_{i_1i_2\ldots i_r}=\bigoplus_{j=1}^r Se_{i_j}\lra S$ sending 
$e_i\mapsto \f_i^m$. 
\item $\nu_0^r=\mu_z$
\item With the map
$$
\Delta:=\Delta_{i_1i_2\ldots t_r}=Diagonal\left(\f_{i_1}, \ldots, \f_{i_r}\right):F_{i_1i_2\ldots i_r} \lra F_{i_1i_2\ldots i_r}
$$ 
we define $\nu_{\bullet}^{m+1}= \nu_{\bullet}^{m}\widetilde{\Delta}$.
So, there are commutative diagrams of maps in $Ch^b\left(Coh(Y)\right)$: 
$$
\diagram
{\CK}_{\bullet}\left(\f_{i_1}^{m+1}, \f_{i_2}^{m+1}, \ldots, \f_{i_r}^{m+1} \right)\ar[r]^{\quad \widetilde{\Delta}}\ar@/_/@{-->}[dr]_{\nu_{\bullet}^{m+1}}&
{\CK}_{\bullet}\left(\f_{i_1}^m, \f_{i_2}^m \ldots, \f_{i_r}^m \right)\ar[d]^{\nu_{\bullet}^m} \\
& {\SF}_{\bullet}\\
\enddiagram
\qquad \forall m\geq r
$$

\eE 
Denote $d_i=\deg(\f_i)$. 
Consider the diagram 
\begin{equation}\label{fzisZETR}
\diagram
S(-d_1) \ar[r]^{\f_1} & S\ar[d]^{\mu_z}\\ 
\Gamma_*\left( Y, {\CB}_0\left(Y, {\SF}_{\bullet}\right)\right) \ar@{^(->}[r]& \Gamma_*\left( Y, {\CZ}_0\left(Y, {\SF}_{\bullet}\right)\right)\\
\enddiagram
\end{equation}
The sheafification, leads to diagram 
$$
\diagram
{\CO}_Y(-d_1)\ar[r]^{\f_1} & {\CO}_Y\ar[d]^{\mu_z} & \\ 
\widetilde{\Gamma_*\left( Y, {\CB}_0\left(Y, {\SF}_{\bullet}\right)\right)}\ar@{^(->}[r]
&\widetilde{\Gamma_*\left( Y, {\CZ}_0\left(Y, {\SF}_{\bullet}\right)\right)}\ar[r]_{\sim}& \iota_*{\CZ}_0\left(Y, {\SF}_{\bullet}\right)\\
\enddiagram
$$
Note, the map $\mu_z\f_1:{\CO}_Y\lra {\CZ}_0\left({\SF}_{\bullet}\right)(d_1)$ is given by $1\mapsto \f_1z$. We have the exact sequence 
$$
\diagram 
0\ar[r] & H^0\left(Y, {\CB}_0\left({\SF}_{\bullet}\right)(d_1) \right) \ar[r] & 
H^0\left(Y, {\CZ}_0\left({\SF}_{\bullet}\right)(d_1)\right) \ar[r] & 
H^0\left(Y,{\CH}_0\left({\SF}_{\bullet}\right) (d_1)\right) \ar[r] & 0\\
\enddiagram 
$$
Since $\f_1\in Ann\left(({\CH}_0\left({\SF}_{\bullet}\right)\right)$, the global
section $\f_1z\in 
H^0\left(Y, {\CB}_0(d_1)\right)$. It follows from the exact sequence above, that the above map factors as follows:
$$
\diagram
{\CO}_Y(-d_1)\ar[r]^{\f_1}\ar@{-->}[d]_{\mu_z\f_1} & {\CO}_Y\ar[d]^{\mu_z}& \\ 
 {\CB}_0\left({\SF}_{\bullet}\right)\ar@{^(->}[r]
&  {\CZ}_0\left({\SF}_{\bullet}\right)\ar@{->>}[r]&{\CH}_0\left({\SF}_{\bullet} \right)\\
\enddiagram
$$
The vertical broken arrow is given by the global section $\f_1z\in H^0\left(Y, {\CB}_0\left({\SF}_{\bullet}\right)(d_1) \right)$. 
We also have the exact sequence
$$
\diagram
0\ar[r] & H^0\left(Y, {\CZ}_1\left({\SF}_{\bullet}\right)(d_1) \right) \ar[r] & 
H^0\left(Y, {\SF}_{1}(d_1)\right) \ar[r]^{\varepsilon\qquad} & 
H^0\left(Y,{\CB}_0\left({\SF}_{\bullet}\right) \right)(d_1) \ar[r] & 0\\
\enddiagram
$$
Now we make a choice $m_1\in H^0\left(X, {\SF}_{1}(d_1)\right)$ such that 
$\varepsilon(m_1)=\f_1z$.  %
This gives a 
map
$$
\diagram 
{\CO}_Y \ar[rr]^{\nu_1}_{1\mapsto m_1}& &{\SF}_1(d_1)\\
\enddiagram 
$$
So, the above diagram lifts as follows:
\begin{equation}\label{6Oct24Modi}
\diagram
&{\CO}_X(-d_1)\ar[r]^{\f_1}\ar@{-->}[d]_{\mu_z\f_1}\ar@/_/@{-->}[dl]_{m_1(d_1)} & {\CO}_X\ar[d]^{\mu_z}& \\ 
{\SF}_1\ar[r]& {\CB}_0\left({\SF}_{\bullet}\right)\ar@{^(->}[r]
&  {\CZ}_0\left({\SF}_{\bullet}\right)\ar@{->>}[r]&{\CH}_0\left({\SF}_{\bullet} \right)\\
\enddiagram
\end{equation} 
Likewise we make one choice $m_i$ for $i=1, 2, \ldots, k$, and do the same.
This completes the proof, when $r=1$. 

Now suppose it is done for any sequence of length $r-1\leq k-1$. In fact, 
we can assume $r=k$, for notational conveniences. 
Before we proceed with the inductive step, we fix some notation
$$
\left\{\begin{array}{l}
F=\bigoplus_{j=1}^k Se_{j}\\
F_{i_1i_2\ldots i_r}= \bigoplus_{j=1}^k Se_{i_j}\\
\hat{F}_i=\bigoplus_{j=1, j\neq i}^k Se_{j}\\
e_{i_1i_2\cdots i_p}=e_{i_1}\wedge e_{i_2} \wedge \cdots \wedge e_{i_p}
\quad \forall 1\leq i_1<i_2<\cdots <i_p\leq k\\
\hat{e}_j=e_{1}\wedge \cdots \wedge e_{i-1} \wedge e_{i+1}\wedge \cdots \wedge e_{k}\\
\e=e_{1}\wedge \cdots \wedge e_{i}\wedge \cdots \wedge e_{k}\\
\end{array}\right. 
$$
We define two sets of maps
$$
\left\{\begin{array}{ll}
\varphi^m:=\varphi^m_{i_1i_2\ldots i_p}: F_{i_1i_2\ldots i_p}\lra S & e_{i_j}\mapsto \f_{i_j}^m\\
\Delta:=\Delta_{i_1i_2\ldots i_p}: F_{i_1i_2\ldots i_p}\lra F_{i_1i_2\ldots i_p} & 
e_{i_j}\mapsto \f_{i_j}e_{i_j}\\
\end{array}\right.
$$
%
%
By inductive hypothesis, for all $1\leq i_1<i_2<\ldots <i_p\leq k$, with $p\leq k-1$ we made choices 
\begin{equation}\label{choice911}
m_{i_1i_2\ldots i_p}\in H^0\left(Y, {\SF}_{r-1}\left(-(d_{i_1}+\cdots+d_{i_p}) \right)\right), \quad 
{\rm Denote}\quad \hat{m}_i=m_{1,\ldots, i-1, i+1, \ldots, k}
\end{equation}
For $j=1, 2, \ldots, k$, such choices lead to  the maps, 
$$
\diagram
 & 0\ar[r] & S\hat{e}_j\ar[d]_{\Lambda^{r-1}{\Delta}} \ar[r]&\Lambda^{r-2}\hat{F}_j \ar[d] \ar[r]& \cdots\ar[r] & \hat{F}_j \ar[r]^{\varphi^{k}} \ar@{-->}[d]_{\Delta} & S \ar[r]\ar[d]^{1}  &0 &\\
 & 0\ar[r] & S\hat{e}_j\ar[d]_{\hat{e}_j\mapsto \hat{m}_j} \ar[r]&\Lambda^{r-2}\hat{F}_j \ar[d] \ar[r]& \cdots\ar[r] & \hat{F}_j \ar[r]^{\varphi^{k-1}} \ar@{-->}[d]_{\nu_1} & S \ar[r]\ar[d]^{\nu_0}  &0 &\\
\cdots \ar[r] & M_r \ar[r] & M_{r-1} \ar[r] & M_{r-2}\ar[r] & \cdots  \ar[r] & M_{0} \ar[r] &M_{-1} \ar[r] &\cdots\\
\enddiagram
$$
where,   
we denote
$$
M_i=\Gamma_*\left(Y, {\SF}_i \right)=\bigoplus_{n\geq 0}H^0\left(Y, {\SF}_i(n)\right)
$$
 These maps agree on the common bases of the respective Koszul complexes. Combining these $k$ maps, we obtain the following map
$$
\diagram
0\ar[r] & S{\e}\ar[r]^{\kappa^k}\ar[d]_{\det(\Delta)} & {\Lambda^{k-1}F}
\ar[d]^{\Lambda^{k-1}{\Delta}} \ar[r]&\Lambda^{k-2}F \ar[d] \ar[r]& \cdots\ar[r] & F \ar[r]^{\varphi^{k}} \ar@{-->}[d]_{\Delta} & S \ar[r]\ar[d]^{1}  &0 &\\
0\ar[r] & S{\e}\ar[r]^{\kappa^{k-1}} &  \Lambda^{k-1}F
\ar[d]_{\hat{e_j}\mapsto \hat{m}_j}^{\nu_{k-1}^{k-1}} \ar[r]&\Lambda^{k-2}F\ar[d] \ar[r]& \cdots\ar[r] & F \ar[r]^{\varphi^{k-1}} \ar@{-->}[d]_{\nu_1} & S \ar[r]\ar[d]^{\nu_0}  &0 &\\
\cdots \ar[r] & M_k \ar[r] & M_{k-1} \ar[r]_{d_{k-1}} & M_{k-2}\ar[r] & \cdots  \ar[r] & M_{1} \ar[r] &M_{0} \ar[r] &\cdots\\
\enddiagram
$$
It follows $d_{k-1}\nu_{k-1}^{k-1}\kappa^{k-1}=0$. Therefore 
$\nu_{k-1}^{k-1}\kappa^{k-1}(\e)\in \ker(d_{k-1})$. Now sheafify this diagram. 
We obtain the diagram
{\scalefont{.7}
$$
\diagram
0\ar[r] & {\CO}_Y(-kd)\ar[r]^{\kappa^k}\ar[d]_{\det(\Delta)} & {\CK}_{k-1}(\varphi^k)\ar[d]^{\Lambda^{k-1}{\Delta}} \ar[r]&
{\CK}_{k-2}(\varphi^k) \ar[d] \ar[r]& \cdots\ar[r] &{\CK}_{1}(\varphi^k)  \ar[r]^{\varphi^{k}} \ar@{-->}[d]_{\Delta} & {\CO}_Y \ar[r]\ar[d]^{1}  &0 &\\
0\ar[r] & {\CO}_Y(-(k-1)d)\ar[r]^{\qquad\kappa^{k-1}} &  {\CK}_{k-1}(\varphi^{k-1})
\ar[d]_{\hat{e_j}\mapsto \hat{m}_j}^{\nu_{k-1}^{k-1}} \ar[r]&{\CK}_{k-2}(\varphi^{k-1}) \ar[d] \ar[r]& \cdots\ar[r] & {\CK}_{1}(\varphi^{k-1})  \ar[r]^{\varphi^{k-1}} \ar@{-->}[d]_{\nu_1} &
 {\CO}_Y \ar[r]\ar[d]^{\nu_0}  &0 &\\
\cdots \ar[r] & {\SF}_k \ar[r] & {\SF}_{k-1} \ar[r]_{d_{k-1}} & {\SF}_{k-2}\ar[r] & \cdots  \ar[r] & {\SF}_{1} \ar[r] &{\SF}_{0} \ar[r] &\cdots\\
\enddiagram
$$
}
where $d=\sum_{i=1}^k d_i$. We  redraw this diagram:
\begin{equation}\label{946geneousdia}
\diagram
0\ar[r] & {\CO}_Y(-kd)\ar[r]^{\kappa^k}\ar[d]_{\det(\Delta)}\ar@/_/@{-->}[dddl]_{\zeta} & {\CK}_{k-1}(\varphi^k) \ar[d]^{\Lambda^{k-1}{\Delta}} \ar[r]&{\CK}_{k-2}(\varphi^k)  \ar[d] \\
0\ar[r] & {\CO}_Y(-(k-1)d)\ar[r]^{\qquad\kappa^{k-1}}\ar@{-->}[rdd]^{\eta} &  {\CK}_{k-1}(\varphi^{k-1})
\ar[d]_{\hat{e_j}\mapsto \hat{m}_j}^{\nu_{k-1}^{k-1}} \ar[r]&{\CK}_{k-1}(\varphi^{k-1})
 \ar[d] \\ 
\cdots \ar[r] & {\SF}_k \ar[r]\ar@{->>}[dl] & {\SF}_{k-1} \ar[r]_{d_{k-1}} & {\SF}_{k-2}\\ 
 {\CB}_{k-1}\left({\SF}_{\bullet}\right) \ar@{^(->}[rr] && {\CZ}_{k-1}\left({\SF}_{\bullet}\right)  \ar@{^(->}[u]  & \\ 
\enddiagram
\end{equation}
It follows $\nu_{k-1}^{k-1}\kappa^{k-1}$ factors through $\eta$. The exact sequence 
$$
\diagram 
0\ar[r] &{\CB}_{r-1}\left({\SF}_{\bullet}\right) \ar[r]& {\CZ}_{r-1}\left({\SF}_{\bullet}\right) \ar[r] &{\CH}_{r-1}\left({\SF}_{\bullet}\right) \ar[r] & 0\\
\enddiagram
$$
induces the exact sequence 
$$
\diagram 
0\ar[r] &\Gamma_*\left(X, {\CB}_{r-1}\right)\ar[r] & \Gamma_*\left(X, {\CZ}_{r-1}\right) \ar[r] &
\Gamma_*\left(X, {\CH}_{r-1}\right) \ar[r] & 0\\
\enddiagram
$$
Note that the map $\det(\Delta)$ is induced by the product $\f_1\f_2\cdots \f_k$. 
Since $\f_i\in I$ 
it follows that the $\det(\Delta)\eta$ factors thorough a map $\zeta:{\CO}_Y(-kd)\lra {\CB}_{k-1}\left({\SF}_{\bullet}\right)$, as in diagram (\ref{946geneousdia}). 
We further have the exact sequence
$$
\diagram 
0\ar[r] & {\CZ}_r\left({\SF}_{\bullet}\right) \ar[r] & {\SF}_r\left({\SF}_{\bullet}\right) \ar[r] & {\CB}_{r-1}\left({\SF}_{\bullet}\right) \ar[r] & 0\\
\enddiagram 
$$
This induces the exact sequence
$$
\diagram 
0\ar[r] &\Gamma_*\left(Y, {\CZ}_{r}\left({\SF}_{\bullet}\right)\right)\ar[r]^{\beta} & \Gamma_*\left(Y, {\SF}_{r}\right) \ar[r]^{\q} &
\Gamma_*\left(Y, {\CB}_{r-1}\left({\SF}_{\bullet}\right)\right) \ar[r] & 0\\
\enddiagram
$$
The map $\zeta$ is induced by a map $\omega:S{\e} \lra H^0\left(Y, {\SB}_{k-1}\left({\SF}_{\bullet}\right)(kd)\right)$. Write
$\omega(\e)=\varepsilon\in H^0\left(Y, {\SB}_{k-1}(kd)\right)$. Now we make a choice 
$$
m_{123\ldots k} \in H^0\left(Y, {\SF}_{k}(kd)\right)\qquad \ni \qquad 
\q(m_{123\ldots k}) =\varepsilon 
$$
Define 
$$
\nu^k_k:{\CO}_Y(-kd) \lra {\SF}_k\qquad {\rm by}   \qquad \nu_k^k(\e) =m_{123\ldots k}
$$
This completes the construction of (\ref{787IndHyp}), as hypothesized.

Note, we assumed ${\CZ}_0\left({\SF}_{\bullet}\right)$ is generated by finitely many global sections 
$z_1, z_2, \ldots, z_t\in H^0\left(Y, {\CZ}_0\left({\SF}_{\bullet}\right)\right)$. 
We construct and define $t$ maps :
$$
\nu_{\bullet}^{k, j}:{\CK}\left(\f_1, \ldots, \f_k\right)\lra {\SF}_{\bullet}\qquad \qquad \nu_{0}^{k, j}(1)=z_j
$$
Taking direct sum 
$$
\nu_{\bullet}=\bigoplus_{j=1}^t\nu_k^{k,j}: {\SE}_{\bullet}:=\bigoplus_{j=1}^t{\CK}\left(\f_1, \ldots, \f_k\right)\lra 
{\SF}_{\bullet}
$$
With $\nu_{\bullet}=\nu_{\bullet}^n$ with $n\geq k$ the proposition is established. To see this, first we have
\begin{equation}\label{1214Didd}
{\CH}_0\left({\SE}_{\bullet}\right)= \left(\frac{{\CO}_X}{{\widetilde{(\f_1^n, \ldots, \f_k^n)}}}\right)^t
\quad {\rm in~this ~case, ~after ~some~twist.} 
\end{equation}
So, (\ref{CH0Kcx}) holds up to some twist. Likewise, (\ref{OneFoxE0}, \ref{TWoFoxBabu}) in this twisted case, and hence they hold after twisting back, as well.

Finally, since $\f_1, \ldots, \f_k$ induce a regular sequence on $X$, it follows from 
(\ref{1214Didd}) that ${\CH}_0\LRf{{\SE}_{\bullet}}_{|X}\in C{\BM}^Z(X)$. So, (\ref{CH0Kcx}) is also established. 
\pic $\eop$
 \vspace{3mm}


  Before we do the quasi projective version of (\ref{projVerFox}),  we insert the following extension of \cite[Ex. II.5.15]{H}.  
\bL\label{HartExtII5p15}{\rm 
Let $Y$ be a noetherian scheme and $\iota: X \hra Y$ be an open subset. Let ${\SG}_{\bullet}$ be a complex coherent ${\CO}_X$-modules, as follows: 
$$
\diagram
{\SG}_{k+1} \ar[r] & {\SG}_k \ar[r]^{\partial_k} \ar[r] & \cdots \ar[r] &
{\SG}_{r} \ar[r]^{\partial_r} & {\SG}_{r-1} \ar[r]
&  \cdots \ar[r] &
{\SG}_{0} \ar[r] & {\SG}_{-1} \\
\enddiagram
$$
Consider the direct image complex $\iota_*{\SG}_{\bullet}$, which is a complex of quasi coherent ${\CO}_Y$-modules. Then there is a sub complex ${\SF}_{\bullet} \subseteq 
\iota_*{\SG}_{\bullet}$ of coherent ${\CO}_Y$-modules, which extends ${\SG}_{\bullet}$.

}
\eL
\pf By \cite[Ex II.1.15]{H}, for $r=-1, 0, \ldots, k, k+1$,  there are coherent ${\CO}_Y$-modules ${F}_r\subseteq \iota_*{\SG}_r$ which extends ${\SG}_r$. However the direct image $\iota_*\partial_r$ may not restrict to a map from ${F}_r$ to $F_{r-1}$. To correct this, inductively we define 
 $$
 {\SF}_{k+1}=F_{k+1}, {\SF}_{k}=F_{k}+\iota_*\partial_{k+1}\LRf{{\SF}_{k+1}},
 \ldots, {\SF}_r=F_{k}+\iota_*\partial_{r+1}\LRf{{\SF}_{r+1}}, \cdots 
 $$
 For open subsets $U\subseteq X$, inductively, we have
 $$
 \left\{\begin{array}{l}
 {\SF}_{k+1}(U)= F_{k+1}(U)={\SG}_{k+1}(U)\\
  {\SF}_{k}(U)= F_{k}(U)+\partial_{k+1}{\SG}_{k+1}(U) ={\SG}_{k}(U)\\
  \cdots\\
    {\SF}_{r}(U)= F_{r}(U)+\partial_{r+1}{\SG}_{k+1}(U) ={\SG}_{r}(U)\\
      \cdots\\
 \end{array}\right.
 $$
 \pic $\eop$ 

\vspace{3mm}
The following is the main result in this section, which is a quasi projective  version of (\ref{projVerFox}).
\bT\label{Them23Foxby}{\rm 
Let $X$ be a quasi projective scheme over  a noetherian affine scheme $\spec{A}$ and let $Z\subseteq X$ be a closed subset, with $grade(Z, X)=k$. 
 Let
\begin{equation}\label{1719Prime}
\diagram 
\cdots \ar[r]& {\SG}_k\ar[r]&  {\SG}_{k-1} \ar[r]& \cdots \ar[r]  & {\SG}_{1} \ar[r]& {\SG}_{0} \ar[r]^{\partial_0}&
{\SG}_{-1}\\
\enddiagram 
\end{equation}
be a complex in $Coh(X)$, such that the homologies ${\CH}_i\left({\SG}_{\bullet}\right)\in Coh^Z(X)~\forall i=0, 1, \ldots, k$. Then there is a sequence  $\f_1, f_2, \ldots, \f_k$ on $X$ such that $Z\subseteq V(\f_1, \f_2, \ldots, \f_k)$ and a map of complexes 
$\nu_{\bullet}: {\SE}_{\bullet}\lra {\SG}_{\bullet}$ such that 
\bE
\item\label{OneFoxQuasiE0} First,  ${\SE}_i\in {\SV}(X)~\forall i$, and ${\SE}_i=0$ unless $0\leq i\leq k$. 
\item\label{TWoFoxQuasi} The map $\nu_0:{\SE}_0\sur {\CZ}_0\LRf{{\SG}_{\bullet}}$ is surjective,  and hence
${\CH}_0(\nu_{\bullet}): {\CH}_0\LRf{{\SE}_{\bullet}} \sur {\CH}_0\LRf{{\SG}_{\bullet}}$ is also surjective. 
\item\label{CH0QuasiKcx} In fact, ${\SE}_{\bullet}$ will be the direct sum of certain Koszul complexes, similar to the map 
(\ref{consi491nubull}). And 
\begin{equation}\label{H0EbullIsDirsom}
{\CH}_0\LRf{{\SE}_{\bullet}}\otimes {\CL}= \LRf{\frac{{\CO}_X}{\widetilde{\LRf{\f_1^m, \f_2^m, \ldots, \f_k^m}_{|X}}}}^t\quad m\gg 0, ~for~some~t
\end{equation} 
and for some line bundle ${\CL}$ (i.e. {\it up to a twist}). 
\item\label{FinInCMZX}  Finally, if $Z=V\LRf{f_1, f_2, \ldots, f_k}\cap X$ complete intersection as in Setup \ref{CIsetUp}, then further 
$
\LRf{{\CH}_0\LRf{{\SE}_{\bullet}}}_{|X}\in C{\BM}^Z(X)
$.
\eE 


}
\eT 
\pf 
We assume ${\SG}_n =0, ~ n\geq k+2$, because there is no loss of generality.
We will complete the prrof, by an application of (\ref{projVerFox}). 
We will deal with (\ref{FinInCMZX}) at the end and prove the rest first. Assume $X\subseteq Y=\proj{S}$ is open and dense, 
where $S=\bigoplus_{n\geq 0}S_n=A[x_1, x_2, \ldots, x_N]$, with $\deg{x_i}=1$.
%
%
 %
 %
 Let $Z=\bigcup_{j=1}^{\ell}V(\wp_j)\cap X$ be an irreducible decomposition,
 with  $\wp_j\in X$ and $J=\bigcap_{j=1}^{\ell}{\wp}_j$. 
 
 By (\ref{HartExtII5p15}) there is a subcomplex 
 ${\SF}_{\bullet} \subseteq \iota_*{\SG}_{\bullet}$ of coherent ${\CO}_Y$-modules such that the restriction 
 ${{\SF}_{\bullet}}_{|X}= {\SG}_{\bullet}$. Since  ${\iota_*{\SG}_{\bullet}}_{|Y -V(J)}=0$, we have ${{\SF}_{\bullet}}_{|Y -V(J)}=0$. It follows that %
 $\supp{\bigoplus {\CH}_i\LRf{{\SF}_{\bullet}}}\subseteq V(J)$.  Let $I\subseteq S$ be a homogeneous ideal such that $\tilde{I}=ann\LRf{\bigoplus {\CH}_i\LRf{{\SF}_{\bullet}}}$. 
 So, $V(I)\subseteq V(J)$ and hence $J^p\subseteq I$, for $p\gg 0$. Further
  ${\SI}:=
 \LRf{\bigoplus {\CH}_i\LRf{{\SG}_{\bullet}}}= \tilde{I}_{|X}$. 
 
 Inductively, we pick $\f_1, \f_2, \ldots, \f_k\in J^p$ that induce an ${\CO}_X$-regular sequences on $X$.
 \bE
 \item 
 Let $\p_1, \ldots, \p_n$ be the associated primes of $S$, such that $\p_i\in X$.
  So, $\p_i$ corresponds to associated 
  primes of $S_{(f_j)}$. Suppose $J\subseteq \p_i$. This is impossible because $grade(\left(J_{(\p_i)}, {\CO}_{X, \p_i}\right)\geq k$. So, we can pick a homogeneous element $\f_1\in J^p-\bigcup \p_i$. 
 \item Now assume $\f_1, \ldots, \f_{u-1}\in J^p$, with $u\leq r$, has been selected so that, any of its subsequence, will induce a 
 regular ${\CO}_X$-sequence. Let 
 $$
 \CP=\left\{{\wp}\in X:\wp\in Ass\left(\frac{S}{(\f_{i_1}, \ldots, \f_{i_v})} \right): v\leq u-1, 1\leq i_1<\cdots<i_v\leq u-1 \right\}
 $$
Suppose $J\subseteq {\wp}\in {\CP}$. This is impossible because 
$grade\left(J_{(\wp)}, S_{(\wp)}\right)\geq k$. So, we  pick 
$$
\f_u\in J^p-\bigcup\left\{\wp: \wp\in {\CP} \right\}\qquad {\rm where~}\f_u~{\rm is~homogeneous.}
$$
 \eE 
 So, the existence of $\f_1, \f_2, \ldots, \f_k\in J^p$ is established. 
 In the case of point (\ref{FinInCMZX}), we do not go through this process, and instead use
 $$
 \f_1:=f_1^p, ~\f_2:=f_2^p, ~\ldots, ~\f_k=f_k^p \in J^p\subseteq  I
 $$
 which is also a regular sequence on $X$, by hypothesis. 
 Apply  (\ref{projVerFox}) to the extended complex ${\SF}_{\bullet}$. This way we obtain a map $\nu_{\bullet}:{\SE}_{\bullet} \lra {\SF}_{\bullet}$, satisfying all the properties in (\ref{projVerFox}). 
  Now the theorem follows by restricting $\nu_{\bullet}$ to $X$.
  \pic 
 $\eop$ 
 
 \vspace{3mm}
 The following small improvement of (\ref{Them23Foxby}) will be some  use in future.
 
 \bC\label{DurgaDoubrr}{\rm 
 Consider the hypotheses of (\ref{Them23Foxby}). Further assume that 
 ${\CZ}_0\LRf{{\SG}_{\bullet}} \in Coh_Z(X)$. Then the map $\nu_{\bullet}:{\SE}_{\bullet} 
 \lra {\SG}_{\bullet}$ in  (\ref{Them23Foxby}), can be constructed with the additional property that  that $\nu_0$ factors as follows:
 $$
\diagram 
{\SE}_0\ar[rr]\ar@{->>}[d]_{\nu_0} && {\CH}_0\LRf{{\SE}_{\bullet}}\ar@/^/[dll]^{\overline{\nu}_0}\\
{\CZ}_0\LRf{{\SG}_{\bullet}} && \\
\enddiagram
$$ 
 
 }
 \eC
 \pf Use the notations as in the proof of (\ref{Them23Foxby}, \ref{projVerFox}). We outline the adjustments needed to proof of (\ref{Them23Foxby}). There is no loss of generality, if we assume 
 ${\CZ}_0\LRf{{\SG}_{\bullet}}= {\SG}_0$. Let $K\subseteq S$ be a homogeneous ideal, such that $\tilde{K}=ann\LRf{{\SF}_0}$. It follows 
 $\supp{{\SF} \oplus \LRf{\bigoplus {\CH}_i\LRf{{\CF}_{\bullet}}}}\subseteq V(J)$. So, 
 $V(K\cap I) \subseteq V(J)$ and hence $J^p\subseteq K\cap I$ for some $p\gg 0$. So, the choices of $\f_1, \f_2, \ldots, \f_k \in J^p\subseteq K\cap I$. 
 
 Now we do the adjustment needed in the proof of (\ref{projVerFox}). In the diagram 
 (\ref{fzisZETR}), we have $\f_1 z=0$. It follows from this that the diagram (\ref{6Oct24Modi})  can be extended to the following
\begin{equation}\label{6Octfgr24Modi}
\diagram
{\CO}_X(-d_1)\ar[r]^{\f_1^m}\ar[d]_{\mu_z\f_1^m} & {\CO}_X\ar[d]^{\mu_z}\ar@{->>}[r]
&\frac{{\CO}_Y}{\widetilde{\LRf{\f_1^m}}} \ar[d]\ar[dl]^{\overline{\mu}_z}\\ 
{\SF}_1\ar[r] 
&  {\CZ}_0\left({\SF}_{\bullet}\right)\ar@{->>}[r]&{\CH}_0\left({\SF}_{\bullet} \right)\\
\enddiagram\quad {\rm commutative~diagram}\qquad \forall m\geq 1
\end{equation} 
Therefore, the diagram 
$$
\diagram
{\CK}_1\LRf{\f_1^m, \ldots, \f_k^m}\ar[r]\ar[d]_{\mu_1} & {\CO}_X\ar[d]^{\mu_0}\ar@{->>}[r]
&\frac{{\CO}_Y}{\widetilde{\LRf{\f_1^m, \ldots, \f_k^m}}} \ar[d]\ar[dl]^{\overline{\mu}_0}\\ 
{\SF}_1\ar[r] 
&  {\CZ}_0\left({\SF}_{\bullet}\right)\ar@{->>}[r]&{\CH}_0\left({\SF}_{\bullet} \right)\\
\enddiagram
\quad {\rm commutes} \quad \forall m\geq k
$$
The map $\nu_{0}^m$ is obtained by taking direct sum of such diagrams. This establishes the corollary. \pic $\eop$



\vspace{3mm}
The following basic result is a consequence of (\ref{Them23Foxby}), which is a refinement of \cite[Lem 22]{M15}. 
\bC\label{112SurCMZ}{\rm 
As in (\ref{Them23Foxby}), let $X$ be a quasi projective scheme over $\spec{A}$,
and $Z\subseteq X$ be as in the setup  (\ref{CIsetUp}). 
Let ${\SF}\in Coh^Z(X)$. Then there is a surjective map
$$
\diagram
 {\SE} \ar@{->>}[r]^{\f} & {\SF}\\
\enddiagram \qquad {\rm with}\qquad {\SE}\in C{\BM}^Z(X).
$$
}
\eC
\pf It is immediate application of (\ref{Them23Foxby}) to the complex 
$\diagram 0 \ar[r] & {\SF} \ar[r] & 0\\ \enddiagram$. \pic $\eop$

\section{Derived Equivalences}\label{equivSec} 

In this section, we prove two equivalence theorems (\ref{32Thm1007}, \ref{EquivKeller}) and summarize them in Theorem \ref{EquiALL}.
Before we proceed we recall the following   from \cite[Lem 3.2]{M15} that is needed to prove the first equivalence Theorem \ref{32Thm1007}.

\bL\label{Lem3p2Obs}{\rm 
Let ${\SA}$ be an abelian category and ${\SF}_{\bullet}, {\SG}_{\bullet}$ be two objects in the bounded derived category ${\bfD}^b({\SA})$. Assume, for some  $n_0\in {\BZ}$, we have
$$
\left\{\begin{array}{ll}
{\CH}_r\left({\SF}_{\bullet}\right)=0 & \forall r\leq n_0-1\\
{\CH}_r\left({\SG}_{\bullet}\right)=0 & \forall r\geq n_0\\
\end{array}\right.
$$
Then $Mor_{{\bfD}^b({\SA})}\left({\SF}_{\bullet}, {\SG}_{\bullet} \right)=0$. Further, if ${\SV}\subseteq {\SA}$ is a 
resolving subcategory, then the same holds for ${\bfD}^b({\SV})$.
}
\eL

\vspace{3mm}
The following is the statement of our first equivalence theorem. 
 \bT\label{32Thm1007}{\rm 
 Suppose $X$ is a quasi projective scheme over $\spec{A}$. 
 Use other notations, as in section \ref{SecVorAus}. 
 So, $X\subseteq Y:=\proj{S}$ is open, as in (\ref{setUP}).
 Let $Z\subseteq X$ be a closed subset.
  Consider the commutative diagram of natural functors
 $$
 \diagram 
 {\bfD}^b\left({\BM}^Z(X)\right)
\ar[r]^{\iota_Z} & {\SD}^b_Z\left({{\BM}}(X)\right) \ar@{^(->}[d]  &{\SD}^b_Z\left({{\SV}}(X)\right)\ar[l]_{\zeta_Z} \ar@{^(->}[d]  \\
&  {\bfD}^b\left({\BM}(X)\right) &  {\bfD}^b\left({\SV}(X)\right)\ar[l]^{\zeta}\\
 \enddiagram
 $$
 Then $\zeta$ and $\zeta_Z$ are  are natural equivalences. 
 Further, if $Z$ is complete intersection with $grade(Z, X)=k$, as in the  setup (\ref{CIsetUp}) then $\iota_Z$ is also a natural equivalence.
 }
 \eT 
 \pf 
 It is well known that the functor $\zeta:{\bfD}^b\left({\SV}(X)\right) \iso {\bfD}^b\left({\BM}(X)\right)$ is an 
 equivalence (e.g. one can use \cite[Sec. 1.5]{K99}). Restricting this equivalence
 to ${\SD}_Z^b\LRf{{\SV}(X)}$, it is clear that 
 $\zeta_Z$ is an equivalence. For future reference we comment that inverse of $\zeta_Z$ is given by the double complex trick. 

So, now we prove that $\iota_Z$ is an equivalence, when $Z$ is a complete intersection as in  (\ref{CIsetUp}).
 We write the proof in the following propositions. See remark \ref{1445DONe}.
 $\eop$

  \vspace{3mm} 
 The following is an analogue of \cite[Prop 3.4]{M15}.
 \bP
 \label{M15kaPrp34}
 {\rm 
 Use the  notations as in Theorem \ref{32Thm1007}. Assume   $Z$ is complete intersection with $grade(Z, X)=k$, as in  (\ref{CIsetUp}). Then, 
 $\iota_Z$ is essentially surjective and full. 
 }
 \eP
 \pf We will write $\iota:=\iota_Z$. Note that $\iota$ is an inclusion, at object level. 
 For ${\SF}_{\bullet} \in Ch^b(Coh(X))$, we say
 $$
 width\left({\SF}_{\bullet}\right) \leq r, \quad if\quad {\CH}_p\left({\SF}_{\bullet}\right)=0, 
 ~~\forall~p\notin [n, n+r], \quad {\rm for~some}\quad n \in {\BZ}. 
 $$
 We use induction on the width $r$, and prove the following:
 \bE
 \item\label{1038One} Let ${\SF}_{\bullet} \in {\SD}^b_Z\LRf{{\BM}(X)}$ and 
 $width\left({\SF}_{\bullet}\right)\leq r$. Then ${\SF}_{\bullet}\cong \iota_Z\left(\tilde{{\SF}}_{\bullet}\right)$, for some $\tilde{{\SF}}_{\bullet} \in Obj\left({\bfD}^b\left({\BM}^Z(X) \right)\right)$.
 
  \item\label{1040Two} Given ${\SF}_{\bullet}, {\SG}_{\bullet} \in Obj\left({\bfD}^b\left({\BM}^Z(X) \right)\right)$, with $width\left({\SF}_{\bullet}\bigoplus {\SG}_{\bullet} \right)\leq r$, the map
 $$
  Mor_{{\bfD}^b\left({\BM}^Z(X) \right)}\left({\SF}_{\bullet}, {\SG}_{\bullet} \right)
  \sur
    Mor_{{\SD}^b_Z\left({\BM}(X) \right)}\left({\SF}_{\bullet}, {\SG}_{\bullet} \right)
 $$
 is surjective. 
 \eE 
 Let ${\SF}_{\bullet} \in {\SD}^b_Z\LRf{{\BM}(X)}$  
 and $width\left({\SF}_{\bullet} \right)\leq r=0$. We can assume $n=0$ and hence 
 ${\CH}_i\left({\SF}_{\bullet}\right)=0,~\forall~i\neq 0$.  By replacing ${\SF}_0$ by
  $Z_0=\ker\left({\SF}_0\lra {\SF}_{-1}\right)\in {\BM}(X)$, we can assume that ${\SF}_i=0,~\forall i\leq -1$. Note, $B_0=Image\left({\SF}_1\lra {\SF}_0\right)\in {\BM}(X)$. So, 
  ${\CH}_0\left({\SF}_{\bullet}\right)=\frac{Z_0}{B_0}\in  {\BM}^Z(X)$. 
  As usual, we treat ${\CH}_0\left({\SF}_{\bullet}\right)$ as a complex in ${\bfD}^b\left({\BM}^Z(X) \right)$,
  concentrated at degree zero. 
  It follows that 
  the map ${\SF}_{\bullet} \lra \iota\left({\CH}_0\left({\SF}_{\bullet}\right)\right)$ is a quasi isomorphism, in $Ch^b\left({\BM}(X)\right)$. So, (\ref{1038One}) is established, when $r=0$.
  
  Now suppose ${\SF}_{\bullet}, {\SG}_{\bullet} \in Obj\left({\bfD}^b\left({\BM}^Z(X) \right)\right)$
  and 
  $width\left({\SF}_{\bullet}\bigoplus {\SG}_{\bullet} \right)\leq r=0$. We can assume 
  $$
  {\CH}_i\left({\SF}_{\bullet}\right)=  {\CH}_i\left({\SG}_{\bullet}\right)=0\quad \forall~i\neq 0
  $$
  Let $f_{\bullet}\in Mor_{{\SD}^b_Z\left({\BM}(X) \right)}\left({\SF}_{\bullet}, {\SG}_{\bullet} \right)$. We can write
  $$
  f_{\bullet}
  =g_{\bullet}t_{\bullet}^{-1}:
  \diagram
{\SF}_{\bullet}  & {\SW}_{\bullet}\ar[l]_{t_{\bullet}} \ar[r]^{g_{\bullet}}  &{\SG}_{\bullet}\\
  \enddiagram 
  \quad {\rm where}\quad t_{\bullet}~~{\rm is~a~quasi ~isomorphism~in}~Ch^b\left(({\BM}^Z(X)\right)
  $$
  As before, we can assume ${\SF}_i={\SW}_i={\SG}_i=0~\forall i\leq -1$.
  Now, the commutative diagram 
 $$
 \diagram 
{\SF}_{\bullet} \ar[d]_{\wr}^{\alpha} & {\SW}_{\bullet}\ar[l]_{t_{\bullet}} \ar[r]^{g_{\bullet}}\ar[d]^{\wr}  &{\SG}_{\bullet}\ar[d]^{\wr}_{\beta}\\
 {\CH}_0\left({\SF}_{\bullet}\right) &  {\CH}_0\left({\SW}_{\bullet}\right)
 \ar[l]^{{\CH}_0\left(t_{\bullet}\right)} \ar[r]_{{\CH}_0\left(g_{\bullet}\right)} &  {\CH}_0\left({\SG}_{\bullet}\right)\\
 \enddiagram 
 $$
 gives a lift of $f_{\bullet}$. 
 In other words, 
 $f_{\bullet}=\iota\left(\beta^{-1}
 {\CH}_0\left(g_{\bullet}\right){\CH}_0\left(t_{\bullet}\right)^{-1}\alpha\right)
 $.
 So, (\ref{1040Two}) is established, when $r=0$.
 
 For the inductive step, assume $r\geq 1$.  Let
 ${\SF}_{\bullet}\in Obj({\SD}^b_Z\left({\BM}(X)\right)$, with $width\left({\SF}_{\bullet}\right)\leq r$. 
 We can assume 
 $$
 {\CH}_i\left({\SF}_{\bullet}\right) = 0 \quad 
 unless ~~r\geq i\geq 0\quad {\rm and}\quad 
 {\CH}_0\left({\SF}_{\bullet}\right)\neq 0, \quad \supp{{\CH}_i\left({\SF}_{\bullet}\right)}\subseteq Z\quad \forall i.
 $$
 As before, we can assume that ${\SF}_i=0, ~\forall~i\leq -1$. 
 Since $Z$ is complete intersection, with $grade(Z, X)=k$, by Theorem \ref{Them23Foxby}(\ref{FinInCMZX}),
 there is a morphism $\nu_{\bullet}:{\SE}_{\bullet}\lra {\SF}_{\bullet}$ such that 
 \bE
 \item the complex ${\SE}_{\bullet}\in{\bfD}^b\left({\SV}(X)\right)$,
 ${\SE}_i=0$ unless  $0\leq i \leq k:=grade(Z, X)$. 
 \item ${\CH}_i\left({\SE}_{\bullet}\right)=0, ~\forall i\neq 0$. 
 \item 
 ${\CH}_0\left({\SE}_{\bullet}\right)\in C{\BM}^Z(X)$. 
 \item Finally, the map ${\CH}_0\left(\nu_{\bullet}\right): {\CH}_0\left({\SE}_{\bullet}\right)
 \sur {\CH}_0\left({\SF}_{\bullet}\right)$ is surjective.
 \eE
 As before, we consider 
  ${\CH}_0\left({\SE}_{\bullet}\right)\in {\bfD}^b\left({\BM}^Z(X)\right)$,
 concentrated at degree zero. Then $\iota\left({\CH}_0\left({\SE}_{\bullet}\right)\right)\cong
 {\SE}_{\bullet}\in {\SD}^b_Z\left({\BM}(X)\right)$. 
 Embed $\nu_{\bullet}$ in an exact 
 triangle 
\begin{equation}\label{1106Tringle}
 \diagram 
 T^{-1}\Delta_{\bullet} \ar[r]^{\quad\mu_{\bullet}} & {\SE}_{\bullet} \ar[r]^{\nu_{\bullet}} & {\SF}_{\bullet} \ar[r]^{\eta_{\bullet}} & \Delta_{\bullet}\\
 \enddiagram \qquad in \quad {\SD}^b_Z\LRf{{\BM}(X)}
 \end{equation}
 The corresponding homology exact sequence yields the following:
 \bE
 \item The sequence
 $$
 \diagram 
 0\ar[r] &{\CH}_1\left({\SF}_{\bullet} \right) \ar[r] &{\CH}_1\left(\Delta_{\bullet} \right) \ar[r] &
 {\CH}_0\left({\SE} _{\bullet} \right) \ar[r]^{{\CH}_0\left(\nu_{\bullet}\right)} &
 {\CH}_0\left({\SF}_{\bullet}\right) \ar[r] &0\\
 \enddiagram 
 $$
 is exact, because ${\CH}_0\left(\nu_{\bullet}\right)$ is surjective.
 \item ${\CH}_i\left(\Delta_{\bullet} \right)=0~\forall  i\leq 0$.
 \item ${\CH}_i\left(\Delta_{\bullet} \right)\cong {\CH}_i\left({\SF}_{\bullet} \right),~\forall i\geq 2$. 
 \eE
 So, it follows $width\left(\Delta_{\bullet}\right)\leq r-1$. 
 By induction, we obtain the following:
 \bE
 \item There is $\widetilde{\Delta}_{\bullet}\in {\bfD}^b\left({\BM}^Z(X)\right)$ such that 
 $\iota\left(\widetilde{\Delta}_{\bullet}\right)\cong \Delta_{\bullet}$ in ${\SD}^b_Z\left({\BM}(X)\right)$. 
 \item As before consider ${\CH}_0\left({\SE}_{\bullet}\right)\in {\bfD}^b\left({\BM}^Z(X)\right)$,
 concentrated at degree zero. Then $width\left({\CH}_0\left({\SE}_{\bullet}\right)\bigoplus 
 T^{-1}\Delta_{\bullet} \right)\leq r-1$. By induction (\ref{1040Two}), there is a  map
 $$
 \diagram 
 T^{-1}\widetilde{\Delta}_{\bullet} \ar[r]^{\varphi_{\bullet}} & {\CH}_0\left({\SE}_{\bullet}\right)\\
 \enddiagram\quad \ni \quad 
 \iota\left(\varphi_{\bullet}\right)=\mu_{\bullet} \qquad {\rm as~in~the~triangle~(\ref{1106Tringle})}.
 $$
  \eE 
 Now, embed $\varphi_{\bullet}$ in an exact triangle
 \begin{equation}\label{4211Tringle}
 \diagram 
 T^{-1}\widetilde{\Delta}_{\bullet} \ar[r]^{\varphi_{\bullet}} & {\CH}_0\left({\SE}_{\bullet}\right) 
\ar[r]^{\psi_{\bullet}} & {\SG}_{\bullet} \ar[r]^{\beta_{\bullet}} & \widetilde{\Delta}_{\bullet}\\
 \enddiagram 
 \qquad {\rm in}\quad {\bfD}^b\left({\BM}^Z(X)\right).
 \end{equation}
 Take the image under $\iota$ and compare these two  (\ref{1106Tringle}, \ref{4211Tringle}) triangles
 $$
 \diagram 
 T^{-1}\Delta_{\bullet} \ar[r]^{\mu_{\bullet}}\ar[d]_{\wr} & {\SE}_{\bullet} \ar[r]^{\nu_{\bullet}} \ar[d]_{\wr} & {\SF}_{\bullet} \ar[r]^{\eta_{\bullet}} \ar@{-->}[d]_{\wr} & \Delta_{\bullet}\ar[d]^{\wr} \\
T^{-1}\iota\left(\widetilde{\Delta}_{\bullet}\right) \ar[r]_{\mu_{\bullet}} & \iota\left({\CH}_0\left({\SE}_{\bullet}\right)\right) 
\ar[r]^{\psi_{\bullet}} & \iota\left({\SG}_{\bullet} \right)\ar[r]^{\iota\left(\beta_{\bullet}\right)} & \iota\left(\widetilde{\Delta}_{\bullet}\right)\\
 \enddiagram 
 $$
 The solid vertical arrows induce the broken arrow. 
 Since the solid vertical arrows are equivalences, so is the broken arrow. This completes the 
 proof of the inductive step of (\ref{1038One}). 
 
To prove the inductive step for (\ref{1040Two}), 
let ${\SF}_{\bullet}$, ${\SG}_{\bullet}$ be complexes in 
${\bfD}^b\left({\BM}^Z(X)\right)$ and let
$f_{\bullet}\in Mor_{{\SD}^b_Z\left({\BM}(X)\right)}\left({\SF}_{\bullet}, {\SG}_{\bullet}\right)$, and 
$width\left({\SF}_{\bullet} \bigoplus {\SG}_{\bullet}\right)\leq r$. 
 Assume ${\SF}_i={\SG}_i=0,~\forall i\leq -1$ and, either ${\CH}_0\left({\SF}_{\bullet}\right)\neq 0$
or  ${\CH}_0\left({\SG}_{\bullet}\right)\neq 0$.
In either case,
by Theorem \ref{Them23Foxby}(\ref{FinInCMZX}),
 there are morphisms $\nu_{\bullet}:{\SE}_{\bullet}\lra {\SF}_{\bullet}$,
 $\mu_{\bullet}:{\SL}_{\bullet}\lra {\SG}_{\bullet}$ of complexes such that 
 \bE
 \item   Unless $0\leq i\leq k=grade(Z, X)$,   
 ${\SE}_i=0$, ${\SL}_i=0$. Further,
 the complexes ${\SE}_{\bullet}, {\SL}_{\bullet}\in{\SD}^b_Z\left({\SV}(X)\right)$,

 \item ${\CH}_i\left({\SE}_{\bullet}\right)=0, {\CH}_i\left({\SL}_{\bullet}\right)=0~\forall i\neq 0$. 
 \item We have
$$
 \left\{ \begin{array}{l}
 {\CH}_0\left({\SE}_{\bullet}\right)\in C{\BM}^Z(X)\\
  {\CH}_0\left({\SL}_{\bullet}\right)\in C{\BM}^Z(X)\\
  \end{array}\right.
 $$. 
 \item Finally, 
 $$
  \left\{ \begin{array}{l}
  {\CH}_0\left(\nu_{\bullet}\right): {\CH}_0\left({\SE}_{\bullet}\right)
 \sur {\CH}_0\left({\SF}_{\bullet}\right)\\
   {\CH}_0\left(\mu_{\bullet}\right): {\CH}_0\left({\SL}_{\bullet}\right)
 \sur {\CH}_0\left({\SG}_{\bullet}\right)\\
   \end{array}\right.
   \quad {\rm are~surjective}.
 $$ 

 \eE
  Consider the following commutative diagrams
$$
\diagram 
{\SE}_1 \ar[r] \ar[d]_{\nu_1}& {\SE}_0 \ar[r] \ar[d]_{\nu_0}& {\CH}_0\left({\SE}_{\bullet} \right)\ar[r] \ar[d]^{{\CH}_0(\nu_{\bullet})}& 0\\
{\SF}_1 \ar[r] & {\SF}_0 \ar[r] & {\CH}_0\left({\SF}_{\bullet} \right) \ar[r]& 0\\
\enddiagram 
\qquad 
\diagram 
{\SL}_1 \ar[r] \ar[d]_{\mu_1}& {\SL}_0 \ar[r] \ar[d]_{\mu_0}& {\CH}_0\left({\SL}_{\bullet} \right)\ar[r] 
\ar[d]^{{\CH}_0(\mu_{\bullet})}& 0\\
{\SG}_1 \ar[r] & {\SG}_0 \ar[r] & {\CH}_0\left({\SG}_{\bullet} \right) \ar[r]& 0\\
\enddiagram 
$$
%
%
It follows from Corollary \ref{DurgaDoubrr} of 
 Theorem \ref{Them23Foxby}, 
the maps $\nu_0$, ${\CH}_0\left({\nu}_{\bullet}\right)$ factors through 
a surjective map
${\CH}_0\left({\SE}_{\bullet}\right) \lra {\SF}_0$ and likewise for $\mu_0$,  ${\CH}_0\left({\mu}_{\bullet}\right)$, as in the diagrams:
%
$$
\diagram 
{\SE}_1 \ar[r] \ar[d]_{\nu_1}& {\SE}_0 \ar[r] \ar[d]_{\nu_0}& {\CH}_0\left({\SE}_{\bullet} \right)\ar[r] \ar[d]^{{\CH}_0(\nu_{\bullet})}\ar@{->>}[dl]^{\overline{\nu}}& 0\\
{\SF}_1 \ar[r] & {\SF}_0 \ar[r] & {\CH}_0\left({\SF}_{\bullet} \right) \ar[r]& 0\\
\enddiagram 
\diagram
{\SL}_1 \ar[r] \ar[d]_{\mu_1}& {\SL}_0 \ar[r] \ar[d]_{\mu_0}& {\CH}_0\left({\SL}_{\bullet} \right)\ar[r] 
\ar[d]^{{\CH}_0(\mu_{\bullet})}\ar@{->>}[dl]^{\overline{\mu}}& 0\\
{\SG}_1 \ar[r] & {\SG}_0 \ar[r] & {\CH}_0\left({\SG}_{\bullet} \right) \ar[r]& 0\\
\enddiagram 
$$
Now treat ${\CH}_0\left({\SE}_{\bullet}\right)$ and ${\CH}_0\left({\SL}_{\bullet}\right)$, 
as two complexes in ${\bfD}^b\left({\BM}^Z(X)\right)$, concentrated at degree zero. The maps $\overline{\nu}, \overline{\mu}$ extend to maps of complexes $\overline{\nu}_{\bullet}, \overline{\mu}_{\bullet}$,
as in the commutative diagrams
$$
\diagram 
& {\CH}_0\left({\SE}_{\bullet}\right) \ar[d]^{\overline{\nu}_{\bullet}}\\
{\SE}_{\bullet} \ar[r]_{\nu_{\bullet}}\ar@/^/[ru] & {\SF}_{\bullet}\\
\enddiagram
\qquad {\rm and}\quad 
\diagram 
& {\CH}_0\left({\SL}_{\bullet}\right) \ar[d]^{\overline{\mu}_{\bullet}}\\
{\SL}_{\bullet} \ar[r]_{\mu_{\bullet}} \ar@/^/[ru]& {\SG}_{\bullet}\\
\enddiagram
\quad {\rm which~ reside~ in}~{\SD}_Z\left({\BM}(X)\right)
$$
By replacing ${\SL}_{\bullet}$ by ${\SE}_{\bullet}\bigoplus {\SL}_{\bullet}$ and $\mu_{\bullet}:=f_{\bullet}\nu_{\bullet}+\mu_{\bullet}$, we can assume that the diagram
$$
\diagram 
{\SE}_{\bullet}\ar[r]^{\nu_{\bullet}} \ar@{^(->}[d]_{\Phi_{\bullet}}& {\SF}_{\bullet}\ar[d]^{f_{\bullet}}\\
{\SL}_{\bullet}\ar[r]_{\mu_{\bullet}} & {\SG}_{\bullet}\\
\enddiagram 
\qquad \qquad {\rm commutes.} 
$$
This extends to the commutative diagram
$$
\diagram 
{\SE}_{\bullet}\ar[rr]^{\nu_{\bullet}} \ar@{^(->}[dd]_{\Phi_{\bullet}}\ar[dr]&& {\SF}_{\bullet}\ar[dd]^{f_{\bullet}}  \\
&{\CH}_0\left({\SE}_{\bullet}\right)  \ar@/_/[ru]_{\overline{\nu}_{\bullet}}\ar[dd]&\\ 
{\SL}_{\bullet}\ar[rr]_{\mu_{\bullet}} \ar[dr]&& {\SG}_{\bullet} \\
&{\CH}_0\left({\SL}_{\bullet}\right) \ar@/_/[ru]_{\overline{\mu}_{\bullet}} &\\ 
\enddiagram 
$$
The middle vertical map is $\varphi_{\bullet}:= {\CH}_0\LRf{\Phi_{\bullet}}$.
Embed $\overline{\nu}_{\bullet}$ and $\overline{\mu}_{\bullet}$ in exact triangles
 in ${\bfD}^b\left({\BM}^Z(X)\right)$, as follows:
$$
\diagram 
T^{-1}\Delta_{\bullet} \ar[r] & {\CH}_0\left({\SE}_{\bullet}\right) \ar[r]^{\quad \overline{\nu}_{\bullet}} 
& {\SF}_{\bullet} \ar[r] 
& \Delta_{\bullet}\\
T^{-1}\Gamma_{\bullet} \ar[r] & {\CH}_0\left({\SL}_{\bullet}\right) \ar[r]_{\quad \overline{\mu}_{\bullet}} & {\SG}_{\bullet} \ar[r] & \Gamma_{\bullet}\\
\enddiagram
$$
Since $f_{\bullet}$ is in ${\SD}^b_Z\left({\BM}(X)\right)$ we extend it, as follows:
\begin{equation}\label{1283OnTrfg}
\diagram 
T^{-1}\Delta_{\bullet} \ar[r] \ar@{-->}[d]_{T^{-1}\eta_{\bullet}}& {\CH}_0\left({\SE}_{\bullet}\right) \ar[d]^{\varphi_{\bullet}}
\ar[r]^{\overline{\nu}_{\bullet}} 
& {\SF}_{\bullet} \ar[r] \ar[d]^{f_{\bullet}}& \Delta_{\bullet}\ar@{-->}[d]^{\eta_{\bullet}}\\
T^{-1}\Gamma_{\bullet} \ar[r] & {\CH}_0\left({\SL}_{\bullet}\right) \ar[r]_{\overline{\mu}_{\bullet}} & {\SG}_{\bullet} \ar[r] & \Gamma_{\bullet}\\
\enddiagram
\qquad {\rm in} \quad {\SD}^b_Z\left({\BM}(X)\right).
\end{equation} 
Here $\eta_{\bullet}$  in ${\SD}^b_Z\left({\BM}(X)\right)$ is obtained, by properties of triangulated categories. 
The induction hypothesis applies to $\eta_{\bullet}$. 
So, there are morphisms
$$
\left\{\begin{array}{l}
\tilde{\varphi}_{\bullet}:=\varphi_{\bullet}
:{\CH}_0\left({\SE}_{\bullet}\right)\lra {\CH}_0\left({\SL}_{\bullet}\right)\\
\tilde{\eta}_{\bullet}:\Delta_{\bullet} \lra \Gamma_{\bullet}\\
\end{array}\right. 
\quad {\rm in}\quad {\bfD}^b\left({\BM}^Z(X)\right) ~~\ni~~
\left\{\begin{array}{l}
\iota\left(\tilde{\varphi}_{\bullet}\right)=\varphi_{\bullet}\\
\iota\left(\tilde{\eta}_{\bullet}\right)=\eta_{\bullet}\\
\end{array}\right. 
$$
Now consider the diagram
\begin{equation}\label{1306Ledfsd}
\diagram 
T^{-1}\Delta_{\bullet} \ar[r] \ar[d]_{T^{-1}\tilde{\eta}}& {\CH}_0\left({\SE}_{\bullet}\right) 
 \ar[d]_{\tilde{\varphi}}\ar[r]^{\overline{\nu}_{\bullet}} 
& {\SF}_{\bullet} \ar[r] \ar@{-->}[d]^{g_{\bullet}}
& \Delta_{\bullet}\ar[d]^{\tilde{\eta}}\\
T^{-1}\Gamma_{\bullet} \ar[r] & {\CH}_0\left({\SL}_{\bullet}\right) \ar[r]_{\overline{\mu}_{\bullet}} & {\SG}_{\bullet} \ar[r] & \Gamma_{\bullet}\\
\enddiagram
\end{equation}
 It is clear that the left hand rectangle commutes in ${\SD}_Z^b\left({\BM}(Z)\right)$. By Lemma \ref{Lem3p2Obs}, it follows that the left hand square commutes in ${\bfD}^b\left({\BM}^Z(X)\right)$. By properties of 
 the triangulated category ${\bfD}^b\left({\BM}^Z(X)\right)$, it follows that there is a map
 $g_{\bullet}:{\SF}_{\bullet} \lra {\SG}_{\bullet}$ in ${\bfD}^b\left({\BM}^Z(X)\right)$, so that the above
 diagram commutes. 
 
 Now we claim that 
 $f_{\bullet}$ is in the image of $\iota$. Write $h=f_{\bullet}-\iota\left(g_{\bullet} \right)$. Apply $\iota$ to the diagram (\ref{1306Ledfsd}) and take the difference with diagram (\ref{1283OnTrfg}). This way we obtain the following diagram:
$$
\diagram 
T^{-1}\Delta_{\bullet} \ar[r] \ar[d]_0& {\CH}_0\left({\SE}_{\bullet}\right) \ar[d]_0
\ar[r]^{\overline{\nu}_{\bullet}} 
& {\SF}_{\bullet} \ar[r]^{\gamma} \ar[d]^{h}\ar@{-->}[dl]_{\omega}& \Delta_{\bullet}\ar@{-->}[d]^0
\ar@{-->}[dl]^{\varepsilon}\\
T^{-1}\Gamma_{\bullet} \ar[r] & {\CH}_0\left({\SL}_{\bullet}\right) \ar[r]_{\overline{\mu}_{\bullet}} & {\SG}_{\bullet} \ar[r] & \Gamma_{\bullet}\\
\enddiagram
$$
The maps $\omega$ and $\varepsilon$ are obtained by weak kernel and weak cokernel properties. 
\bE
\item Consider the case, ${\CH}_0\left({\SF}_{\bullet}\right)=0$. By Lemma 
\ref{Lem3p2Obs} $Mor_{{\SD}_Z^b\left({\BM}(X)\right)}
\left({\SF}_{\bullet}, {\CH}_0\left({\CL}_{\bullet}\right)\right)=0$. So, $\omega=0$ and hence $h=0$. Therefore, $\iota\LRf{g_{\bullet}}=f_{\bullet}$. 
So, the proposition is established whenever $width\left({\SF}\bigoplus {\SG}\right)\leq r$ and 
 ${\CH}_0\left({\SF}_{\bullet}\right)=0$.

\item Consider the case, ${\CH}_0\left({\SF}_{\bullet}\right)\neq 0$.
 In this case, ${\CH}_0\left(\Delta_{\bullet}\right)=0$,
and $width\left(\Delta_{\bullet}\bigoplus  {\SG}_{\bullet} \right)\leq r$. So, by the above case,
$$
Mor_{{\bfD}^b\left({\BM}^Z(X)\right)}\left(\Delta_{\bullet}, {\SG}_{\bullet}\right) \sur 
Mor_{{\SD}^b_Z\left({\BM}^Z(X)\right)}\left(\Delta_{\bullet}, {\SG}_{\bullet}\right) 
\qquad {\rm is~surjective.} 
$$
Therefore $\iota(\tilde{\varepsilon})= \varepsilon$,
for some $\tilde{\varepsilon}$. 
 Hence $h=\iota(\tilde{\varepsilon}\gamma)$ %
  is in the image of $\iota$. Hence $f_{\bullet}=\iota(g_{\bullet})+h=\iota\left((g_{\bullet}+\tilde{\varepsilon}\gamma\right)$ is in the image of $\iota$.

\eE 
\pic $\eop$

 \vspace{3mm}
Now we prove the faithfulness of $\iota_z$.
  
  \bP
  \label{030Faithful}{\rm 
 Use the notations, as in  (\ref{32Thm1007}), and let $Z$ is complete intersection as in (\ref{CIsetUp}). Then the functor 
 $\iota:=\iota_Z$ is faithful. 
 }
 \eP
\pf Let $f_{\bullet}:{\SF}_{\bullet} \lra {\SG}_{\bullet}$ is a map in ${\bfD}^b\left({\BM}^Z(X)\right)$,
and $\iota_Z\left(f_{\bullet}\right)=0$. We prove $f_{\bullet}=0$. We can assume that 
$f_{\bullet}$  is denominator free. So,  $f_{\bullet}$ is a morphism in $Ch^b\left({\BM}^Z(X)\right)$. Embed $f_{\bullet}$ in
an exact triangle
$$
\diagram
T^{-1}\Delta_{\bullet} \ar[r]^{g_{\bullet}} & {\SF}_{\bullet} \ar[r]^{f_{\bullet}} & {\SG}_{\bullet} \ar[r]^{h_{\bullet}}
& \Delta_{\bullet}\\
\enddiagram
$$
Applying $\iota_Z$, we obtain the exact triangle
$$
\diagram
T^{-1}\Delta_{\bullet} \ar[r]^{g_{\bullet}} & {\SF}_{\bullet} \ar[r]^{0} & {\SG}_{\bullet} \ar[r]^{h_{\bullet}}
& \Delta_{\bullet}\\
\enddiagram
\quad {\rm in}\quad {\SD}^b_Z\left({\BM}(X)\right).
$$
Using weak kernel property, $g_{\bullet}$ splits, as follows:
  $$
\diagram
&{\SF}_{\bullet}\ar[d]^1\ar@/_/[dl]_{\eta_{\bullet}} &&\\
T^{-1}\Delta_{\bullet} \ar[r]_{g_{\bullet}} & {\SF}_{\bullet} \ar[r]_{0} & {\SG_{\bullet}} \ar[r]_{h_{\bullet}}
& \Delta_{\bullet}\\
\enddiagram
\quad {\rm in}\quad {\SD}^b_Z\left({\BM}(X)\right).
$$
  By fullness (\ref{M15kaPrp34}), there is 
  $$
  \tilde{\eta}_{\bullet}\in Mor_{{\bfD}^b\left({\BM}^Z(X) \right)}\left({\SF}_{\bullet}, T^{-1}\Delta_{\bullet} \right)
  \qquad \ni \qquad \iota_Z\left( \tilde{\eta}_{\bullet}\right)=\eta_{\bullet}.
$$
Now we embed $g_{\bullet}\tilde{\eta}_{\bullet}$ in an exact triangle
$$
\diagram 
{\SF}_{\bullet}\ar[r]^{g_{\bullet}\tilde{\eta}_{\bullet}} & {\SF}_{\bullet} \ar[r] & \Gamma_{\bullet}\ar[r] &T{\SF}_{\bullet}\\
\enddiagram
\qquad {\rm in}\quad {\bfD}^b\left({\BM}^Z(X)\right).
$$
 Under $\iota_Z$, this exact triangle maps to 
 $$
\diagram 
{\SF}_{\bullet}\ar[r]^{1} & {\SF}_{\bullet} \ar[r] & \Gamma_{\bullet}\ar[r] &T{\SF}_{\bullet}\\
\enddiagram
\quad {\rm in}\quad {\SD}^b_Z\left({\BM}(X)\right).
$$
This means that $\Gamma_{\bullet}\cong 0$ is acyclic in ${\SD}^b_Z\left({\BM}(X)\right)$. In other words, 
it is exact, and hence is acyclic (exact) in  ${\bfD}^b\left({\BM}^Z(X)\right)$, as well. It follows now that $g_{\bullet}\tilde{\eta}_{\bullet}$ is an isomorphism in ${\SD}^b_Z\left({\BM}(X)\right)$.
Therefore $g_{\bullet}\left(\tilde{\eta}_{\bullet}\left(g_{\bullet}\tilde{\eta}_{\bullet}\right)^{-1}\right)=1_{{\SF}_{\bullet}}$. So, $g_{\bullet}$ splits in  ${\SD}^b_Z\left({\BM}(X)\right)$. Since 
$f_{\bullet}g_{\bullet}=0$, it follows 
$f_{\bullet}=0$ in  ${\SD}^b_Z\left({\BM}(X)\right)$.
\pic $\eop$ 


\vspace{3mm}
\begin{remark}\label{1445DONe}{\rm 
This completes the proof of Theorem \ref{32Thm1007}. 
 }
 \end{remark} 
 

\subsection{The category $C{\BM}^Z(X)$ of perfect modules} \label{SecCmZxCaT}
The importance of various categories of perfect modules in $K$-theory was already evident in \cite{M15, M17, M23} and others. In this article, we further work with that category of perfect modules $C{\BM}^Z(X)$, with support in $Z$. Its importance in hermitian theory is particularly crucial, because it  $C{\BM}^Z(X)$  has a natural duality.
In this subsection, we prove the equivalence 
 ${\bfD}^b\left(C{\BM}^Z(X) \right)\cong {\bfD}^b\left({\BM}^Z(X) \right)$
 of the derived categories, under the setup (\ref{CIsetUp}).
 We augment the equivalence in Theorem \ref{32Thm1007}, with this further equivalence.
%
 First, we establish some preparatory work.


 \bL
 \label{CWks24four}{\rm 
 Suppose $X$ is a quasi projective scheme over $\spec{A}$. 
 Use other notations, as in (\ref{nota}, \ref{setUP}). 
 So, $X\subseteq Y:=\proj{S}$ is open.
 Let $Z\subseteq X$ be a closed subset and $grade(Z, X)=k$.  Consider the exact sequence 
 $$
 \diagram 
 0\ar[r] & {\SK} \ar[r] & {\SE} \ar[r] & {\SF} \ar[r] & 0\\
 \enddiagram
 $$
 Suppose ${\SE}\in C{\BM}^Z\left(X\right)$ and ${\SF}\in {\BM}^Z\left(X\right)$. 
 Then $\PDV{{\SK}}\leq \max\left\{k, \PDV{{\SF}}-1\right\}$. In case $\max\left\{k, \PDV{{\SF}}-1\right\}=\PDV{{\SF}}-1$ then $\PDV{{\SK}}=\PDV{{\SF}}-1$. 
 }
 \eL
\pf Use Tor argument, locally. $\eop$ 

\vspace{3mm}
We record some basic properties of the category $C{\BM}^Z(X)$, as follows.
\bL
\label{BasicSMZX}{\rm 
 Suppose $X$ is a quasi projective scheme over $\spec{A}$. 
 Use other notations, as in (\ref{nota}, \ref{setUP}). 
 So, $X\subseteq Y:=\proj{S}$ is open.
 Let $Z\subseteq X$ be a closed subset and $grade(Z, X)=k$. 
 \bE
 \item\label{1550One} Suppose ${\SE}\in C{\BM}^Z(X)$ and ${\SL}\in {\SV}(X)$. 
 Then  ${\SE}\otimes {\SL}\in C{\BM}^Z(X)$.
 \item\label{1551TWO} Suppose ${\SF}\in Coh^Z(X)$ and $Z$ is complete intersection as in (\ref{CIsetUp}). Then there is a surjective map
 $\varphi:{\SE} \sur {\SF}$ where ${\SE}\in C{\BM}^Z(X)$.
 \eE 
 }
 \eL
 \pf First, (\ref{1550One}) is obvious. Now, (\ref{1551TWO}) is a restatement of Corollary \ref{112SurCMZ}.
 $\eop$

\vspace{3mm}

\bL
\label{CWks24four}{\rm 
 Suppose $X$ is a quasi projective scheme over $\spec{A}$. 
 Use other notations, as in (\ref{nota}, \ref{setUP}). 
 So, $X\subseteq Y:=\proj{S}$ is open, for some $Y$.
 Let $Z\subseteq X$ be a complete intersection closed subset as in (\ref{CIsetUp}).
 Let ${\SF}\in {\BM}^Z(X)$. Then there is a resolution
 \begin{equation}\label{CMka25Res}
 \diagram 
 0\ar[r] &{\SL}_n \ar[r] &{\SL}_{n-1}\ar[r]^{d_{n-1}}& \cdots \ar[r]& {\SL}_1\ar[r] &{\SL}_0\ar[r] &
 {\SF} \ar[r] & 0\\
 \enddiagram
 \quad \ni ~~ {\SL}_i\in C{\BM}^Z(X)
 \end{equation}
 where $n=\PDV{{\SF}}-k$. 
 }
 \eL
\pf Write $m=\PDV({\SF})$. Then $k=grade(Z, X) \leq grade({\SF})\leq m$. If $m=k$ then ${\SF}\in 
C{\BM}^Z(X)$. Therefore, assume $k\leq m-1$. 
By (\ref{BasicSMZX}) there is a surjective map $d_0:{\SL}_0 \sur {\SF}$, with 
${\SL}_0\in C{\BM}^Z(X)$. Let ${\SE}_0=\ker(d_0)$. So, we have the exact sequence
$$
\diagram
0\ar[r] & {\SE}_0\ar[r] & {\SL}_0 \ar[r]^{d_0} & {\SF} \ar[r] & 0\\
\enddiagram\qquad {\rm with}\qquad {\SE}_0\in {\BM}^Z(X).
$$
It follows, $\PDV{{\SE}_0}=m-1$. Again, by 
(\ref{BasicSMZX}) there is a surjective map $d_1:{\SL}_1 \sur {\SE}_0$. So, we have an exact sequence
$$
\diagram
0\ar[r] & {\SE}_1\ar[r]& {\SL}_1\ar[r]^{d_1} & {\SL}_0 \ar[r]^{d_0} & {\SF} \ar[r] & 0\\
\enddiagram\qquad {\rm with}\qquad {\SE}_1\in {\BM}^Z(X).
$$
Now the proof is complete, inductively. $\eop$

\vspace{3mm}
\bL
\label{CWks24four}{\rm 
 Suppose $X$ is a quasi projective scheme over $\spec{A}$. 
 Use other notations, as in (\ref{nota}, \ref{setUP}). 
 So, $X\subseteq Y:=\proj{S}$ is open for some $Y$.
 Let $Z\subseteq X$ be a closed subset and $grade(Z, X)=k$. 
 Consider an exact sequence
 $$
\diagram
0\ar[r] & {\SK}\ar[r] & {\SF} \ar[r] & {\SG} \ar[r] & 0\\
\enddiagram
\qquad \qquad {\rm in}\qquad {\BM}^Z(X). 
$$
\bE
\item\label{744One} If ${\SK}, {\SG}\in C{\BM}^Z(X)$ then ${\SF}\in C{\BM}^Z(X)$.
\item\label{174Two} If ${\SF}, {\SG}\in C{\BM}^Z(X)$ then ${\SK}\in C{\BM}^Z(X)$.
\eE 
 }
 \eL
 \pf Suppose ${\SK}, {\SG}\in C{\BM}^Z(X)$. It is immediate ${\SF}\in {\BM}^Z(X)$ and $k=\PDV({\SF})$. So, ${\SF}\in C{\BM}^Z(X)$. This settles (\ref{744One}).
 
 Now ${\SF}, {\SG}\in C{\BM}^Z(X)$. So, ${\SK}\in {\BM}^Z(X)$. It is also clear $\PDV{{\SK}}\leq k$.
 Since ${\SK}\in {\BM}^Z(X)$, $grade({\SF})\geq k$. So, $k\leq grade({\SF}) \leq \PDV{{\SF}}\leq k$.
 So, $\PDV{{\SF}}=k$. This settles (\ref{174Two}).
 \pic $\eop$ 
 
 \vspace{3mm}
The following casts the above in the mold of resolving subcategories.
\vspace{3mm}
\bC
\label{CWks24four}{\rm 
 Suppose $X$ is a quasi projective scheme over $\spec{A}$. 
 Use other notations, as in (\ref{nota}), (\ref{setUP}). 
 So, $X\subseteq Y:=\proj{S}$ is open.
 Let $Z\subseteq X$ be a complete intersection closed subset, as in (\ref{CIsetUp}). 
 
 Then $C{\BM}^Z(X) \subseteq Coh^Z(X)$ is a resolving subcategory (in the sense of \cite[Def 3.1]{M15d}). 
 Further, 
 $$
 {\BM}^Z(X) =\left\{{\SF}\in Coh^Z(X): \dim_{C{\BM}^Z(X)}{\SF}<\infty \right\}
 $$
  }
 \eC
 \pf Follows directly from (\ref{CWks24four}, \ref{BasicSMZX}). $\eop$



\vspace{3mm}
Now, we are ready to state and prove the  second equivalence theorem.
\bT\label{EquivKeller}{\rm 
Suppose $X$ is a quasi projective scheme over $\spec{A}$. 
 Use other notations, as in (\ref{nota}, \ref{setUP}). 
 So, $X\subseteq Y:=\proj{S}$ is open, for some $Y$. 
  Let $Z\subseteq X$ be a complete intersection closed subset, as in (\ref{CIsetUp}).
 Then the natural functor
 \begin{equation}\label{7892ndEqniso}
 \diagram 
 {\bfD}^b\left(C{\BM}^Z(X) \right)\ar[r]^{\sim} &  {\bfD}^b\left({\BM}^Z(X) \right)\\
 \enddiagram
 \quad {\rm is~an~equivalence.} 
 \end{equation}

 }
 \eT
\pf Consider the commutative diagram
\begin{equation}\label{807KleHy}
\diagram 
0\ar[r] & {\SK}_1 \ar@{-->}[r] \ar@{-->}[d]_{\q_2}&{\SE}_1 \ar@{-->}[r]^{\varepsilon\q}\ar@{-->>}[d]^{\q_1}& {\SF} \ar[r]\ar@{=}[d] & 0&\\
0\ar[r] & {\SK}_2 \ar[r] &{\SE}_2 \ar[r]_{\varepsilon} & {\SF} \ar[r] & 0& \quad {\rm in}\quad {\BM}^Z(X)\\
\enddiagram
\end{equation}
Here the second line is a given exact sequence in ${\BM}^Z(X)$, and ${\SF}\in C{\BM}^Z(X)$. 
By  (\ref{BasicSMZX}), there is a surjective map $\q_1: {\SE}_1 \sur {\SE}_2$, with 
${\SE}_1\in C{\BM}^Z(X)$. Let ${\SK}_1=\ker(\varepsilon\q_1)$. Then $\q_1$ induces a 
 map $\q_2:{\SK}_1\sur {\SK}_2$. In fact, $\q_2$ is surjective (use Snake Lemma). Further, 
 ${\SK}_1\in C{\BM}^Z(X)$, by (\ref{CWks24four}).

This shows that  \cite[Cond b, \S 1.5]{K99} is satisfied by the functor 
$\zeta: C{\BM}^Z(X) \lra {\BM}^Z(X)$. By (\ref{CWks24four}), any object ${\SF}\in {\BM}^Z(X)$ has a finite resolution by objects in $C{\BM}^Z(X)$. In other words,  \cite[Cond a, \S 1.5]{K99}
 is also satisfied by the functor 
$\zeta$. 

By \cite[Last para of \S 1.5]{K99} (the dual version), it follows that the map (\ref{7892ndEqniso}),
which is induced by $\zeta$, is an equivalence.
 \pic $\eop$ 
%


\vspace{3mm}
The following is the summary of the equivalences that we have established.

\bT\label{EquiALL}{\rm 
Suppose $X$ is a quasi projective scheme over $\spec{A}$. 
 Use other notations, as in (\ref{nota}, \ref{setUP}). 
 So, $X\subseteq Y:=\proj{S}$ is open, for some $Y$.
 Let $Z\subseteq X$ be a closed complete intersection subset, as in  (\ref{CIsetUp}).
 Then we have the following  
 equivalences of Derived categories:
  \begin{equation}\label{3rreeEnnEquiv}
 \diagram 
 {\bfD}^b\left(C{\BM}^Z(X)\right) \ar[r]_{\sim}^{\alpha} &  {\bfD}^b\left({\BM}^Z(X)\right) \ar[r]_{\sim}^{\beta}&
 {\SD}_Z^b\left({\BM}(X) \right) &  {\SD}_Z^b\left({\SV}(X)\right)\ar[l]^{\sim}_{\eta}  \\
 \enddiagram 
 \end{equation}
 All the functor $\alpha$, $\beta$ and $\eta$ are defined naturally. The inverse of $\eta$ is given by double complex trick. For future reference (\ref{d079verEqui4p6}), we note that there is a natural equivalence
 \begin{equation}\label{longEQUI74}
  \diagram 
 {\bfD}^b\left(C{\BM}^Z(X)\right) \ar[rr]_{\sim}^{\eta^{-1}\beta\alpha}
&&
 {\SD}_Z^b\left({\SV}(X)\right)  \\
 \enddiagram 
 \end{equation}
An object ${\SF}\in C{\BM}^Z(X)$ is treated as a complex concentrated at degree zero. We have
$$
\eta^{-1}\beta\alpha({\SF})={\rm a~resolution~of}~{\SF}
$$
 }
 \eT
 \pf This follows from Theorems \ref{EquivKeller} and \ref{32Thm1007}.
 \pic $\eop$ 
 
%

\section{The ${\BK}$-Theory Localization sequence} \label{SECBKTHEY}
In this section we discuss
the main results on ${\BK}$-Theory exact sequences, in this article. By 
${\BK}$-Theory, we mean non connective version of $K$-theory, which takes values in 
the category $\Sp$ of spectrum of pointed spaces \cite[pp 402]{M23}. We will not discuss connective $K$-theory. However, for idempotent complete exact categories ${\SE}$, the nonnegative groups $K_n\LRf{{\SE}}\cong {\BK}_n\LRf{{\SE}}~\forall n\geq 0$ are isomorphic \cite[pp 406]{M23}. 
Subsequently, for an exact category ${\SE}$ the ${\BK}$-theory spectra of ${\SE}$ will be denoted by 
${\BK}({\SE})$  \cite[Section 2.4.3]{S11}. Unless stated otherwise, we talk about ${\BK}$-theory of exact categories only. 

\bT\label{Equi1909ALL}{\rm 
Suppose $X$ is a quasi projective scheme over $\spec{A}$. 
 Use other notations, as in (\ref{nota}, \ref{setUP}). 
 Let $Z\subseteq X$ be a complete intersection closed subset, as in (\ref{CIsetUp}).  and 
 Write $U=X-Z$.
 Then the sequence
 \begin{equation}\label{190MainHSTH}
 \diagram 
  {\BK}\left(C{\BM}^Z(X)\right) \ar[r] &  {\BK}\left({\SV}(X)\right)  \ar[r]  &  
  {\BK}\left({\SV}(U)\right)  \\
 \enddiagram
 \qquad in \qquad \Sp
 \end{equation}
 is a homotopy  fibratiion of  spectra. 
 Consequently, $\forall n\in {\BZ}$, there is a long exact sequence of ${\BK}$-groups:
 \begin{equation}\label{1911HomFFibr}
 \diagram
 \cdots \ar[r] & {\BK}_n\left(C{\BM}^Z(X) \right)  \ar[r] & {\BK}_n\left({\SV}(X) \right)  \ar[r] 
 & {\BK}_n\left({\SV}(U) \right)  \ar[r] & {\BK}_{n-1}\left(C{\BM}^Z(X) \right)\ar[r] &\cdots\\
 \enddiagram
 \end{equation} 
 We further have homotopy equivalences:
 \begin{equation}\label{1918HomoEqui}
 \left\{\begin{array}{l}
 {\BK}\left(C{\BM}^Z(X) \right) \iso {\BK}\left({\BM}^Z(X) \right)\\
  {\BK}\left({\SV}(X)\right)\iso {\BK}\left({\BM}(X) \right)\\
 \end{array}\right.
 \end{equation} 
 }
 \eT 
\pf As in (\ref{nota}), for any exact subcategory ${\SB}(X)\subseteq Coh(X)$, denote the subcategories
$$
Ch^b_Z\left({\SB}(X)\right)=\left\{{\SE}_{\bullet}\in Ch^b\left({\SB}(X)\right): {\CH}_i\left({\SE}_{\bullet}\right)\in Coh^Z(X) ~\forall~i \right\}
$$
We use this notation for ${\SB}(X)={\SV}(X),~ {\BM}(X)$ and others. These are treated as complicial exact categories with weak equivalences (quasi isomorphisms). 
By Agreement theorem \cite[2.4.3, 3.2.30]{S11}, \cite[pp 409]{M23} it follows 
$$
\left\{\begin{array}{l}
{\BK}\left({\BM}\left(X\right) \right)\iso {\BK}\left(Ch^b\left({\BM}(X)\right)\right)\\
{\BK}\left({\BM}\left(U\right) \right)\iso {\BK}\left(Ch^b\left({\BM}(U)\right)\right)\\
{\BK}\left({\BM}^Z\left(X\right) \right)\iso {\BK}\left(Ch^b\left({\BM}^Z(X)\right)\right)\\
{\BK}\left(C{\BM}^Z\left(X\right) \right)\iso {\BK}\left(Ch^b\left(C{\BM}^Z(X)\right)\right)\\
{\BK}\left({\SV}\left(X\right) \right)\iso {\BK}\left(Ch^b\left({\SV}(X)\right)\right)\\
{\BK}\left({\SV}\left(U\right) \right)\iso {\BK}\left(Ch^b\left({\SV}(U)\right)\right)\\
\end{array}\right.
\quad 
{\rm are~ homotopy~ equivalences~in } ~\Sp.
$$
Now consider the sequence of functors
$$
\diagram 
Ch^b\left(C{\BM}^Z(X) \right)\ar[r] &
Ch^b\left({\BM}^Z(X) \right)\ar[r] & Ch^b_Z\left({\BM}(X)\right)\\
\enddiagram 
$$
of complicial exact categories. The functors of associated exact categories 
$$
\diagram 
{\bfD}^b\left(C{\BM}^Z(X) \right)\ar[r]^{\sim} &
{\bfD}^b\left({\BM}^Z(X) \right)\ar[r]^{\sim} & {\SD}^b_Z\left({\BM}(X)\right)\\
\enddiagram 
$$
is an equivalences, by Theorem \ref{EquiALL}. 
By Triangular equivalence (\cite[Thm 3.2.29]{S11}, \cite[pp 409]{M23}), it  follows that the following maps of ${\BK}$-theory spectra 
$$
\diagram 
{\BK}\left(Ch^b\left(C{\BM}^Z(X) \right)\right)\ar[r] & 
{\BK}\left(Ch^b\left({\BM}^Z(X) \right)\right)\ar[r] & {\BK}\left(Ch^b_Z\left({\BM}(X)\right)\right)\\
\enddiagram \quad in \quad \Sp
$$
are homotopy equivalences. Combining with the above, we have a sequence of homotopy equivalences
$$
\diagram
{\BK}\left({\BM}^Z(X\right)  \ar[r] &{\BK}\left(Ch^b_Z\left({\BM}(X)\right)\right)\ar[r] & {\BK}\left(Ch^b_Z\left({\BM}(X)\right)\right)\\
\enddiagram  \quad in \quad \Sp
$$
of ${\BK}$-theory spectra. 
Further,
the functor of complicial exact categories, with weak equivalences (i.e. the quasi isomorphisms)
$$
Ch^b_Z\left({\SV}(X)\right)\lra Ch^b_Z\left({\BM}(X)\right)\qquad {\rm induce~the~functor}\quad
\zeta: {\SD}^b_Z\left({\SV}(X)\right)\lra {\SD}^b_Z\left({\BM}(X)\right).
$$
The latter functor is also an equivalence, by Theorem \ref{32Thm1007}. Hence by the triangular equivalence theorem (\cite[Thm 3.2.29]{S11}, \cite[pp 409]{M23}), it follows
the map
$$
{\BK}\left(Ch^b_Z\left({\SV}(X)\right)\right)\iso {\BK}\left(Ch^b_Z\left({\BM}(X)\right)\right)
\qquad in \qquad \Sp
$$
is a homotopy equivalence. Putting all these together, 
$$
\diagram
{\BK}\left({\BM}^Z(X\right)  \ar[r] &{\BK}\left(Ch^b_Z\left({\BM}(X)\right)\right)\ar[r] & {\BK}\left(Ch^b_Z\left({\BM}(X)\right)\right) & {\BK}\left(Ch^b_Z\left({\BM}(X)\right)\right)\ar[l]\\
\enddiagram
$$
is a zig-zag  sequence of homotopy equivalences of spectra.
 Of course, homotopy inverse of the last map is also given by a functor,
by double complex. 

Now consider the sequences for functors
$$
\diagram 
Ch^b_Z\left({\SV}(X)\right)\ar[r] \ar[d]& Ch^b\left({\SV}(X)\right)\ar[r] \ar[d]&Ch^b\left({\SV}(U)\right)\ar[d]\\
Ch^b_Z\left({\BM}(X)\right)\ar[r] & Ch^b\left({\BM}(X)\right)\ar[r] &Ch^b\left({\BM}(U)\right)\\
\enddiagram 
$$
The induced diagram of the associated triangulated categories is as follows:
$$
\diagram 
{\SD}^b_Z\left({\SV}(X)\right)\ar[r] \ar[d]_{\wr}& {\bfD}^b\left({\SV}(X)\right)\ar[r] \ar[d]^{\wr}&{\bfD}^b\left({\SV}(U)\right)\ar[d]^{\wr}\\
{\SD}^b_Z\left({\BM}(X)\right)\ar[r] & {\bfD}^b\left({\BM}(X)\right)\ar[r] &{\bfD}^b\left({\BM}(U)\right)\\
\enddiagram 
$$
The lower line is exact up to factor (direct summand), by \cite[Prop A.4.7]{S11}. Since the vertical maps are equivalences, the first line is also  exact up to factor. 
Consider the diagram 
$$
\diagram 
&{\BK}\left({\SV}(X)\right)\ar[d]_{\wr}\ar[r]&{\BK}\left({\SV}(U)\right)\ar[d]^{\wr}\\
& {\BK}\left(Ch^b\left({\SV}(X)\right)\right)\ar[r] \ar[d]^{\wr}&{\BK}\left(Ch^b\left({\SV}(U)\right)\right)\ar[d]^{\wr}\\
{\BK}\left(Ch^b_Z\left({\BM}(X)\right)\right)\ar[r] & {\BK}\left(Ch^b\left({\BM}(X)\right)\right)\ar[r] &{\BK}\left(Ch^b\left({\BM}(U)\right)\right)\\
{\BK}\left({\BM}^Z\left(X\right) \right)\ar[u]^{\wr}&&\\
\enddiagram 
$$
of maps of spectra.
By Waldhausen-Thomason localization (\cite[3.2.27]{S11}, \cite[pp 408]{M23}) it follows that the third line is 
a homotopy fibration. Further the vertical maps are homotopy equivalences. Therefore, 
the sequence 
\begin{equation}\label{25SptHus}
 \diagram 
  {\BK}\left({\BM}^Z(X)\right) \ar[r] &  {\BK}\left({\SV}(X)\right)  \ar[r]  &  
  {\BK}\left({\SV}(U)\right)  \\
 \enddiagram
 \qquad in \qquad \Sp
\end{equation}
 is a homotopy fibration. 
By similar argument, using the Agreement theorem \cite[3.2.30]{S11} and  \cite[Thm 3.2.29]{S11}, 
 the following  are equivalences of homotopy spectra:
$$
\diagram 
{\BK}\left({\SV}(X)\right) \ar[d]_{\wr} &{\BK}\left({\BM}(X)\right) \ar[d]^{\wr}  \\
{\BK}\left(Ch^b\left({\SV}(X)\right)\right) \ar[r]_{\sim}& {\BK}\left(Ch^b\left({\BM}(X)\right)\right) \\
\enddiagram
\qquad 
\diagram 
{\BK}\left(C{\BM}^Z(X)\right) \ar[d]_{\wr} &{\BK}\left({\BM}^Z(X)\right) \ar[d]^{\wr}  \\
{\BK}\left(Ch^b\left(C{\BM}^Z(X)\right)\right) \ar[r]_{\sim}& {\BK}\left(Ch^b\left({\BM}^Z(X)\right)\right) \\
\enddiagram
$$
This establishes the equivalences (\ref{1918HomoEqui}). Combining with (\ref{25SptHus}) it follows that the sequence (\ref{190MainHSTH})
 is also a homotopy fibration. 
\pic $\eop$

\section{The ${\BG}W$-Theory Localization sequence}\label{SEcBGWSpSec}
In this section we develop exact sequence of ${\BG}W$-groups, similar to (\ref{190MainHSTH}, \ref{1911HomFFibr}). We refer to section \ref{SecVorAus} for notations, in particular $X$ denotes a quasi projective scheme over a noetherian 
affine scheme $\spec{A}$. We assume $1/2\in A$, in this section and whenever we talk about hermitian ($GW$) theory.
 Further refer to sections \ref{BISP4Exact}, \ref{subSecKTandGWSP}, \ref{dgCatSECGWsp}, \ref{subSecKaroBISP} for background on $GW$ theory, Karroubi
 ${\BG}W$ bispectrum $\BiSp$ and others. In particular, we distinguish ${\BG}W$ bispectrum ${\BG}W\LRf{{\SE}}^{\pm}\in \BiSp$ of an exact category ${\SE}$ with duality and 
 the bispectrum ${\BG}W^{[n]}\LRf{{\bf dg}{\SE}}\in \BiSp$ of the dg category 
 ${\bf dg}{\SE}$ of its bounded chain complexes. We will only discuss Karoubi ${\BG}W$ bispetrum, outside the background given in section \ref{BISP4Exact}. 
 

The  
${\BG}W$-theory differs from ${\BK}$-theory, in that that, for quasi projective schemes $X$ over a noetherian affine scheme $\spec{A}$,  the category ${\BM}^Z(X)$ does not have a natural duality. So, to do ${\BG}W$-theory (or $G{\CW}$, $\g{W}$ theory), our only option remains is to work with the category
$C{\BM}^Z(X)$ of perfect modules, which has a natural duality. Because of this our results in the section are sharper than the results on ${\BK}$-theory in section \ref{SECBKTHEY}.
\bD\label{dualDefExt}{\rm 
Suppose $X$ is a quasi projective scheme over $\spec{A}$. 
 Use other notations, as in (\ref{nota}, \ref{setUP}). 
 Let $Z\subseteq X$ be a closed subset and $grade(Z, X)=k$. 
  Fix  an invertible sheaf ${\SL}$ on $X$. Define the functor
 \begin{equation}\label{DualFunc}
^{\vee}:= \left(-\right)^{\vee}_{\SL}: C{\BM}^Z(X) \iso C{\BM}^Z(X)\quad {\rm by}\quad 
 {\SF}^{\vee}= {\SE}xt^k\left({\SF}, {\SL}\right)
 \end{equation} 
 Usually, we ignore to mention the label ${\SL}$. It is easy to check that $^{\vee}$ defines a duality.
 
}
\eD 


Following are some more notations 
and some comments.
\begin{notation}\label{79GWNoa}{\rm
Suppose $X$ is a quasi projective scheme over a noetherian affine scheme $\spec{A}$. 
 Assume $1/2\in A$.
 Let $Z\subseteq X$ be a closed subset and $grade(Z, X)=k$. We also fix an invertible sheaf 
 ${\SL}$ on $X$.
\bE
\item Recall the notations ${\bf P}erf(X)$, ${\bf P}erf_Z(X)$, ${\bf dg}{\SE}(X)$ and others from section \ref{SecVorAus}. We will be mindful to add prefix ${\bf dg}$ to an exact category, when we want to distinguish these two. 

%
%
%

\item Given $X$, as above, we fix a minimal injective resolution ${\SI}_{\bullet}$ of ${\CO}_X$:
\begin{equation}\label{injRes113}
\diagram 
0\ar[r] & {\CO}_X\ar[r] & {\SI}_0\ar[r] & {\SI}_{-1} \ar[r] & {\SI}_{-2} \ar[r]& \cdots\\
\enddiagram 
\end{equation}
Clearly, ${\SI}_{\bullet}\in {\SP}erf(X)$. 

Consequently, ${\SL} {\SI}_{\bullet}:={\SL}\otimes {\SI}_{\bullet}\in {\SP}erf(X)$, is an injective resolution of ${\SL}$. 

For ${\SF}_{\bullet}\in {\bf P}erf(X)$, we denote 
\begin{equation}\label{112DefSFveeD}
{\SF}^{\vee}:={\SF}^{\vee}_{\SL}={\CH}om\left({\SF}_{\bullet}, {\SL}{\SI}_{\bullet} \right).
\end{equation} 
\item ({\bf Lemma:}) It is a standard fact that %
the association ${\SF}\mapsto {\SF}^{\vee}$
endows $\left({\bf dg}{\bf P}erf(X), \q\right)$ with a structure of a dg category with duality and weak equivalences. 
Here $\q$ denotes the subcategory of quasi isomorphisms in  ${\bf P}erf(X)$.
For a proof see  \cite[Lem 4.8]{S17}.
\eE 
}
\end{notation}

Now we formulate a duality version of the main equivalence Theorem \ref{EquiALL}, as follows.

\bP
\label{d079verEqui4p6}{\rm 
Suppose $X$ is a quasi projective scheme over $\spec{A}$, with $1/2\in A$. 
 Use other notations, as in (\ref{nota}, \ref{setUP}). 
 Let $Z\subseteq X$ be a complete intersection closed subset, as in  (\ref{CIsetUp}), and $grade(Z, X)=k$. We also fix an invertible sheaf 
 ${\SL}$ on $X$.
 Then there is a natural duality preserving equivalence
 \begin{equation}\label{2186TkequDul}
 \diagram 
 {\bfD}^b\left(C{\BM}^Z(X)\right) \ar[r] & T^k{\SD}_Z^b\left({\SV}(X) \right)\\
 \enddiagram
 \end{equation}
 of derived categories. Here $T$ denotes the shift functor, on either side. The duality on the left side is induced by the duality given in Definition \ref{dualDefExt}. The duality $\#$ on right side is also induced by that on ${\SV}(X)$, namely $\#:=T^k{\CH}om\left(-, {\SL}\right)$. 
 }
 \eP
 \pf As was commented in Theorem
 \ref{EquiALL}, Equation \ref{longEQUI74}, there is a 
 natural equivalence
 $$
\omega:{\bfD}^b\left(C{\BM}^Z(X)\right) \lra {\SD}^b_Z\left({\SV}(X)\right)
 $$
 of derived categories, given by resolutions. For ${\SF}\in C{\BM}^Z(X)$, we have $\PDV{{\SF}}=k$.
 It follows from this that $\omega$ is actually an equivalence, as required (\ref{2186TkequDul}).
 \pic $\eop$
 
 \vspace{3mm}

Readers are referred to \cite[\S 1.10, pp. 1747]{S17}, \cite[Def 10.2.4, 10.2.15]{M23} for the jargon on dg categories with duality and weak equivalences. As usual, in chain complex categories,  the quasi isomorphisms $\q$ are the weak equivalences, by default.

\bL
\label{withDg145}{\rm
Suppose $X$ is a quasi projective scheme over $\spec{A}$ with $1/2\in A$. 
 Use other notations, as in (\ref{nota}, \ref{setUP}). 
 Let $Z\subseteq X$ be a closed subset and $grade(Z, X)=k$. We also fix an invertible sheaf 
 ${\SL}$ on $X$.  Then
 \bE
 \item\label{161One} The inclusion  $\zeta: {\bf dg}_Z\left({\SV}(X) \right)\lra {\bf dg}\left({\bf P}erf_Z(X)\right)$
 is a duality preserving morphism 
   of pointed dg categories with weak equivalences and duality. Here the duality on ${\SV}(X)$ is given by ${\SF}\mapsto {\CH}om\left({\SF}, {\SL}\right)$, and duality on right hand side is as in (\ref{112DefSFveeD}). 
 
 Further, the   functor
 $$
 \T\left({\bf dg}_Z\left({\SV}(X)\right) \right)\lra \T\left({\bf dg}\left({\bf P}erf_Z(X)\right)\right)
 $$
 of the associated triangulated categories, is a duality preserving  equivalence. 
 
 \item\label{216Two} The inclusion $\beta: {\bf dg}\left(C{\BM}^Z(X)\right) \lra T^k\left({\bf dg}\left({\bf P}erf_Z(X) \right)\right)$
 is a duality preserving morphism, of pointed dg categories with weak equivalences and duality, where $T$ denotes the shift. The duality on  ${\bf dg}\left(C{\BM}^Z(X)\right)$ is induced by (\ref{DualFunc}).

 Further,  assume that $Z$ is complete intersection, as in (\ref{CIsetUp}). Then the  functor
 $$
 \T\left({\bf dg}\left(C{\BM}^Z(X)\right) \right)
 \lra T^k\T\left({\bf dg}\left({\bf P}erf_Z(X)\right)\right)
 $$
of the associated triangulated categories is a duality preserving equivalence. 
 
 \eE 
 
 }
 \eL
 \pf First, for ${\SF}_{\bullet}\in Ch^b({\SV})$ the map 
 ${\CH}om\left({\SF}_{\bullet}, {\SL} \right)\iso {\CH}om\left({\SF}_{\bullet}, {\SL}{\SI}_{\bullet} \right)$
 is a quasi isomorphism. %
 By definition, $\zeta$ preserves duality. 
 Now, 
 we have the associated triangulated categories
 $$
\left\{ \begin{array}{l}
  \T\left({\bf dg}_Z\left({\SV}(X)\right) \right)={\SD}_Z^b({\SV}(X)\\
  \T\left({\bf dg}\left({\bf P}erf_Z(X)\right)\right)={\SD}_Z\left({\bf P}erf(X)\right)\\
  \end{array}\right. 
 $$
 In fact, the natural inclusion functor ${\bfD}\left({\SV}(X)\right) \iso {\bfD}\left({\bfP}erf(X)\right)$ is an equivalence, by \cite[Lem 3.8]{TT90}, \cite[Prop 3.4.8]{S11}. 
 Restricting this equivalence, to the complexes with support on $Z$, we obtain that
  $$
 {\SD}_Z\left({\SV}(X)\right)\lra  {\SD}_Z\left({\bfP}erf(X)\right)
 $$
 is an equivalence. This completes the proof of (\ref{161One}). 
 
 Now we proceed to prove (\ref{216Two}). We use the generic notation $^{\vee}$ to denote the dualities,
 for all categories. First, we prove that $\beta$ preserves duality.
    For ${\SF}_{\bullet}\in Ch^b\left(C{\BM}^Z(X)\right)$, its dual is given as follows:
  $$
  {\SF}_{\bullet}^{\vee}:\quad 
\diagram 
&&\deg=n\ar@{..>}[d]&&\\
\cdots\ar[r] &  {\CE}xt^k\left({\SF}_{-n-1}, {\SL} \right) \ar[r] &  {\CE}xt^k\left({\SF}_{-n}, {\SL} \right) \ar[r]  &{\CE}xt^k\left({\SF}_{-n+1}, {\SL} \right) \ar[r] & \cdots\\
 \enddiagram
  $$
 We are required  to prove that 
 $$
{\SF}^{\vee}_{\bullet} \iso {\CH}om\left({\SF}_{\bullet}, {\SL}{\SI}_{\bullet}\right)
 {\rm is ~a ~quasi~ isomorphism}
 $$
Since ${\SF}_i\in C{\BM}^Z(X)\subseteq C{\BM}^k(X)$, we have $k=grade\left({\SF}_i\right)=\PDV{{\SF}_i}$. Therefore 
 $$
 {\SE}xt^r\left({\SF}_i, {\SL}\right) =0 \qquad \forall i \quad {\rm and}\quad r\neq k
 $$
 %
 %
With this summary, for ${\SF}_{\bullet}\in Ch^b\left(C{\BM}^Z(X)\right)$, we obtain the double complex:
$$
\diagram 
&&\deg=-n\ar@{..>}[d]&&\\
&0\ar[d]&0\ar[d]&0\ar[d]&\\
\cdots\ar[r] &  {\CE}xt^k\left({\SF}_{n-1}, {\SL} \right) \ar[r]\ar[d] &  {\CE}xt^k\left({\SF}_{n}, {\SL} \right) \ar[r] \ar[d] &{\CE}xt^k\left({\SF}_{n+1}, {\SL} \right) \ar[r] \ar[d]& \cdots\\
\cdots\ar[r] & {\CH}om\left({\SF}_{n-1}, {\SL}{\SI}_{-k} \right) \ar[r] \ar[d]&  {\CH}om\left({\SF}_{n}, {\SL}{\SI}_{-k} \right) \ar[r] \ar[d] &{\CH}om\left({\SF}_{n+1}, {\SL}{\SI}_{-k} \right) \ar[r] \ar[d]& \cdots\\
\cdots\ar[r] & {\CH}om\left({\SF}_{n-1}, {\SL}{\SI}_{-(k+1)} \right) \ar[r] \ar[d]&  {\CH}om\left({\SF}_{n}, {\SL}{\SI}_{-(k+1)} \right) \ar[r] \ar[d] &{\CH}om\left({\SF}_{n+1}, {\SL}{\SI}_{-(k+1)} \right) \ar[r] \ar[d]& \cdots\\
\cdots\ar[r] & {\CH}om\left({\SF}_{n-1}, {\SL}{\SI}_{-(k+2)} \right) \ar[r] \ar[d]&  {\CH}om\left({\SF}_{n}, {\SL}{\SI}_{-(k+2)} \right) \ar[r] \ar[d] &{\CH}om\left({\SF}_{n+1}, {\SL}{\SI}_{-(k+2)} \right) \ar[r] \ar[d]& \cdots\\
&\cdots&\cdots &\cdots &\\
 \enddiagram
$$ 
Again, since ${\SF}_n\in C{\BM}^Z(X)$ and $k=grade(Z, X)$,  all the vertical lines above are exact. 
From this double complex, we obtain a natural transformation:
$$
\diagram
\beta\left({\SF}_{\bullet}^{\vee}\right) \ar[r] & \beta\left({\SF}_{\bullet}\right)^{\vee}=:
{\CH}om\left({\SF}_{\bullet}, {\SI}_{\bullet} \right) \\
\enddiagram 
$$
which is a quasi isomorphism. This establishes that $\beta$ preserves duality. 
Further, the associated triangulated categories
$$
\left\{\begin{array}{l}
\T\left({\bf dg}\left(C{\BM}^Z(X)\right) \right)= {\bfD}^b\left(C{\BM}^Z(X)\right)\\
\T\left({\bf dg}\left({\bf P}erf_Z(X)\right)\right)=T^k {\bfD}_Z\left({\bf P}erf(X)\right)\\
\end{array}\right.
$$
Now we have a commutative diagram of functors:
$$
\diagram 
 {\bfD}^b\left(C{\BM}^Z(X)\right) \ar[r]^{\sim} \ar@/_/@{-->}[rd]_{\sim}& T^k{\SD}^b_Z\left({\SV}(X) \right)\ar[d]^{\wr}\\
 & T^k {\bfD}_Z\left({\bf P}erf(X)\right)\\
\enddiagram 
$$
When $Z$ is complete intersection the horizontal arrow is an equivalence by (\ref{EquiALL}). The vertical arrow is also an equivalence by (\ref{161One}). 
This completes the proof of (\ref{216Two}). $\eop$

 \subsection{Homotopy fibrations of ${\BG}W$ Theory bi-spectra ${\BiSp}$}
 
 Now we state the equivalence theorem of ${\BG}W$-theory bi-Spectra.
 Refer to section \ref{BISP4Exact}, \ref{SecVorAus} for notations and definitions on $GW$-theory, spectra and bispectra.
For a dg category ${\SA}=\LRf{{\SA}, \q, ^{\vee}, \varpi}$ with duality and weak equivalences, the ${\BG}W$ bispectrum
 ${\BG}W\LRf{{\SA}}\in\BiSp$ takes values in 
 the category $\BiSp$  of bi-spectra (see \S \ref{subSecKaroBISP}). 
 The ${\BG}W$-groups  
 ${\BG}W_p^{[n]}\LRf{{\bf dg}{\SE}}$ are defined as the homotopy groups $\forall n, p\in {\BZ}$. 
 Similar notations are used for exact categories with duality and weak equivalences
 (see \S \ref{BISP4Exact}). 
\bT\label{BiSpectrEquivho}{\rm 
Suppose $X$ is a quasi projective scheme over $\spec{A}$. 
 Use other notations, as in (\ref{nota}, \ref{setUP}). 
 Let $Z\subseteq X$ be a complete intersection closed subset, as in (\ref{CIsetUp}),
  and $grade(Z, X)=k$. 
  We also fix an invertible sheaf 
 ${\SL}$ on $X$.  We assume $1/2\in A$. Then the following zig-zag sequence 
 $$
 \diagram
 {\BG}W^{[r]}\left({\bf dg}C{\BM}^Z(X) \right) \ar[r]_{\sim\qquad}
 &  {\BG}W^{[r+k]}\left({\bf dg}\left({\bf P}erf_Z(X)\right) \right)
 &  {\BG}W^{[r+k]}\left({\bf dg}_Z{\SV}(X)\right)\ar[l]^{\qquad\sim}  \\
 \enddiagram 
 $$
 are stable homotopy equivalences 
   in the category $\mathrm{BiSp}$ of bispectra, $\forall r\in {\BZ}$. 
%
 Consequently, the ${\BG}W$-theory groups 
 $$
 {\BG}W^{[r]}_p\left({\bf dg}C{\BM}^Z(X)\right) \cong  {\BG}W_p^{[r+k]}\left({\bf dg}_Z{\SV}(X)\right) \quad {\rm are~isomorphic}~\forall ~r, p \in {\BZ}.
 $$
}
\eT
\pf Follows directly from Theorem \ref{withDg145}, and \cite[Thm 8.9]{S17},
\cite[pp 466]{M23}. \pic $\eop$

\vspace{3mm}
Now we are ready to state the theorem on  homotopy fibration of ${\BG}W$-bispectra.

\bT\label{BiSpectrEHmoFib}{\rm 
Suppose $X$ is a quasi projective scheme over $\spec{A}$. 
 Use other notations, as in (\ref{nota}, \ref{setUP}). 
 Let $Z\subseteq X$ be a closed complete intersection subset, as in (\ref{CIsetUp}),  with $grade(Z, X)=k$ and let $U=X-Z$. We also fix an invertible sheaf 
 ${\SL}$ on $X$.  We assume $1/2\in A$. 
 Then the sequence
 \begin{equation}\label{190MainBosTH}
 \diagram 
  {\BG}W^{[r]}\left({\bf dg}C{\BM}^Z(X)\right) \ar[r] &  {\BG}W^{[r+k]}\left({\bf dg}{\SV}(X)\right)  \ar[r]  &  
  {\BG}W^{[r+k]}\left({\bf dg}{\SV}(U)\right)  \\
 \enddiagram
 \end{equation}
 is a homotopy fibration in the category  $\mathrm{BiSp}$ of  bispectra $\forall r\in {\BZ}$. 
 The first map here is the zig-zag map given in (\ref{BiSpectrEquivho}). 
 Consequently, there is a long exact sequence 
  \begin{equation}\label{1911HomBOSFFibr}
 \diagram
 \cdots \ar[r] & {\BG}W_p^{[r]}\left({\bf dg}\left(C{\BM}^Z(X)\right) \right)  \ar[r] & {\BG}W_p^{[r+k]}\left({\bf dg}{\SV}(X) \right)  \ar[r] 
 & {\BG}W_p^{[r+k]}\left({\bf dg}{\SV}(U) \right) \\
 \ar[r] & {\BG}W_{p-1}^{[r]}\left({\bf dg}\left(C{\BM}^Z(X)\right) \right)\ar[r] &\cdots&\forall p, r\in {\BZ}\\
 \enddiagram
 \end{equation} 
 }
 \eT
 \pf Consider the sequence of dg categories with dualities and weak equivalences (quasi isomorphism) 
 $$
 \diagram 
 {\bf dg}\left(C{\BM}^Z(X)\right) \ar[r] & T^k{\bf dg}{\SV}(X)\ar[r] & T^k{\bf dg}{\SV}(U)\\
 \enddiagram
 $$
 The corresponding sequence of associated triangulated categories is as follows:
 $$
 \diagram 
 {\bfD}^b\left(C{\BM}^Z(X)\right) \ar[r] & T^k{\bfD}^b\left({\SV}(X)\right)\ar[r] & T^k{\bfD}^b\left({\SV}(U)\right)\\
 \enddiagram
 $$
 This sequence is exact up to factor (direct summand), by \cite[Prop A.4.7]{S11},
 because ${\bfD}^b\left(C{\BM}^Z\LRf{{\SV}(X)}\right)  \iso {\SD}^b_Z({\SV}(X))$ is an equivalence by 
 (\ref{EquiALL}). 
 Now the proof follows from 
 \cite[Theorem 8.10]{S17} that the sequence (\ref{190MainBosTH}) is a homotopy fibration
 in $\BiSp$. 
 \pic $\eop$

\begin{remark}{\rm 
Let $A$ be a commutative ring and ${\SE}:=\LRf{{\SE}, ^\vee, \varpi}$ be a $A$-linear exact category with duality. Assume $1/2\in A$. 
For the purpose of making statements on ${\BG}W$-bispectrum of ${\SE}$, we 
defined ${\BG}W({\SE})\in \BiSp$ directly, in section \ref{BISP4Exact}, without any reference to the definition of
${\BG}W({\bf dg}{\SE})\in \BiSp$ as a dg category. Our sense of esthetics and completeness of the theory demands this. Refer to  section \ref{BISP4Exact} for further discussions. We point out that  two ${\BG}W$ bispectra 
$$
\left\{\begin{array}{l}
{\BG}W\LRf{{\SE}}:={\BG}W\LRf{{\SE}}^+:={\BG}W\LRf{{\SE}, ^{\vee}, \varpi}\\
{\BG}W\LRf{{\SE}}^-:={\BG}W\LRf{{\SE}, ^{\vee}, -\varpi}\\
\end{array}\right. \qquad \in {\BiSp}
$$
 are defined, while 
${\BG}W^{[n]}\LRf{{\bf dg}{\SE}}\in \BiSp$ has four periodicity. 
}
\end{remark} 
%
Now, we interpret Theorem \ref{BiSpectrEquivho}  the above in terms of ${\BG}W$-bispetrum of exact categories. 
\bT\label{BiSpeEXACCuivho}{\rm 
Suppose $X$ is a quasi projective scheme over $\spec{A}$. 
 Use other notations, as in (\ref{nota}, \ref{setUP}). 
 Let $Z\subseteq X$ be a 
 closed subset, 
 and $grade(Z, X)=k$. We also fix an invertible sheaf 
 ${\SL}$ on $X$.  We assume $1/2\in A$. Then the following 
 $$
 \left\{\begin{array}{l}
  {\BG}W\left(C{\BM}^Z(X)\right) \cong  {\BG}W^{[k]}\left({\bf dg}_Z{\SV}(X)\right) \\
    {\BG}W\left(C{\BM}^Z(X)\right)^- \cong  {\BG}W^{[k+2]}\left({\bf dg}_Z{\SV}(X)\right) \\
 \end{array}\right.
 \quad {\rm are ~stable~equivalences~in}~\BiSp
 $$
of  zig-zag maps. 
}
\eT
\pf Follows from (\ref{BiSpectrEquivho}), with $r=0, 2$. \pic $\eop$

\vspace{3mm}
The following is reinterpretation of (\ref{BiSpectrEHmoFib}), for even shifts.

\bT\label{BiExacttrEHmoFib}{\rm 
Suppose $X$ is a quasi projective scheme over $\spec{A}$. 
 Use other notations, as in (\ref{nota}, \ref{setUP}). 
 Let $Z\subseteq X$ be a complete intersection closed subset,
 as in  (\ref{CIsetUp}),   $grade(Z, X)=k$.  Write  $U=X-Z$. We also fix an invertible sheaf 
 ${\SL}$ on $X$.  We assume $1/2\in A$. 
 Then the sequences
  \begin{equation}\label{extMainduiNomBos}
 \left\{\begin{array}{ll}
  \diagram 
  {\BG}W\left(C{\BM}^Z(X)\right) \ar[r] &  {\BG}W\left({\SV}(X)\right)  \ar[r]  &  
  {\BG}W\left({\SV}(U)\right)  \\
 \enddiagram & in~case~ k=0~mod~4\\
  \diagram 
  {\BG}W\left(C{\BM}^Z(X)\right) \ar[r] &  {\BG}W\left({\SV}(X)\right)^-  \ar[r]  &  
  {\BG}W\left({\SV}(U)\right)^-  \\
 \enddiagram &  in~case~ k=2~mod~4\\
   \diagram 
  {\BG}W\left(C{\BM}^Z(X)\right)^- \ar[r] &  {\BG}W\left({\SV}(X)\right)^-  \ar[r]  &  
  {\BG}W\left({\SV}(U)\right)^-  \\
 \enddiagram &  in~case~ k=0~mod~4\\
  \diagram 
  {\BG}W\left(C{\BM}^Z(X)\right)^- \ar[r] &  {\BG}W\left({\SV}(X)\right)  \ar[r]  &  
  {\BG}W\left({\SV}(U)\right)  \\
 \enddiagram  &  in~case~ k=2~mod~4\\
 \end{array}\right.
 \end{equation}
 are  homotopy fibratiion in the category  $\mathrm{BiSp}$ of bispectra. 
 Consequently, the long exact sequences of the ${\BG}W$ groups, corresponding to (\ref{extMainduiNomBos}) can be written down.
  }
 \eT
 \pf Follows from  (\ref{BiSpectrEHmoFib}), by taking $r=0, 2$. See Definition Remark \ref{GW444Ed} 
 and equation (\ref{DefbhseExc}). \pic $\eop$ 

%

\section{Applications} \label{SecAPPP}
Localization theorems, such as Theorems \ref{Equi1909ALL},  \ref{BiSpectrEquivho}, 
\ref{BiExacttrEHmoFib}, usually would have some  consequences. In this section, we list some of them. 
Usually, we consider only idempotent complete the exact categories ${\SE}$ and exact categories with duality 
$\LRf{{\SE}, ^{\vee}, \varpi}$. 
For such  idempotent complete exact categories ${\SE}$, we use the notations:
$$
\left\{\begin{array}{ll}
{\BK}_i\LRf{{\SE}}=: {K}_n\LRf{{\SE}} & \forall i\geq 0\\
{\BG}W^{[n]}_i\LRf{{\bf dg}{\SE}, \q, ^{\vee}, \varpi}=:{G}W^{[n]}_i\LRf{{\bf dg}{\SE}, \q^{\vee}, \varpi}%
& \forall i\geq 0, n\in {\BZ}\\
\end{array}\right.
$$
where $\q$ denotes the subcategory of quasi isomorphisms.
First, we state the application to the punctured schemes.
\bT\label{puncute}{\rm 
Let $(A, \m)$ be a Cohen Macaulay local ring and $X=\spec{A}$, 
 $Z=V(\m)$  and $U=X-Z$. Assume $1/2\in A$ and write $d=\dim A$.
Then the sequence 
$$
\diagram
{\BK}\LRf{C{\BM}^Z(X)} \ar[r] & {\BK}\LRf{{\SV}(X)} \ar[r]& {\BK}\LRf{{\SV}(U)}\\
\enddiagram\qquad in\qquad \Sp
$$
is a homotopy fibration. Consequently, there is a long exact sequence of ${\BK}$-groups.
Further, %
$$
{\BK}_r\LRf{C{\BM}^Z(X)}=
\left\{\begin{array}{ll}
 K_r\LRf{C{\BM}^Z(X)}&\forall  r\geq 0\\
 0 &otherwise\\
\end{array}\right.
$$
Consequently, 
$$
\forall r\leq -1 \quad {\rm the~map}\quad {\BK}_r\LRf{{\SV(X)}}\iso {\BK}_r\LRf{{\SV(U)}}
\quad {\rm is ~an~isomorphism.}
$$
In particular, if $U$ is regular (i.e. $A$ is regular in codimension $d-1$) then 
$$
 {\BK}_r\LRf{{\SV(X)}} =0 \qquad \qquad \forall r\leq -1
 $$
}
\eT
\pf Let $d=\dim(A)$. Note $Z=V(\f_1, \ldots, f_d)$ is complete intersection. Now, the first statement is an immediate consequence of the localization theorem \ref{Equi1909ALL} of ${\BK}$-theory spectra. 
 Note, $C{\BM}^Z(X)$ represents the
category of  $A$-modules, with finite length and finite projective dimension. 
Further note that 
$C{\BM}^Z(X)$ is a noetherian abelian category. Now it follows from \cite[Thm 7]{S06} that 
$$
\forall r\leq -1 \quad {\BK}_r\LRf{C{\BM}^Z(X)}=0
$$
Now, it follows from the long exact sequence that ${\BK}_r\LRf{{\SV}(X)}
\cong {\BK}_r\LRf{{\SV}(U)}~\forall r\leq -1$. If $U$ is regular then 
${\BK}_r\LRf{{\SV}(U)}=0~\forall r\leq -1$ \cite[pp 182]{S11}. So, the last statement follows. 
\pic $\eop$

\vspace{3mm}
The above (\ref{puncute}) can be improved, as follows. 
\bT\label{nFutos}{\rm 
Let $X$ be a quasi projective scheme over a noetherian affine scheme $\spec{A}$.
Let $Z=\LRs{\m_1, \m_2, \ldots, \m_{\ell}}\subseteq X$ be finite set of closed points.
Further, assume   that $Z$ is complete intersection (\ref{CIsetUp}) and $U=X-Z$. 
Then the sequence 
$$
\diagram
{\BK}\LRf{C{\BM}^Z(X)} \ar[r] & {\BK}\LRf{{\SV}(X)} \ar[r]& {\BK}\LRf{{\SV}(U)}\\
\enddiagram\qquad in\qquad \Sp
$$
is a homotopy fibration. Consequently, there is a long exact sequence of ${\BK}$-groups.
Further, %
\begin{equation}\label{SplitGHJK}
{\BK}_r\LRf{C{\BM}^Z(X)}=
\left\{\begin{array}{ll}
 \bigoplus_{i=1}^{\ell} K_r\LRf{C{\BM}^{V(\m_i)}\LRf{\spec{{\CO}_{X, \m_i}}}}
 &\forall  r\geq 0\\
 0 &otherwise\\
\end{array}\right.
\end{equation} 
Consequently, 
$$
\forall r\leq -1 \quad {\rm the~map}\quad {\BK}_r\LRf{{\SV(X)}}\iso {\BK}_r\LRf{{\SV(U)}}
\quad {\rm is ~an~isomorphism.}
$$
In particular, if $U$ is regular,   
then 
$$
 {\BK}_r\LRf{{\SV(X)}} =0 \qquad \qquad \forall r\leq -1
 $$
 
}
\eT 
\pf  Note that $C{\BM}^Z(X)$ is a noetherian abelian category. Now, there is an equivalence $\Psi: C{\BM}^Z(X) \iso \prod_i C{\BM}^{V(\m_i)}(X)$. 
 To define $\Psi$ let $M\in C{\BM}^Z(X)$.  For $i=1, \ldots, \ell$ define $M_i\in {\BM}^{V(\m_i)}(X)$, as follows
 \begin{equation}\label{PSIfun26}
 \left\{\begin{array}{l}
 (M_i)_{|X-\{\m_1, \ldots \hat{\m}_i\ldots, \m_{\ell}\}}=M_{|X-\{\m_1, \ldots \hat{\m}_i\ldots, \m_{\ell}\}}\\
 (M_i)_{|X-\{\m_i\}}=0\\
 \end{array}\right.
 \quad Define\quad \Psi(M) =\LRf{M_1, \ldots, M_{\ell}}
 \end{equation} 
The inverse functor is defined by direct sum.
Finally, there are equivalences 
 \begin{equation}\label{PSIfDusoun26}
 \Phi: C{\BM}^{V(\m_i)}\LRf{X} \iso C{\BM}^{V(\m_i)}\LRf{\spec{{\CO}_{X, \m_i}}}
 \quad M\mapsto M_{\m_i}
  \end{equation} 
Now the proof is completed, as in  the proof of (\ref{puncute}). 
\pic $\eop$

\begin{remark}{\rm 
If $\m:=\m_i$ is a regular point, then the equation can be further simplified, as follows. Suppose 
$(A, \m, \kappa)$ is a regular local ring. Then  by D\'{e}vissage theorem 
\cite[pp 111]{M23}
$$
{\BK}_rC{\BM}^{V(\m)}\LRf{\spec{A}}
= \left\{\begin{array}{ll}
K_r\LRf{\kappa} & r\geq 0\\
0 & r\leq -1\\
\end{array}\right.
$$
}
\end{remark}

We leave it to the reader to write
a ${\BG}W$-version of (\ref{puncute}). However, negative ${\BG}W$ groups noetherian abelian categories do not necessarily vanish, as in ${\BK}$-theory. 
 We write down the ${\BG}W$-version of (\ref{nFutos}), as follows.

\bT\label{BGWnFutos}{\rm 
Let $X$ be a quasi projective scheme over a noetherian affine scheme $\spec{A}$,
with $1/2\in A$.
Let $Z=\LRs{\m_1, \m_2, \ldots, \m_{\ell}}\subseteq X$ be finite set of closed points.
Further, assume   that $Z$ is complete intersection (\ref{CIsetUp}), with 
$d=grade(Z, X)$ and $U=X-Z$. 
Then $\forall r\in {\BZ}$, the sequence 
\begin{equation}\label{2737dgKiExct}
\diagram
{\BG}W^{[r]}\LRf{{\bf dg}C{\BM}^Z(X)} \ar[r] & {\BG}W^{[r+d]}\LRf{{\bf dg}{\SV}(X)} \ar[r]& {\BG}W^{[r+d]}\LRf{{\bf dg}{\SV}(U)}\\
\enddiagram~~ in~~ \BiSp
\end{equation}
is a homotopy fibration in $\BiSp$. %
Further, %
\begin{equation}\label{SplitBIGHJK}
{\BG}W^{[r]}\LRf{{\bf dg}C{\BM}^Z(X)} \iso \prod_{i=0}^{\ell}
{\BG}W^{[r]}\LRf{{\bf dg}C{\BM}^{V(\m_i)}\LRf{X_{\m_i}}}
\end{equation} 
is a stable equivalence in $\BiSp$, where $X_{\m_i}=\spec{{\CO}_{X, \m_i}}$. 

Consequently, the 
sequence (\ref{2737dgKiExct}), leads to  long exact sequences of
 ${\BG}W$-groups.
Further, one can reinterpret (\ref{2737dgKiExct}) in terms of ${\BG}W$-spectra of the exact categories, as in (\ref{extMainduiNomBos}), which we skip.
 
}
\eT 
\pf Same as the proof of (\ref{nFutos}), while we use the theorem \ref{BiSpectrEHmoFib}.
Also note that the functors $\Psi$ and $\Phi$  (\ref{PSIfun26}, \ref{PSIfDusoun26})
are duality preserving equivalences. 
\pic $\eop$

\vspace{3mm} 
When the points $\m_i$ in (\ref{BGWnFutos}) are regular, the D\'{e}vissage theorems
 for
${\BG}W$-groups or bispectrum applies \cite{M24}. Thus the groups 
${\BG}W^{[r]}_n\LRf{C{\BM}^{V(\m_i)}\LRf{\spec{{\CO}_{X, \m_i}}}}$ have better descriptions, as follows. 

\vspace{3mm}
%
 \bT\label{DivCMmXpcmi}{\rm 
 Let $\LRf{R, \m, \kappa}$ be a noetherian regular local ring, $X=\spec{R}$, $Z=V(\m)$, $d=\dim R$. Note
 $C{\BM}^{Z}\LRf{X}$ is the category of  $R$-modules $M$ with finite length 
 ({\it and $\PDV{M}=d$}). Then 
 $C{\BM}^{Z}\LRf{X}$  has a natural duality $M^{\vee}=Ext^d\LRf{M, R}$. 
 Let $\varpi$ denote the double dual identification in $C{\BM}^{Z}\LRf{X}$.
 As usual, $\SV(\kappa)$ denotes the category of finite dimensional $\kappa$-vector spaces.  Note $\SV(\kappa)$ also has a natural duality, defined by $M^{\vee}=Hom\LRf{M, \kappa}$. Let $\e$ denote the double dual identification in ${\SV}(\kappa)$,
 given by evaluation. Then there is a natural equivalence, of bi spectra:
 $$
 {\BG}W^{[n]}\LRf{{\bf dg}C{\BM}^Z(X),  \q,  ^{\vee},\varpi} \cong   {\BG}W^{[n]}\LRf{{\bf dg}{\SV}(\kappa), \q, ^{\vee},  \e} \qquad \forall n\in {\BZ}.
 $$
 Here, on left side, $\varpi$ 
 continues to denote the double dual identification induced by $\varpi$, 
 and $\q$ denotes quasi isomorphism. 
 We write
 $$
 \left\{\begin{array}{ll}
 W(\kappa)=W\LRf{{\SV}(\kappa), ^{\vee}, \e} & W(\kappa)^-=W\LRf{{\SV}(\kappa), ^{\vee}, -\e} \\
  GW_p(\kappa)=GW_p\LRf{{\SV}(\kappa), ^{\vee}, \e} & GW_p(\kappa)^-=GW\LRf{{\SV}(\kappa), ^{\vee}, -\e} \\
 \end{array}\right.
 $$
 %
%
  %
Consequently,
 $$
  {\BG}{W}^{[n]}_p\LRf{{\bf dg}C{\BM}^Z(X), ^{\vee},  \q, \varpi} \cong 
  \left\{\begin{array}{ll}
 %
  W(\kappa)&  p\leq -1, n-p-d=0~mod~4\\ 
  0 &  p\leq -1, ~n-p-d=1, 2, 3~mod~4\\ 
   G{W}_p(\kappa)& n=0~mod ~4,~ p\geq 0\\
   G{W}_p(\kappa)^-& n=2~mod ~4, p\geq 0\\
    \end{array}\right.
 $$
 For $p\geq 0$ and $n=2m+1$  odd, we have an exact sequence (see \cite{M24}).

 }
 \eT
 \pf See \cite{M24}. $\eop$

 The following is an immediate corollary.
 \bC\label{2518CorexLocal}{\rm
 Consider the setup and hypotheses of (\ref{DivCMmXpcmi}).  Then 
 $$
  {\BG}{W}_p\LRf{C{\BM}^Z(X), ^{\vee}, \varpi} \cong 
  \left\{\begin{array}{ll}
  W(\kappa)&  p\leq -1, p+d=0~mod~4\\ 
  0 &  p\leq -1, ~p+d=1, 2, 3~mod~4\\ 
  GW_p(\kappa)&  p\geq 0\\
    \end{array}\right.
 $$
and 
 $$
  {\BG}{W}_p\LRf{C{\BM}^Z(X), ^{\vee}, \varpi}^- \cong 
  \left\{\begin{array}{ll}
 %
  W(\kappa)&  p\leq -1, p+d=2~mod~4\\ 
  0 &  p\leq -1, ~p+d=1, 2, 3~mod~4\\ 
  G{W}_p(\kappa)^-
   &  p\geq 0\\
    \end{array}\right.
 $$
 }
 \eC
 \pf Follows from (\ref{DivCMmXpcmi}) and section \ref{BISP4Exact}.
 $\eop$
 
 The following is the exact category version of (\ref{BGWnFutos}).
 \bC\label{exctnfuto}{\rm 
 Under the hypothesis of  (\ref{BGWnFutos}),  the following sequences:
 $$
\left\{ \begin{array}{ll}
 \diagram 
 {\BG}W\LRf{C{\BM}^Z(X)} \ar[r] & {\BG}W\LRf{{\SV}(X)} \ar[r]& {\BG}W\LRf{{\SV}(U)}
 \enddiagram & r=0, d= 0 ~mod~4\\
 \diagram 
 {\BG}W\LRf{C{\BM}^Z(X)} \ar[r] & {\BG}W\LRf{{\SV}(X)} ^-\ar[r]& {\BG}W\LRf{{\SV}(U)}^-
 \enddiagram & r=0, d= 2 ~mod~4\\
 \end{array}\right.
 $$
 are stable homotopy fibration in $\BiSp$.  
 
 Now assume that $\m_i$ are regular points,  and let $\kappa(\m_i)$ denote the residue field of $\m_i$.
 In the first case, $r=0$ and $d=0~mod~ 4$, it follows from the  exact sequence of ${\BG}W$-gropus that 
 \bE
 \item The long exact sequence of non negative groups terminates as follows, at degree zero:
 $$
  \diagram 
 \cdots  \ar[r]&
\bigoplus_{i=1}^{\ell} GW_0\LRf{\kappa(\m_i)}\ar[r] & GW_0\LRf{{\SV}(X)} \ar[r] & GW_0\LRf{{\SV}(U)}\ar[r] &0\\
 \enddiagram
 $$
 \item We have isomorphisms of groups
 $$
 \left\{\begin{array}{l}
 \diagram 
 GW_{-1-4n}\LRf{{\SV}(X)} \ar[r]^{\sim} & GW_{-1-4n}\LRf{{\SV}(U)}\\
 \enddiagram\\
  \diagram 
 GW_{-2-4n}\LRf{{\SV}(X)} \ar[r]^{\sim} & GW_{-2-4n}\LRf{{\SV}(U)}\\
 \enddiagram\\
 \end{array}\right.
 \qquad \qquad \forall n\geq 0
 $$ 
 \item There are exact sequences
 $$
  \diagram 
0\ar[r] & GW_{-3-4n}\LRf{{\SV}(X)} \ar[r] & GW_{-3-4n}\LRf{{\SV}(U)}\ar[r] &\\
\bigoplus_{i=1}^{\ell} W\LRf{\kappa(\m_i)}\ar[r] & GW_{-4-4n}\LRf{{\SV}(X)} \ar[r] & GW_{-4-4n}\LRf{{\SV}(U)}\ar[r] &0\\
 \enddiagram
 \quad \forall n\geq 0.
 $$
 \eE 
Similar analysis is possible in other case. 
 
 }
 \eC
 \pf Follows from (\ref{BGWnFutos}) and the descriptions in (\ref{2518CorexLocal}).
  $\eop$


\vspace{3mm}
The following is another example exploiting regularity.
\begin{example}\label{RegualrU}{\rm 
Let $X$ be a quasi projective scheme over a noetherian affine scheme $\spec{A}$. 
Let $U\subseteq X$ be an open subset of $X$, and $Z=X-U$ be complete intersection,
as in (\ref{CIsetUp}). Assume that $U$ is regular. Then 
$$
{\BK}_p\LRf{C{\BM}^Z(X)} \cong {\BK}_p\LRf{{\SV}(X)} \qquad \forall p\leq -2
$$
and the long exact sequence on the left terminates at the term ${\BK}_{-1}({\SV}(X))$.
}
\end{example} 
\pf By the  localization theorem \ref{Equi1909ALL}, the sequence 
$$
\diagram
{\BK}\LRf{C{\BM}^Z(X)} \ar[r] & {\BK}\LRf{{\SV}(X)} \ar[r]& {\BK}\LRf{{\SV}(U)}\\
\enddiagram\qquad in \qquad \Sp
$$
is a homotopy fibration, of ${\BK}$-theory spectra. 
By
 \cite[Thm 7]{S06},
we have  ${\BK}_k({\SV}(U)) ={\BK}_k(Coh(U)) = 0~\forall k\leq -1$ . The rest follows from the long exact sequence. $\eop$

\vspace{3mm}
We use Theorem \ref{BiExacttrEHmoFib} to obtain a ${\BG}W$ analogue of 
 example \ref{RegualrU}. 
 
\begin{example}\label{RegualrU}{\rm 
Let $X$ be a quasi projective scheme over a noetherian affine scheme $\spec{A}$. 
Let $U\subseteq X$ be an open subset of $X$, and $Z=X-U$ be complete intersection,
as in (\ref{CIsetUp}). Assume that $U$ is regular and  $grade(Z, X)=k$. Also, let ${\CL}$ be a fixed line bundle on $X$ and $1/2\in A$. Then 
corresponding to each of the four homotopy fibrations in (\ref{extMainduiNomBos}), we obtain an exact sequence of ${\BG}W$ groups. We only write down the exact sequence corresponding to the first homotopy fibration in (\ref{extMainduiNomBos}),
which applies when $k=0~mod~4$, as follows:
$$
\diagram
\cdots \ar[r] &{\bf GW}_0\LRf{C{\BM}^Z(X)}\ar[r] & {\bf GW}_0\LRf{{\SV}(X)} \ar[r] & {\bf GW}_0\LRf{{\SV}(U)} \\
 \ar[r] &{\BG}W_{-1}\LRf{C{\BM}^Z(X)}\ar[r] & {\BG}W_{-1}\LRf{{\SV}(X)} \ar[r] & 
 W^1\LRf{{\SV}(U)} \\
  \ar[r] &{\BG}W_{-2}\LRf{C{\BM}^Z(X)}\ar[r] & {\BG}W_{-2}\LRf{{\SV}(X)} \ar[r] & W\LRf{{\SV}(U)}^- \\
    \ar[r] &{\BG}W_{-3}\LRf{C{\BM}^Z(X)}\ar[r] & {\BG}W_{-3}\LRf{{\SV}(X)} \ar[r] & W^3\LRf{{\SV}(U)} \\
       \ar[r] &{\BG}W_{-4}\LRf{C{\BM}^Z(X)}\ar[r] & {\BG}W_{-4}\LRf{{\SV}(X)} \ar[r] &W\LRf{{\SV}(U)}\\
  \ar[r] & \cdots &&\\
\enddiagram
$$
}
\end{example} 
\pf We consider the case $k=4k_0$ only. We start with the homotopy fibration
$$
  \diagram 
  {\BG}W\left(C{\BM}^Z(X)\right) \ar[r] &  {\BG}W\left({\SV}(X)\right)  \ar[r]  &  
  {\BG}W\left({\SV}(U)\right)  \\
 \enddiagram
 $$
 All the non negative groups are the homotopy groups of the $GW$-spaces \cite[450, 466]{M23}. Note ${\BK}_p\LRf{{\SV}(U)}=0~\forall p\leq -1$. It follows from this that 
 \cite[pp 1805]{S17}
 $$
 \G{W}_p({\SV}(U)) \cong  G{\CW}_p({\SV}(U)) \cong {\BG}W_p({\SV}(U))~
 \qquad \qquad \forall p\leq -1
 $$
  Therefore
 $\forall p\leq -1$ we have 
 $$
 {\BG}W_p({\SV}(U))\cong {G}{\CW}_p({\bf dg}{\SV}(U))
 \cong W^{-p}({\SV}(U))
 =
 \left\{\begin{array}{ll}
 W({\SV}(U)) & if~p=0~mod~4\\
  W({\SV}(U))^- & if~p=2~mod~4\\
 \end{array}\right.
 $$
\pic $\eop$ 


\vspace{3mm}
Thomason \cite{TT90}   considered the ${\BK}$-theory of perfect complexes, and likewise Schlichting \cite{S11, S17} considered ${\BG}W$-theory of perfect complexes. We differ from \cite{TT90, S11, S17} in that that we grind down  the theory to the  level of exact subcategories
of $Coh(X)$ (see Example \ref{DelineExample}). For notations ${\bf P}erf(X)$ and ${\bf P}erf_Z(X)$ refer to (\ref{nota}). Notationally, we denote Thomason ${\BK}$ theory by ${\BK}\LRf{{\bf P}erf(X)}$ and likewise, to distinguish them from our theory on exact categories.
\bT[Agreement]\label{AgrBkt}{\rm 
Let $X$ be a quasi projective scheme over a noetherian affine scheme $\spec{A}$, $U\subseteq X$ be open and $Z=X-U$. Assume that $Z$ is complete intersection, as in 
(\ref{CIsetUp}). 
Then the natural degree zero inclusion $C{\BM}^Z(X) \lra {\bf P}erf_Z(X)$ induces, a homotopy equivalence: 
$$
{\BK}\LRf{C{\BM}^Z(X)}\iso {\BK}\LRf{{\bf P}erf_Z(X)} \qquad in \qquad \Sp. 
$$
}
\eT
\pf We have commutative diagram 
$$
\diagram 
{\BK}\LRf{C{\BM}^Z(X)} \ar[r]\ar[d] & {\BK}\LRf{{\SV}(X)}\ar[r] \ar[d]&  {\BK}\LRf{{\SV}(U)}\ar[d]\\
{\BK}\LRf{{\bf P}erf_Z(X)} \ar[r] & {\BK}\LRf{{\bf P}erf(X)}\ar[r] &  {\BK}\LRf{{\bf P}erf(U)}\\
\enddiagram
$$
of Homotopy fibration. For the lower fibration see \cite{TT90, S11} or \cite[pp 411]{M23}, and the upper fibration is an application of (\ref{Equi1909ALL}). It is known that the $2^{nd}$ and $3^{rd}$ vertical maps are homotopy equivalence (see \cite{TT90, S11} or \cite[pp 410]{M23}. Now one can write down the commutative diagram of the homotopy groups and use Five Lemma. \pic $\eop$

%

\vspace{3mm} 
Schlichting \cite[pp 1804]{S17} considers ${\BG}W$-theory of ${\bf dg}{\SV}(X)$ ({\it  denoted by} $sPerf(X)$ {\rm in} \cite{S17}). We  have a ${\BG}W$-bispectra analog to (\ref{AgrBkt}). 

\bT[Agreement]\label{AgrBktBGW}{\rm 
Let $X$ be a quasi projective scheme over a noetherian affine scheme $\spec{A}$, $U\subseteq X$ be open and $Z=X-U$. Assume $1/2\in A$ and $Z$ is complete intersection, as in 
(\ref{CIsetUp}), with $k=grade(Z, X)$. Let $n\in {\BZ}$,  and ${\CL}$ be a line bundle on $X$.  
Then   the natural map
$$
{\BG}W^{[n-k]}\LRf{{\bf dg}C{\BM}^Z(X), \q, ^{\vee}, \varpi}\iso {\BG}W^{[n]}\LRf{{\bf dg}_Z{\SV}(X), \q, ^{\vee {\CL}}, \varpi} \qquad in \qquad \BiSp
$$
is a stable homotopy equivalence in the category $\BiSp$ of bispectra.  Here, on 
${\bf dg}C{\BM}^Z(X)$, the duality $^{\vee}$ is induced by $M\mapsto {\CE}xt^k(M, {\CL})$ and on ${\bf dg}{\SV}(X)$ the duality is induced by $P\mapsto {\CH}om(M, {\CL})$. Further,  on both sides, $\varpi$ denotes the double dual identifications, and $\q$ denotes quasi isomorphisms. 
Further, there are equivalences in $\BiSp$, as follows 
$$
{\BG}W^{[n]}\LRf{{\bf dg}_Z{\SV}(X)} 
 \cong 
\left\{\begin{array}{ll}
{\BG}W\LRf{C{\BM}^Z(X), ^{\vee}, \varpi}  
& in~case~n-k=0~mod~4\\
{\BG}W\LRf{C{\BM}^Z(X), ^{\vee}, \varpi}^- 
& in~case~n-k=2~mod~4\\
\end{array}\right.
\qquad \forall n\in {\BZ}
$$
Here on right side, the duality on the exact categories  is given by $M\mapsto 
{\SE}xt^k\LRf{M, {\CL}}$.
}
\eT
\pf The latter part follows from former, using  discussion in section \ref{BISP4Exact}.
To prove the first part, consider the  commutative diagram 
$$
\diagram 
{\BG}W^{[n-k]} \LRf{{\bf dg}C{\BM}^Z(X)} \ar[r]\ar[d] & {\BG}W^{[n]} \LRf{{\bf dg}{\SV}(X)}\ar[r] \ar[d]&  {\BG}W^{[n]} \LRf{{\bf dg}{\SV}(U)}\ar[d]\\
{\BG}W^{[n]} \LRf{{\bf dg}_Z{\SV}(X)} \ar[r] & {\BG}W^{[n]} \LRf{{\bf dg}{\SV}(X)}\ar[r] &  {\BG}W^{[n]} \LRf{{\bf dg}{\SV}(U)}\\
\enddiagram
$$
of Homotopy fibration in $\BiSp$. 
 The top left map here is the zig-zag map given in (\ref{BiSpectrEquivho}). 
The upper fibration follows from (\ref{BiSpectrEquivho}), and the lower fibration follows from \cite[pp 1807]{S17}. The $2^{nd}$ and $3^{rd}$ vertical maps are identity. Hence the first vertical map is a homotopy equivalence in $\BiSp$. \pic $\eop$

\vspace{3mm}


\vspace{3mm}
The Excision and Mayer-Vietories theorems would be routine consequence, as follows.
\bT[Excision]\label{exiKThr}{\rm 
Let $X$ be a quasi projective scheme over a noetherian affine scheme $\spec{A}$. Let $V\subseteq X$ be an open subset, and $Z\subseteq X$ be a closed subset such that 
$Z\subseteq V$. Assume that $Z$ is complete intersection, as in (\ref{CIsetUp}).
Then the restriction map
$$
{\BK}\LRf{C{\BM}^Z(X)}\lra {\BK}\LRf{C{\BM}^Z(V)}\qquad in \qquad \Sp
$$
is a homotopy equivalence of spectra. Consequently, the induced maps
$$
{\BK}_p\LRf{C{\BM}^Z(X)}\iso {\BK}_p\LRf{C{\BM}^Z(V)}\quad
{\rm are~ isomorphisms}~\forall p\in {\BZ}.
$$
}
\eT
\pf One can use the Agreement theorem \ref{AgrBkt}, and refer to the Excision theorem \cite[Thm 3.4.10]{S11}. For a more direct proof, note that
by agreement  and equivalence theorems \cite[pp 409]{M23}, we are reduced to proving that 
the restriction functor of derived categories ${\bfD}^b\LRf{C{\BM}^Z(X))} \iso {\bfD}^b\LRf{C{\BM}^Z(V))}$  is an equivalence. However, the functor 
$$
C{\BM}^Z(X) \lra C{\BM}^Z(V)\quad {\rm is~an~equivalence, in~fact, ~is~a~bijection}.
$$
The inverse functor is given by extending $M\in C{\BM}^Z(V)$ by zero, on $X-Z$. 
\pic $\eop$ 


\vspace{3mm}
Same proof leads to a ${\BG}W$-analogue of (\ref{exiKThr}).
\bT[Excision]\label{BgWExci}{\rm 
Let $X$ be a quasi projective scheme over a noetherian affine scheme $\spec{A}$. Let $V\subseteq X$ be an open subset, and $Z\subseteq X$ be a closed subset such the 
$Z\subseteq V$. Assume $1/2\in A$ and ${\CL}$ is a line bundle on $X$. 
Further assume that $Z$ is complete intersection on $X$, as in (\ref{CIsetUp}). 
Then the restriction map
$$
{\BG}W^{[n]}\LRf{{\bf dg}C{\BM}^Z(X), {\CL}}\lra {\BG}W^{[n]}\LRf{{\bf dg}C{\BM}^Z(V),
{\CL}_{|V}}
\qquad in \qquad \BiSp
$$
is a homotopy equivalence in the category $\BiSp$ of bispectra. Consequently, the induced maps
$$
{\BG}W^{[n]}_r\LRf{{\bf dg}C{\BM}^Z(X), {\CL}}\iso {\BG}W^{[n]}_r\LRf{{\bf dg}C{\BM}^Z(V), {\CL}_{|V}}\quad
{\rm are~ isomorphisms}~\forall n, r\in {\BZ}.
$$
Here on either side, with $k=grade(Z, X)$, the duality is given by
$$
\left\{\begin{array}{ll}
M^{\vee}={\SE}xt^k\LRf{M, {\CL}}& M\in C{\BM}^Z(X)\\
M^{\vee}={\SE}xt^k\LRf{M, {\CL}_{|V}}& M\in C{\BM}^Z(V)\\
\end{array}\right.
$$
}
\eT
\pf Again the map $C{\BM}^Z\LRf{X} \lra C{\BM}^Z\LRf{V}$ is an equivalence
(In fact, a bijection). Hence the it induced functor ${\bfD}^b\LRf{C{\BM}^Z\LRf{X}} 
\lra {\bfD}^b\LRf{C{\BM}^Z\LRf{V}}$ is a duality preserving equivalence. Therefore, the proof is complete by the Invariance theorem \cite[pp 1800]{S17}, \cite[pp 466]{M23}. $\eop$ 

\vspace{3mm}
The following are exact category version of the above (\ref{BgWExci}).
\bC[Excision]\label{2364ExciEcat}{\rm 
Consider the setup and the hypotheses of (\ref{BgWExci}), then the following maps
$$
\left\{\begin{array}{l}
{\BG}W\LRf{C{\BM}^Z(X)}\lra {\BG}W\LRf{C{\BM}^Z(V)}\\
{\BG}W\LRf{C{\BM}^Z(X)}^-\lra {\BG}W\LRf{C{\BM}^Z(V)}^-\\
\end{array}\right.
\qquad in \qquad \BiSp
$$
are homotopy equivalences in $\BiSp$.
}
\eC
\pf Immediate from (\ref{BgWExci}), by taking $n=0, 2$. \pic $\eop$ 

\vspace{3mm} 
The Mayer-Vietories Theorem is, in fact, a routine consequence of excision theorem. So, as a consequence of (\ref{exiKThr}), we obtain the following.
\bT[Mayer-Vietoris]\label{MaViBK}{\rm 
Let $X$ be a quasi projective scheme over a noetherian affine scheme $\spec{A}$, and let $X=U \cup V$ be an open cover. Assume $Z=X-U$ is complete intersection in $X$,
as in (\ref{CIsetUp}).
Then the diagram 
$$
\diagram 
{\BK}({\SV}(X)) \ar[r] \ar[d] & {\BK}({\SV}(U)) \ar[d]\\
{\BK}({\SV}(V)) \ar[r] & {\BK}({\SV}(U\cap V)) \\
\enddiagram
\quad {\rm is ~a~cartesian ~square ~in~the~category}~\Sp~{\rm of ~spectra}.
$$
Consequently, there is a long exact sequence
$$
\diagram 
\cdots \ar[r] & {\BK}_n({\SV}(X))\ar[r] & {\BK}_n({\SV}(U))\bigoplus {\BK}_n({\SV}(V))
\ar[r] & 
{\BK}_n({\SV}(U\cap V))\\
\ar[r] & {\BK}_{n-1}({\SV}(X))\ar[r] & \cdots&\forall n\in {\BZ}\\
\enddiagram
$$
of ${\BK}$-groups. 
}
\eT
\pf Write $Z=X-U=V-U\cap V\subseteq V$. Consider the diagram
$$
\diagram 
{\BK}(C{\BM}^Z(X)) \ar[r] \ar[d]_{\wr} & {\BK}({\SV}(X)) \ar[r] \ar[d] & {\BK}({\SV}(U)) \ar[d]\\
{\BK}(C{\BM}^Z(V)) \ar[r]  &{\BK}({\SV}(V)) \ar[r] & {\BK}({\SV}(U\cap V)) \\
\enddiagram
\qquad in \qquad \Sp.
$$
Both the horizontal lines are homotopy fibrations in $\Sp$, by Theorem \ref{Equi1909ALL}. The $1^{st}$ vertical map is an isomorphism in $\Sp$, by Excision theorem \ref{exiKThr}. This completes the proof that the given square is cartesian. 
The long exact sequence above is established routinely, by comparing two exact sequences:
$$
\diagram 
{\BK}_n(C{\BM}^Z(X)) \ar[r] \ar[d]_{\wr} & {\BK}_n({\SV}(X)) \ar[r] \ar[d] & {\BK}_n({\SV}(U)) \ar[d]\ar[r] &{\BK}_{n-1}(C{\BM}^Z(X))  \ar[d]_{\wr} \\
{\BK}_n(C{\BM}^Z(V)) \ar[r]  &{\BK}_n({\SV}(V)) \ar[r] & {\BK}_n({\SV}(U\cap V))\ar[r] & {\BK}_{n-1}(C{\BM}^Z(V)) \\
\enddiagram
$$
\pic $\eop$


\vspace{3mm} 
The same proof of (\ref{MaViBK}) gives a Mayer-Vietories theorem for ${\BG}W$ theory, as follows.

\bT[Mayer-Vietoris]\label{MaViBK}{\rm 
Let $X$ be a quasi projective scheme over a noetherian affine scheme $\spec{A}$, and let $X=U \cup V$ be an open cover. Assume $1/2\in A$ and ${\CL}$ be a line bundle on $X$. Assume $Z=X-U$ is complete intersection in $X$,
as in (\ref{CIsetUp}).  Then the diagram 
$$
\diagram 
{\BG}W^{[n]}({\bf dg}{\SV}(X), {\CL}) \ar[r] \ar[d] & {\BG}W^{[n]}({\bf dg}{\SV}(U), {\CL}_{|U}) \ar[d]\\
{\BG}W^{[n]}({\bf dg}{\SV}(V), {\CL}_{|U}) \ar[r] & {\BG}W^{[n]}({\bf dg}{\SV}(U\cap V), {\CL}_{|U\cap V}) \\
\enddiagram
\qquad in \qquad \BiSp
$$
is a cartesian square in the category $\BiSp$ of bispectra $\forall n\in {\BZ}$. 
Consequently, there are long exact sequences $\forall n\in {\BZ}$ of ${\BG}W$-groups
(we ignore to include the ${\CL}$-coordinate):
$$
\diagram 
\cdots \ar[r] & {\BG}W^{[n]}_p({\SV}(X))\ar[r] & {\BG}W^{[n]}_p({\SV}(U) )\bigoplus {\BG}W_n({\SV}(V))\ar[r] &\\
 {\BG}W^{[n]}_p({\SV}(U\cap V))\ar[r] & {\BG}W^{[n]}_{p-1}({\SV}(X))\ar[r] & \cdots &
\\
\enddiagram
$$
}
\eT
\pf 
Write $Z=X-U=V-U\cap V\subseteq V$, and $k=grade(Z, X)$. Consider the diagram
$$
\diagram 
{\BG}W^{[n]}({\bf dg}C{\BM}^Z(X)) \ar[r] \ar[d]_{\wr} & {\BG}W^{[n]}({\bf dg}{\SV}(X)) \ar[r] \ar[d] & {\BG}W^{[n]}({\bf dg}{\SV}(U)) \ar[d]\\
{\BG}W^{[n]}({\bf dg}C{\BM}^Z(V)) \ar[r]  &{\BG}W^{[n]}({\bf dg}{\SV}(V)) \ar[r] & {\BG}W^{[n]}({\bf dg}{\SV}(U\cap V)) \\
\enddiagram
$$
Both the horizontal lines are homotopy fibrations in $\BiSp$, by Theorem \ref{BiSpectrEquivho}. The $1^{st}$ vertical map is an isomorphism in $\BiSp$, by Excision  theorem \ref{BgWExci}. This completes the proof that the given square is cartesian. Again, the exact sequences are established by routine methods. \pic $\eop$

The following is the exact category version of (\ref{MaViBK}).
\bT[Mayer-Vietoris]\label{MxaVbinBK}{\rm 
Consider the setup and the hypotheses of (\ref{MaViBK}).
 With $\delta=+, -$,  the diagrams 
$$
\diagram 
{\BG}W({\SV}(X), {\CL})^{\delta} \ar[r] \ar[d] & {\BG}W({\SV}(U), {\CL}_{|U})^{\delta} \ar[d]\\
{\BG}W({\SV}(V), {\CL}_{|U})^{\delta} \ar[r] & {\BG}W({\SV}(U\cap V), {\CL}_{|U\cap V})^{\delta} \\
\enddiagram
\qquad in \qquad \BiSp
$$
are cartesian square in the category $\BiSp$ of bispectra. 
Consequently, there are long exact sequences $\forall p\in {\BZ}$ of ${\BG}W$-groups (with $\delta=+, -$): 
$$
\diagram 
\cdots \ar[r] &{\BG}W_p\LRf{{\SV}(X)}^{\delta} \ar[r] & {\BG}W_p\LRf{{\SV}(U)}^{\delta} \bigoplus {\BG}W_p\LRf{{\SV}(V)}^{\delta} \ar[r] & {\BG}W_p\LRf{{\SV}(U \cap V)}^{\delta}\\
 \ar[r] &{\BG}W_{p-1}\LRf{{\SV}(X)}^{\delta} \ar[r]&\cdots& \\
\enddiagram
$$

%
}
\eT
\pf Follows from (\ref{MaViBK}), with $n=0, 2$. \pic $\eop$ 

 
\appendix


\section{Grade and annihilator of  coherent modules} \label{SecGrade}
 While grade of coherent modules over affine schemes are fairly standard, 
 in this section, we discuss grade of coherent modules over quasi projective schemes.
First, we recall the definition of grade and summarize the basics on the same from \cite[pp. 103]{Mh}.
\bD\label{gradeMatr103}{\rm 
Suppose $A$ is a noetherian ring, $X=\spec{A}$ and $M\in Coh(A)$. 
\bE
\item Define
$$
grade(M)=\inf\left\{r:Ext^r(M, A)\neq 0 \right\}
$$
Consequently, for an ideal $I$, define
\begin{equation}\label{grdOfIii}
grade\left(\frac{A}{I}\right)= \inf\left\{r:Ext^r\LRf{\frac{A}{I}, A}\neq 0 \right\}
\end{equation}
Sometimes, one writes 
$$
grade(I):=
grade\left(\frac{A}{I}\right)
$$
This is not to be confused with the grade of $I$, as an $A$-module, which would be zero if $I$ contains a non zero divisor. 
\item 
It follows, for $M\in Coh(A)$ and $I =Ann(M)$, we have 
\begin{equation}\label{grdDrei396}
grade(M)=depth_I(A)=grade\left(\frac{A}{I}\right) 
\end{equation}
\item It is obvious from definition 
$$
grade(M) \leq \PDV(M)
$$
\eE
}
\eD
\bC\label{grdRAD}{\rm 
For an ideal $I$ of $A$, we have 
$$
grade(I)=grade\left(\sqrt{I}\right)
$$
}
\eC
\pf Let $r=grade(I)$ and $s=grade\left(\sqrt{I}\right)$.
First,  $r=depth_I(A)\leq depth_{\sqrt{I}}(A)=s$. 
 If $a_1, \ldots a_s\in \sqrt{I}$ is a regular sequence then so 
 is $a_1^m, \ldots a_s^m\in I$ for $m\gg 0$,  \cite[Thm 26, pp. 96]{Mh}. So, $s\leq r$. So, $s=r$. 
 \pic $\eop$ 

\vspace{3mm}
We define grade of a closed subset of $\spec{A}$. 
\bD\label{GrdClsdsset}{\rm 
Let $A$ be a noetherian ring, and $X=\spec{A}$ Let $Z=V(I)\subseteq X$ be a closed subset. Define
$$
grade(Z, X)=grade(I)=grade\left(\frac{A}{I}\right)
$$
By (\ref{grdRAD}), it is well defined and is independent of the choice of $I$. 
}
\eD

We record the following non local version of \cite[pp. 108, Thm 31]{Mh}. 
\bL[pp. 108, Thm 31 \cite{Mh}]{\rm 
Suppose $A$ is a Cohen Macaulay ring (non local) and $M\in Coh(A)$. Then 
$$
grade(M)=height(I) \qquad {\rm where} \qquad I=Ann(M).
$$
}
\eL 

\vspace{3mm}
Subsequently, we wish to make some sense out of the idea of annihilators of coherent sheaves. 
\bD\label{BDAnniSF}{\rm 
Let $X$ be a noetherian scheme and ${\SF}\in Coh(X)$. 
Define the annihilator
$$
ann({\SF})(U)=\left\{t\in \Gamma\left(U, {\CO}_X\right): 
t\cdot \Gamma\left(U, {\SF}\right)=0\right\} \subseteq \Gamma\left(U, {\CO}_X\right)
$$
The definition makes sense.
However, $ann({\SF})$ does not seem to make sense, as a (pre)sheaf of ideals.
This is because the restriction maps 
$\Gamma\left(U, {\CO}_X\right) \lra \Gamma\left(V, {\CO}_X\right)$ 
may not give a restriction map on $ann\left({\SF}\right)$.

It is a standard fact that it makes sense when $X=\spec{A}$ is affine. 
Consequently,  in a noetherian scheme $X$, on the affine pieces,
$U_i=\spec{A_i}=\spec{\Gamma\left(U_i, {\CO}_X\right)}$ and 
$M_i=\Gamma\left(U_i, {\SF}\right)$ the annihilators 
$I_i=ann(M_i)$ make sense. 

Subsequently, we prove that $ann({\SF})$ makes sense when $X$ is quasi projective. 
}
\eD

\subsection{(Quasi) projective schemes} 

 \bL\label{theANNhProj}{\rm Suppose $Y=\proj{S}$ with $S=A[x_0, x_1, \ldots, x_N]$ and $\deg(x_i)=1$. Assume $A$ is noetherian. 
Let ${\SF}\in Coh(Y)$ and ${\SF}=\widetilde{M}$ 
for some finitely generated graded $S$-module $M$. 
({\it Refer to \cite[\S II.5, pp116]{H} for this $\widetilde{M}$ construction}.) 
Let $I=ann(M)$. Then $\widetilde{I}=ann({\SF})$.
}
\eL 
\pf First we need to prove that $ann({\SF})$ is a sheaf. From the affine case it is clear that $ann({\SF})_{|D(g)}= \widetilde{ann\left(M_{(g)}\right)}~\forall~g$ homogeneous. Now suppose $V\subseteq U$.
We need to give a map 
$$
ann\left({\SF}_{|U}\right) \lra ann\left({\SF}_{|V}\right) 
$$
First assume $V=\spec{S_{(f)}}$, for some homogeneous $f\in S$. 
Write $U=\bigcup_{i=1}^n\spec{S_{(g_i)}}$. There are maps 
$$
\left\{\begin{array}{l}
S_{(g_i)}\lra S_{(fg_i)}\\
 M_{(g_i)}\lra M_{(fg_i)}\\
\end{array}\right. 
\Lra  \left\{\begin{array}{l}
ann\left(M_{(g_i)}\right)\lra ann\left(M_{(fg_i)}\right) \\
{\rm sending}\quad \frac{a}{g_i^k} \mapsto  \frac{a}{g_i^k}\quad with \quad \deg(a)=\deg(g_i^k) \\
\end{array}\right. 
$$
In deed, the commuting square
$$
\diagram 
ann\left(M_{(g_i)}\right)\ar[r]\ar@{^(->}[d]  & ann\left(M_{(fg_i)}\right)\ar@{^(->}[d]   \\
S_{(g_i)}\ar[r] & S_{(fg_i)}\\
\enddiagram\Lra
\diagram 
ann\left({\SF}\right)_{|D(g_i)}\ar[r]\ar@{^(->}[d]  & ann\left({\SF}\right)_{|D(fg_i)}\ar@{^(->}[d]   \\
{\CO}_{|D(g_i)}\ar[r] & {\CO}_{|D(fg_i)}\\
\enddiagram\quad {\rm commutes.}
$$
Let $\lambda \in ann\LRf{{\SF}_{|D(g_i)}}$.
It follows
$\lambda_{|D(fg_i)}\in ann\LRf{{\SF}_{|D(fg_i)}}~\forall i$. 
So, given $m\in \Gamma\LRf{U, {\SF}}$ we have $(\lambda\cdot m)_{|D(fg_i)}=0 \forall i$. So, $\lambda_{|D(f)}
\in ann\LRf{{\SF}_{|D(f)}}$.
%
These patch to give a map
\begin{equation}\label{440ANN}
ann\left({\SF}_{|U} \right)\lra ann\left({\SF}_{|D(f)}\right)
\end{equation}
Now suppose 
$V=\bigcup_{j=1}^t D(f_j)$. Given $\lambda \in ann\LRf{\Gamma(U, {\SF)}}$ and 
$m\in \Gamma(V, {\SF})$, it follows from above that 
$(\lambda\cdot m)_{|D(f_j)}=0~\forall j$. So, $(\lambda\cdot m)_{|V}=0$. So, the map
$\Gamma\LRf{U, {\CO}_X}\lra \Gamma\LRf{V, {\CO}_X}$ restricts, as follows:
$$
\diagram 
ann\LRf{\Gamma(U, {\SF})} \ar[r]\ar@{^(->}[d] & ann\LRf{\Gamma(V, {\SF})}\ar@{^(->}[d] \\
\Gamma(U, {\CO}_X) \ar[r] & \Gamma(V, {\CO}_X)\\
\enddiagram
$$
So, $ann({\SF})\subseteq {\CO}_X$ is an ideal sheaf. Now consider the natural map
$$
\iota: \widetilde{I} \lra ann\LRf{{\SF}}. 
$$
If follows from the affine case that the restriction $\iota_{|D(g)}$ is an isomorphism, for all homogeneous $g\in S$. \pic $\eop$

\vspace{3mm}
 \bC\label{CorAnnGrQP}{\rm Suppose $Y=\proj{S}$ with $S=A[x_0, x_1, \ldots, x_N]$ and $\deg(x_i)=1$. Assume $A$ is noetherian. Let $X\subseteq Y$ be open and
 ${\SF}\in Coh(X)$. Let ${\SG}:=\widetilde{M}\in Coh(Y)$ be an extension of ${\SF}$ 
 (see \cite[Ex II.5.15]{H}) and 
 $I =ann(M)$. Then $ann({\SF})=\widetilde{I}_{|X}$.
 
}
\eC
\pf Follows from (\ref{theANNhProj}). $\eop$

 \bD{\rm 
Let $X$ be a noetherian scheme and ${\SF}\in Coh(X)$. Define 
$$
grade({\SF})=\inf\left\{r:{\SE}xt^r\left({\SF}, {\CO}_X\right)\neq 0 \right\}
$$
}
\eD


\vspace{3mm}
 \bP\label{PSgrdNht}{\rm Suppose $Y=\proj{S}$ with $S=A[x_0, x_1, \ldots, x_N]$ and $\deg(x_i)=1$. Assume $A$ is noetherian. Let $X\subseteq Y$ be open,
 ${\SF}\in Coh(X)$ and ${\SG}=\widetilde{M}$ where $M$ is a finitely generated graded $S$-module, such that ${\SG}_{|X}={\SF}$ (see \cite[Ex II.5.15]{H}).
Let $I =ann(M)$ and ${\SI}=\widetilde{I}=ann({\SG})\subseteq {\CO}_Y$.  Then 
$$
grade({\SF}) = grade\left(\left(\frac{{\CO}_Y}{{\SI}} \right)_{|X}\right)=
grade\left(\frac{{\CO}_X}{ann{\SF}} \right)
$$
}
\eP
\pf First note that, for any ${\SE}, {\SG}\in Coh(X)$ we have  (see \cite[Prop III.6.8]{H})
$$
{\SE}xt^i\LRf{{\CE}, {\CG}}_x = Ext^i\LRf{{\SE}_x, {\SG}_x}\qquad \forall x\in X, i\geq 0.
$$
Therefore,
$$
\begin{array}{l}
grade({\SF})=\inf\left\{grade\left({\SF}_y\right): y\in V(I)\cap X  \right\}\\
=\inf\left\{grade\left(\frac{{\CO}_{Y}}{{\SI}}\right)_y: y\in V(I)\cap X  \right\}=
grade\left(\left(\frac{{\CO}_Y}{{\SI}} \right)_{|X}\right)\\
\end{array}
$$
This middle equality follows from (\ref{grdDrei396}). $\eop$ 

\vspace{3mm}
 \bP\label{PSgRegSeq}{\rm Suppose $Y=\proj{S}$ with $S=A[x_0, x_1, \ldots, x_N]$ and 
 $\deg(x_i)=1$. Assume $A$ is noetherian. Let $X\subseteq Y$ be open and
 ${\SF}\in Coh(X)$. Let ${\SG}=\widetilde{M}$ where $M$ is a finitely generated graded $S$-module, such that ${\SG}_{|X}={\SF}$ (see \cite[Ex II.5.15]{H}). 
Let $I =ann(M)$ and ${\SI}=\widetilde{I}=ann({\SG})\subseteq {\CO}_Y$. Then the following are equivalent:
\bE
\item \label{592ONE}
$
grade({\SF}) \geq r. 
$
\item\label{596TWO} There is a sequence $f_1, f_2, \ldots, f_r\in I$, of homogeneous elements, such that they induce a regular sequence on $X$. 
 ({\it This can be read as, if ${\wp}\in X$ and $x_i\notin {\wp}$, then
$$
\frac{f_1}{x_i^{\deg(f_1)}}, \frac{f_2}{x_i^{\deg(f_2)}},, \ldots, \frac{f_r}{x_i^{\deg(f_r)}}
$$ 
is a $S_{(\wp)}$-regular sequence. })

\eE 
}
\eP
\pf ($(\ref{596TWO})\Lra (\ref{592ONE})$:) it follows from (\ref{596TWO}) and (Equation \ref{grdDrei396}) that $grade({\SF}_x)=depth_{I_x}S_{({\wp})}\geq r$ and hence 
$Ext^i\left({\SF}_x, {\CO}_x\right)=0 ~\forall x\in X, ~i\leq r-1$. So, 
${\SE}xt^i({\SF}, {\CO}_X)=0 ~\forall i\leq r-1$. Hence $grade({\SF}) \geq r$. So, (\ref{592ONE})
is established. 

\noindent($(\ref{592ONE})\Lra (\ref{596TWO})$:) 
Suppose $f_1, f_2, \ldots, f_p\in I$ is a sequence of homogeneous elements, with $p\leq r-1$, that induce a $S_{({\wp})}$-regular sequence. Let $J=(f_1, f_2, \ldots, f_p)$. Let
$$
Ass\left(\frac{S}{J} \right)\cap X=\left\{{\wp}_1, \ldots, {\wp}_q\right\}
$$
Suppose $I\subseteq \bigcup_{i=1}^q{\wp}_i$. Then $I\subseteq {\wp}={\wp}_i$. Then $I_{(\wp})
\subseteq {\wp}S_{(\wp)}$. Assume $x_0\notin {\wp}$. Since $depth_{I_{({\wp})}}S_{(\wp)}$ is well defined, we obtain $depth_{I_{({\wp})}}S_{(\wp)}=p$. Hence $grade\left({\SF}_{\wp}\right)=depth_{I_{({\wp})}}S_{(\wp)}=p$. 
Hence $Ext^p\left({\SF}_{\wp}, {\CO}_{X, {\wp}}\right)\neq 0$. Therefore 
${\SE}xt^p\left({\SF}, {\CO}_{X,}\right)\neq 0$, and $grade({\SF})\leq p\leq r-1$, which is a contradiction. 

So, $I\not\subseteq \bigcup_{i=1}^q{\wp}_i$. Let  
$f_{p+1}\in I -\left( \bigcup_{i=1}^q{\wp}_i\right)$. Then $f_1, \ldots, f_p, f_{p+1}\in I$ induce a 
regular sequence, as required.
\pic $\eop$ 

\vspace{3mm}
 \bP\label{CohenHtt}{\rm Suppose $Y=\proj{S}$ with $S=A[x_0, x_1, \ldots, x_N]$ and 
 $\deg(x_i)=1$. Assume $A$ is noetherian. Let $X\subseteq Y$ be open and
 ${\SF}\in Coh(X)$. Let ${\SG}=\widetilde{M}$ where $M$ is a finitely generated graded $S$-module, such that ${\SG}_{|X}={\SF}$ (see \cite[Ex II.5.15]{H}).
%
Assume $X$ is Cohen-Macaualy. Then 
$$
grade({\SF}) =co\dim(V(I)\cap X, X) \qquad \qquad {\rm where}\quad I=ann(M).
$$
}
\eP

 \vspace{3mm}
 We extend the definition (\ref{GrdClsdsset}) and  define the grade of a closed subset 
 of a quasi projective scheme $X$. 
\bD\label{GradCloQA}{\rm 
Let $X$ be a quasi projective scheme of over a noetherian affine scheme $\spec{A}$ and $Z\subseteq X$ be a closed subset of $X$. Further assume 
 $Y=\proj{S}$ with $S=A[x_0, x_1, \ldots, x_N]$,  
 $\deg(x_i)=1$, and $X\subseteq Y$ is an open subset. 
 Define
 $$
grade(Z):= grade(Z, X) =grade\left(\frac{{\CO}_X}{\tilde{I}_{|X}}\right)
 \qquad {\rm where}
 \left\{\begin{array}{l}
  Z=V(I)\cap X,~I\subseteq S \\
  {\rm is~a~homogenous~ideal.}\\
  \end{array}\right. 
 $$
 This definition is independent of choice of $Y$ and $I$. 
 }
 \eD
\pf Let $Z=V(I)\cap X=V(J)\cap X$. Considering  the irreducible components of $Z$, we have
$$
X\cap \min(I)=X\cap \min(J)=\left\{{\wp}_1, \ldots, {\wp}_r \right\}
$$
So, we can assume $J=\cap_{i=1}^r{\wp}_i$ is reduced and that $I \subseteq J$.
Now
$$
grade\left(\frac{{\CO}_X}{\tilde{I}_{|X}}\right)=\inf\left\{r: {\SE}xt^r\left(\frac{{\CO}_X}{\tilde{I}_{|X}}, {\CO}_X\right)\neq 0 \right\}
=
\inf\left\{r: {\SE}xt^r\left(\frac{{\CO}_x}{\tilde{I}_{x}}, {\CO}_x\right)\neq 0, x\in X \right\}
$$
Using (\ref{grdRAD}) or \cite[Thm 26, pp. 98]{Mh}, the above
$$
=\inf\left\{r: {\SE}xt^r\left(\frac{{\CO}_x}{\tilde{J}_{x}}, {\CO}_x\right)\neq 0, x\in X \right\}=
grade\left(\frac{{\CO}_X}{\tilde{J}_{|X}}\right)
$$
\pic $\eop$ 

We recall the definition of perfect modules, before we close this section.
\bD\label{defPerModd}{\rm 
Let $X$ be a noetherian scheme and $Z\subseteq X$ be a closed subset. For ${\SF}\in Coh(X)$, we say that ${\SF}$ is a {\bf perfect module}, if $grade\LRf{{\SF}}=\PDV{{\SF}}$. Denote the full subcategory $C{\BM}^Z(X)=\LRs{{\CF}\in Coh_Z(X): \PDV{{\SF}}=grade(Z, X)}$. Note  $C{\BM}^Z(X)$ is a subcategory of the category of perfect modules. 
}\eD
\section{Hermitian theory and Karoubi ${\BG}W$ bispectrum}\label{BISP4Exact}

The purpose of this section is to
 provide a brief background on jargon of homotopy approach to  $K$ theory, in the post Quillen \cite{Q73} era. 
 Hermitian $K$ theory (also called  $GW$-theory) is to be seen as part of Algebraic $K$-theory. However, we assume readers are more familiar with $K$-theory, and we will discuss more about $GW$-theory,  in this section.
 We assume that the readers have some familiarity with the background  \cite{Q73, W83, TT90} on $K$-theory, or are referred to \cite{M23}. 
 We also assume that the readers have some familiarity with the category $\Sp$ of spectra of pointed topological spaces. 
Quillen established Algebraic $K$-Theory on the foundation of Exact categories. 
 In our view, for the esthetics of the narrative, Algebraic $K$-Theory should be a story of exact categories. Formulation of Theorem \ref{GrayQuillen} of Quillen supports the same sense of esthetics. 
We like to make clean statements within the universe of exact categories, when possible.
However, after \cite{Q73}, the literature seemed to have drifted away from $K$-theory  of exact categories ${\SE}$, to that of the category of bounded chain complexes $Ch^b\LRf{{\SE}}$ \cite{W83, TT90}. The theory further drifted into other derivatives of exact categories, like complecial exact categories \cite{S11} and dg categories \cite{S17}. While these have enriched the theory, they also have obscured the theory for the broader readership. 
The literature misses any common umbrella for all these different derivatives, while $K$-theory of exact categories (inclusive of $GW$-theory) remains the common denominator  or, of the fundamental interest. 
For this reason, in this section, we  recapitulate some of the background on $K$-theory,
mainly of $GW$ theory,
and reinterpret them in the context of exact categories. The main goal of this section, is to explain the 
Karoubi ${\BG}W$ bispectrum ${\BG}W\LRf{{\SE}_d}\in \BiSp$ of exact categoeies ${\SE}_d:=\LRf{{\SE}, ^{\vee}, \varpi}$ with duality. 
 
 Let $\eCat$  (resp. $\widetilde{\eCat}$) denote the category of small
(resp. idempotent complete)  exact categories. Likewise,  ${\eCat}D$
(resp.  $\widetilde{\eCat}D$) denote the 
  category of small (resp. idempotent complete) exact categories with duality, over ${\BZ}(1/2)$.
 We start with the following definition.
 \bD{\rm
 Let ${\SE}_d:=\LRf{{\SE}, ^{\vee}, \varpi}\in {\eCat}D$. 
 Consider the following homotopy fibration: 
$$
\diagram
{\bf GW}\LRf{{\SE}_d} \ar[r] & {\BB}{\BQ}^h\LRf{{\SE}_d} \ar[r] & {\BB}{\BQ}\LRf{{\SE}}\\
\enddiagram{\rm ,\quad and}~\quad 
\diagram
 {\bfK}\LRf{{\SE}}\ar[r]&
{\bf GW}\LRf{{\SE}_d} \ar[r] & {\BB}{\BQ}^h\LRf{{\SE}}\\
\enddiagram
$$
Here $ {\bfK}\LRf{{\SE}}$ denotes the $K$-theory space of ${\SE}$ and  
${\bf GW}\LRf{{\SE}_d}$ denotes the GW space of ${\SE}_d$. Ideally, we should denote ${\bf GW}\LRf{{\SE}_d}:= {\bf GW}\LRf{{\SE}, ^{\vee}, \varpi}$. The latter sequence is a rotation of the first one. 
 Higher $K$-groups and $GW$-groups are defined by 
$$
\left\{\begin{array}{l}
K_n\LRf{{\SE}}=\pi_n\LRf{{\bfK}\LRf{{\SE}}}\\
GW_n\LRf{{\SE}_d}=\pi_n\LRf{{\bf GW}\LRf{{\SE}}}\\
\end{array}\right.
\qquad \forall n\geq 0.
$$
}
\eD
\subsection{$K$-theory and $GW$-theory spectrum}\label{subSecKTandGWSP}
The work of Waldhausen \cite{W83} considered $K$-theory as a spectrum. 
We will not discuss the connective $K$-theory spectrum. In this section we give a overview of non connective $K$-theory and $GW$-theory spectrum. 
 Cone functors 
  $$
  \left\{\begin{array}{l}
  \C_K: \widetilde{\eCat} \lra \eCat\\
  \C_{GW}: \widetilde{\eCat}D \lra {\eCat}D\\
  \end{array}\right.
  \quad {\rm together~with}\quad 
   \left\{\begin{array}{ll}
   {\SE} \hra \C_K\LRf{{\SE}} & \forall {\SE}\in \widetilde{\eCat}\\
      {\SE} \hra \C_{GW}\LRf{{\SE}} & \forall {\SE}\in \widetilde{\eCat}D\\
     \end{array}\right.
  $$
 fully faithful exact functors (duality preserving in the latter case),  were defined.
 Further, for ${\SE}\in \widetilde{\eCat}$ and ${\SE}_d:=\LRf{{\SE}, ^{\vee}, \varpi} \in \widetilde{\eCat}W$ 
 \bE
 \item  The respective quotients 
  $\frac{\C_K{\SE}}{{\SE}}\in  {\eCat}$ and 
 $\frac{\C_{GW}{\SE}_d}{{\SE}_d}\in  {\eCat}D$ are defined. 
 \item The spaces
  \begin{equation}\label{CotBle2289}
  \left\{\begin{array}{l}
  {\bfK}\LRf{\C_K({\SE})}\cong \star\\
    {\bf GW}\LRf{\C_{GW}\LRf{{\SE}_d}}\cong \star\\
  \end{array}\right.~{\rm are~contractible}.
  \end{equation} 
  \item The sequences,
  \begin{equation}\label{DonFibrs}
  \begin{array}{l}
  \diagram
{\bfK}\LRf{{\SE}} \ar[r] & {\bfK}\LRf{\C_K{{\SE}}} \ar[r] & {\bfK}\LRf{\frac{\C_K{{\SE}}}{{\SE}}}  \\
\enddiagram \\
 \diagram
{\bf GW}\LRf{{\SE}} \ar[r] & {\bf GW}\LRf{\C_{GW}{{\SE}_d}} \ar[r] & 
{\bf GW}\LRf{\frac{\C_{GW}{{\SE}_d}}{{\SE}_d}}  \\
\enddiagram \\
  \end{array}
  \quad {\rm are~homotopy~fibrations}.
  \end{equation} 
  \eE 
We define suspension functors, as follows.
\bD{\rm 
Let ${\SE}\in \eCat$ and ${\SE}_d:=\LRf{{\SE}, ^{\vee}, \varpi}\in {\eCat}D$. Let 
$\widetilde{{\SE}}$ and $\widetilde{{\SE}_d}$  denote the respective idempotent completion.
Define suspensions
$$
\left\{\begin{array}{l}
{\BS}_K\LRf{{\SE}}:= \frac{\C_{K}{\widetilde{\SE}}}{\widetilde{\SE}}\\
{\BS}_{GW}\LRf{{\SE}_d}:= \frac{\C_{K}{\widetilde{\SE}_d}}{\widetilde{\SE}_d}\\
\end{array}\right. 
\quad {\rm Inductively,~define}\quad 
\left\{\begin{array}{l}
{\BS}_K^{n+1}\LRf{{\SE}}:={\BS}_K\LRf{\widetilde{{\BS}_K^n({\SE})}}\\ 
{\BS}_{GW}^{n+1}\LRf{{\SE}_d}:= {\BS}_{GW}^n\LRf{\widetilde{{\BS}_{GW}^n({\SE}_d)}}\\
\end{array}\right. 
$$
We define two non connective spectrums, as follows:
\bE
\item 
The non connective ${\BK}$-theory spectrum of ${\SE}$ is defined as
$$
{\BK}\LRf{{\SE}}=\LRs{
{\BK}_n\LRf{{\SE}}: n\geq 0
},\quad {\rm where}\quad {\BK}_n\LRf{{\SE}}={\bfK}\LRf{\widetilde{{\BS}_K^n\LRf{{\SE}}}}
$$
The ${\BK}$-groups are defined by 
$$
{\BK}_p\LRf{{\SE}}= \pi_p\LRf{{\BK}{{\SE}}}\qquad \forall p\in {\BZ}.
$$
\item Define the non connective $GW$-theory spectrum as 
$$
\g{W}\LRf{{\SE}_d}=\LRs{
\g{W}_n\LRf{{\SE}_d}: n\geq 0
},\quad {\rm where}\quad \g{W}_n\LRf{{\SE}_d}=\g{W}\LRf{\widetilde{{\BS}_{GW}^n\LRf{{\SE}_d}}}
$$
The $\g{W}$-groups are defined by 
$$
\g{W}_p\LRf{{\SE}_d}= \pi_p\LRf{\g{W}{{\SE}_d}}\qquad \forall p\in {\BZ}.
$$

\item It follows from (\ref{CotBle2289}, \ref{DonFibrs}) that both ${\BK}\LRf{{\SE}}\in \Sp$ and 
$\g{W}\LRf{{\SE}_d}\in \Sp$ are $\Omega$-spectra. We we have defined two functors:
$$
\left\{\begin{array}{l}
{\BK}: {\eCat} \lra \Sp\\
{\g}W: {\bf eCat}D \lra \Sp\\
\end{array}\right. 
$$
\eE 
}
\eD 
Thus far, in this section, all the above concerns only exact categories.

\subsection{$G{\CW}$ spectrum of dg categories} \label{dgCatSECGWsp}
Thomason \cite{TT90}  started considering $K$-theory of the category of chain complexes of exact categories. This led to consideration of a variety of concepts like BiWaldhausen categories, complicial exact categories, and dg categories. 
All these culminated to the $GW$-theory of dg categories in \cite{S17}. It is worth pointing out that $GW$-theory encapsulates, both $K$-theory and Witt theory. For a dg category
${\SA}$ the $K$-theory of ${\SA}$ is given by the $GW$-theory of the hyperbolic category ${\CH}\LRf{{\SA}}$. 


Given an exact category $({\SE}, ^{\vee},  \varpi)$ with duality $^{\vee}$ and double dual identification $\varpi$, there is another 
exact category $({\SE}, ^{\vee},  -\varpi)$, usually referred to as the skew duality. So,  any theory for $({\SE}, ^{\vee},  \varpi)$, automatically produces a theory for $({\SE}, ^{\vee},  -\varpi)$. After the work of Thomason \cite{TT90}, preceded by work of Walhausen \cite{W83}, attention shifted to the $GW$-Theory of the category  of bounded chain complexes $Ch^b({\SE})$, for exact categories $({\SE}, ^{\vee}, \varpi)$ with duality. Because of the shift functor, one has one $GW$-theory  $GW^{[n]}(Ch^b({\SE}))$, for each 
$n\in {\BZ}$. However, these theories are four periodic 
$GW^{[n]}(Ch^b({\SE}))=GW^{[n+4]}(Ch^b({\SE}))$ for all $n\in {\BZ}$. So, we obtain four $GW$-theory of 
$Ch^b({\SE})$, while we have only two $GW$-Theory for $({\SE}, ^{\vee}, \varpi)$. 

The category $Ch^b({\SE})$ of exact categories $\SE$ can be endowed with the structure of a dg category, denoted by ${\bf dg}{\SE}$. This provides justification for
 developing $K$-theory and 
$GW$-theory of dg categories with weak equivalence and duality. Unless stated otherwise, we work with dg categories over $k={\BZ}(1/2)$, whenever we work on $GW$-theory. 
In addition to the above notations, let
${\bf dgCat}WD$ denote the category of small  pointed dg $k$-categories ${\SA}:=\LRf{{\SA}, \q, ^{\vee}, \varpi}$ with weak equivalences $\q$, 
 duality $^{\vee}$ and double dual identification $\varpi$.  In our case of ${\bf dg}{\SE}$,  the quasi isomorphisms are defined to be  the weak equivalences. 

For ${\SA}:=\LRf{{\SA},  \q, ^{\vee},\varpi}\in {\bf dgCat}WD$, in \cite{S17}, a spectrum $G{\CW}({\SA})\in \Sp$ was associated in the category $\Sp$ of spectrum
of pointed spaces. So, this is a functor
$$
G{\CW}:=\{G{\CW}(-)_n: n\geq 0 \}: {\bf dgCat}WD \lra \Sp
$$
In fact,  $G{\CW}({\SA})_n$ was defined as the geometric realization of an explicitly defined simplicial category. In fact, $G{\CW}\LRf{\SA}$ is a positive $\Omega$-spectrum, which means $G{\CW}\LRf{\SA}_n\cong \Omega G{\CW}\LRf{\SA}_{n+1}$,
$\forall n\geq 1$. 
However,  for ${\SE}\in {\bf eCat}D$, $G{\CW}\LRf{{\bf dg}{\SE}}$ actually a $\Omega$-spectrum, and 
$$
{\bf GW}\LRf{{\SE}}\cong
{\bf GW}\LRf{{\bf dg}{\SE}}\cong \Omega G{\CW}\LRf{{\bf dg}{\SE}}_{1}
$$ 
({\it See, \cite[Thm 10.1.17, pp 431]{M23}, \cite[Eqn 10.31, pp 449]{M23}).}
Now, we have the following agreement proposition. 

\bP[Agreement]\label{AghrThm2GG}{\rm 
With notations as above, for ${\SE}:=\LRf{{\SE}, ^{\vee}, \varpi}\in {\eCat}D$, with $1/2\in {\SE}$, there is a
\begin{equation}\label{3531Homennnh}
{\g}W({\SE}) \cong G{\CW}\LRf{{\bf dg}{\SE}}\quad homotopy~equivalence ~in ~\Sp.
\end{equation} 
 ({\it We caution that left side is a $\Omega$ spectrum and right side is a positive $\Omega$ spectrum.})
}
\eP
\pf It follows from uniqueness of $\Omega$ spectrum. Recall the formula \cite[pp 107]{W78}, 
$\LRt{{\BS}^1\wedge X, Y}\cong \LRt{X, \Omega Y}$, for $X, Y\in {\bf Top}_{\bullet}$. \pic $\eop$

\subsection{Karoubi ${\BG}W$ bispectrum}\label{subSecKaroBISP}
For ${\SE}\in {\eCat}D$, we provided two descriptions (\ref{3531Homennnh}) of the $GW$
spectrum ${\g}W({\SE}) \cong G{\CW}\LRf{{\bf dg}{\SE}}\in {\Sp}$. In this section, for 
${\SA}=\LRf{{\SA}, \q, ^{\vee}, \varpi}\in {\dgCat}WD$, we associate the Karaiubi Grothendieck-Witt  bispectrum ${\BG}W\LRf{{\SA}}$, which takes values in the category 
$\BiSp:= \BiSp\LRf{{\BS}^1-\tilde{\BS}^1_{\Sp}}$ of 
 ${\BS}^1-\tilde{{\BS}}^1_{\Sp}$-bispectrum. 
 Here ${\BS}^1:=\frac{\LRt{0, 1}}{0\sim 1}$ denotes the $1$-sphere.
 
 We refer to \cite[pp 458-]{M23}, \cite{S17} for the background on the category $\BiSp$, including definitions of symmetric $K$-spectrum, ${\Sp}\LRf{{\SC}, K}$ for  monoidal model categories ${\SC}:=\LRf{{\SC},\otimes,  {\mathbbm 1} }$. Among the original sources of all these may be the work of Hovey \cite{H01}.
We recall some of the jargon and notations:
\bE
\item We used the notation ${\Sp}$ for the category of spectrum of pointed topological spaces. Denote the subcategory of symmetric spectra, by $sym{\Sp}$. 
\item We let $\widetilde{\Sp}:=\LRf{sym{\Sp}, \wedge_{{\BS}}, {\BS}}$ denote the category symmetric spectrum, where ${\BS}=\LRs{\wedge^n{\BS}^1:n\geq 0}$ is the ring spectrum, which acts as the unit.
With ${\SC}=\LRf{{\bf Top}_{\bullet}, \wedge, {\BS}^0}$, and  $K={\BS}^1$ we have 
$\widetilde{\Sp}:={\Sp}\LRf{{\SC}, K}$ (see \cite[pp 458]{M23}). Note $\widetilde{\Sp}$ is a symmetric monoidal model category \cite[pp 460]{M23}.
\item Now, we repeat the same process, starting with 
${\SC}:=\widetilde{\Sp}:=\LRf{sym{\Sp}, \wedge_{{\BS}}, {\BS}}$. Let ${\BS}_{{\Sp}}^1=
\LRs{{\BS}^1, {\BS}^2, {\BS}^3, \ldots}\in sym{\Sp}$ be the sphere spectrum, and let 
$K=\widetilde{{\BS}_{{\Sp}}^1}\in sym{\Sp}$ be the cofibrant replacement of 
${\BS}_{{\Sp}}^1$. This means $K$ is cofibrant and $K\lra {\BS}_{{\Sp}}^1$ is trivial fibration. Now define the monoidal model category of bispctra (see \cite[pp 461]{M23})
$$
\BiSp:= \BiSp\LRf{{\BS}^1-\tilde{\BS}^1_{\Sp}}:={\Sp}\LRf{\widetilde{{\Sp}}, K}
$$
So, the objects of $\BiSp$ are symmetric $\widetilde{{\Sp}}$-spectrum, $\LRs{{\CX}_0, {\CX}_1, \ldots }$, where ${\CX}_i\in sym{\Sp}$. 
\eE 

Write $K=\tilde{\BS}^1_{\Sp}$. Given ${\bf X}=\LRs{{\CX}_0, {\CX}_1, \ldots }\in \BiSp$ following homotopy groups are defined:
\begin{equation}\label{BiGroips}
\left\{\begin{array}{ll}
\pi_n\LRf{{\bfX}}= \lim_k \LRt{K^{n+k}, {\CX}_k}_{Ho\Sp^{\Sigma}} & Naive~homotopy\\
\pi_n^t\LRf{{\bfX}}= \lim_k \LRt{\Sigma K^{n}, {\bfX}}_{Ho\BiSp} & True~homotopy\\
\end{array}\right.
\qquad \forall n\in {\BZ}
\end{equation}
where $\Sigma K^n=\LRs{K^{n+k}: k=0, 1, 2, \ldots} \in \BiSp$. 
Under suitable conditions, these two groups coincide \cite[pp 1840]{S17}.

%

Now we proceed to define Karoubi ${\BG}W$ bispectrum. 
\bD{\rm 
Let $k$ be a commutative ring, with $1/2\in k$. 
Let ${\dgCat}WD_k$ denote the category of dg $k$-categories with weak equivalences and duality. 
$$
\left\{\begin{array}{ll}
{\rm A~cone~functor} & {\bf C}: {\bf dg}CatWD_k \lra {\bf dg}CatWD_k~and\\
{\rm a~suspension~functor} & \s: {\bf dg}CatWD_k \lra {\bf dg}CatWD_k\\
\end{array}\right.
$$
were constructed, on the category ${\bf dg}CatWD_k$.
For $\LRf{{\SA}, ^{\vee}, \q, \varpi}\in {\bf dg}CatWD_k$, and integer $n\geq 0$, 
one defines 
\begin{equation}\label{382Deflll}
{\BG}W({\SA})_n:=G{\CW}\LRf{\s^n{\SA}}\in \widetilde{\Sp}
\end{equation}
 in the category 
$ \widetilde{\Sp}$ of symmetric spectra. The bispectrum is defined, as a functor
$$
{\BG}W:{\bf dgCat}WD_k\lra \BiSp \quad sending \quad 
{\SA} \mapsto {\BG}W({\SA}):=\LRs{
{\BG}W({\SA})_n: n\geq 0
}.
$$
Consequently, for ${\SE}\in {\bf eCat}D_k$, one can associate
$$
{\SE} \mapsto {\BG}W\LRf{{\bf dg}{\SE}}\in \BiSp.
\quad {\rm Note}\quad {\BG}W\LRf{{\bf dg}{\SE}}=\LRs{G{\CW}\LRf{\s^n{\bf dg}{\SE}}:n\geq 0}
$$
}
\eD
Following are some observation on $GW$ invariants of exact categories with duality.
\begin{remark}\label{GW444Ed}{\rm 
Let ${\SE}=\LRf{{\SE}, ^{\vee}, \varpi}\in \eCat$. 
\bE
\item The $\g W\LRf{{\SE}} \in \widetilde{\Sp}=\LRf{sym{\Sp}, \wedge_{{\BS}}, {\BS}}$. This follows from the fact the spectrum $\g W\LRf{{\SE}}$ is defined by iterated suspension. 
\item ${\SE}$ can be considered as a dg category with weak equivalence (isomorphisms) and duality, concentrated at degree zero \cite[pp 435]{M23}.
\item Since $\s$ is ring,  it is easy to see that 
$\s^n\otimes {\SE}$ is defined, and is an exact category with duality, where split exact sequences are defined to be exact. Define 
$$
{\G}W\LRf{{\SE}}_n= \g{W}\LRf{\s^n\otimes {\SE}}, \quad and \quad 
{\G}W\LRf{{\SE}}=\LRs{{\G}W\LRf{{\SE}}_n: n\geq 0}
$$
\eE 
}
\end{remark}
The following reconciles two bispectrum of exact categories.
\bP{\rm 
Let ${\SE}=\LRf{{\SE}, ^{\vee}, \varpi}\in \eCat$. 
Then 
\bE
\item ${\G}W\LRf{{\SE}}\in \BiSp$ is a bispectrum.
\item There is zig zag stable equivalence 
$$
\diagram
{\G}W\LRf{{\SE}}\ar[r]^{\sim} & {\BG}W\LRf{{\bf dg}{\SE}}
\enddiagram
\qquad in \quad \BiSp. 
$$

\eE
}
\eP
\pf The proof is elementary, which we skip. Main point is that an exact category ${\SE}$ can be considered as a dg category, concentrated at degree zero. We plan to write a complete exposition somewhere else \cite{M24p}. $\eop$



 %
We summarize some of the above notations and introduce some more, for $n\geq 0$, as follows.
 \begin{equation}\label{DefbhseExc}
 \left\{\begin{array}{l}
 {\G}W({\SE})_n:={\g}W\LRf{\s^n{\SE}}\in sym\Sp, \quad and\\
    {\G}W({\SE}):=\LRs{ {\G}W({\SE})_n: n\geq 0}\in \BiSp\\
     {\BG}W({\SE})_n:= {\G}W({\SE})_n, ~{\BG}W({\SE})^+:= {\G}W({\SE})\in \BiSp\\
     {\BG}W({\SE})^{-}:={\BG}W\LRf{{\SE}, ^{\vee}, -\varpi}\in \BiSp\\
    \end{array}\right. 
 \end{equation}

 %
It also follows that 
$$
\left\{\begin{array}{l}
{\BG}W({\SE}) \cong{\BG}W^{[4n]}({\bf dg}{\SE})\\
{\BG}W({\SE})^- \cong {\BG}W^{[4n+2]}({\bf dg}{\SE})\\
\end{array}\right.
\qquad \forall n\in {\BZ}.
$$
While we obtain two bispectrum ${\BG}W({\SE})$, ${\BG}W({\SE})^-$ directly from an exact category, ${\SE}=\LRf{{\SE}, ^{\vee}, \varpi}$, the shifted bispectrum  
${\BG}W^{[n]}({\bf dg}{\SE})\cong {\BG}W^{[n+4]}({\bf dg}{\SE})$ is four periodic. 
Following (\ref{BiGroips}),
define the Karoubi ${\BG}W$-groups
$$
\left\{\begin{array}{ll}
{\BG}W_p\LRf{{\SA}}= \pi_p^t\LRf{{\BG}W\LRf{{\SA}}} & \forall{\SA}\in {\dgCat}WD_k \\
{\BG}W_p\LRf{{\SE}}= \pi_p^t\LRf{{\BG}W\LRf{{\SE}}} \cong
 \pi_p^t\LRf{{\BG}W^{[4n]}\LRf{{\bf dg}{\SE}}} & \forall{\SE}\in {\eCat}D \\
{\BG}W_p\LRf{{\SE}}^{-}= \pi_p^t\LRf{{\BG}W\LRf{{\SE}}^{-}} \cong \pi_p^t\LRf{{\BG}W^{[4n+2]}\LRf{{\bf dg}{\SE}}} & \forall{\SE}\in {\eCat}D \\
\end{array}\right.
\quad \forall p, n\in {\BZ}.
$$
More explicitly, we have the following.
\bC\label{Coolastgr}{\rm 
Let ${\SC}=\LRf{{\SC}, ^{\vee}, \varpi}\in {\bf eCat}D_{\BZ}$. As usual, denote
$\forall k\in {\BZ},~{\BG}W_k({\SE})=\pi_k\LRf{{\BG}W_k({\SE})}$ denotes the homotopy group. 
Then 
$$
{\BG}W_p({\SE}) =\left\{\begin{array}{ll}
{\bf GW}_p\LRf{{\SE}}& p\geq 1\\
{\bf GW}_0\LRf{\widetilde{{\bf D}^b\LRf{{\SE}}}}& p=0\\
{\bf GW}_0\LRf{\widetilde{{\bf D}^b\LRf{\s^{-p}{\SE}}}}& p\leq -1\\
\end{array}\right.
$$
Here $\forall k\geq 0, ~{\bf GW}_k\LRf{{\SE}}=\pi_k\LRf{{\bf GW}\LRf{{\SE}}}$ denote the $GW$-groups,
while ${\bf GW}\LRf{{\SE}}$ denotes the ${GW}$-space. 
}
\eC
\pf It follows directly from \cite[pp 1799]{S17}, \cite{M23}.

\section{Notations}\label{SecVorAus}


In this section we establish some standard, and some not-so-standard notations. We freely use some of the notations from \cite{H}, which I consider standard by now. 

\begin{notation}\label{nota}{\rm 
Let $X$  be a quasi projective scheme over a noetherian affine scheme $A$, and 
$Z\subseteq X$ be a closed subset and $U=X-Z$. We assume that $grade(Z)=k$. In the context of ${\BG}W$-theory, ${\SL}$ will denote a fixed invertible sheaf. As is standard, 
$QCoh(X)$, $Coh(X)$, ${\SV}(X)$, respectively, denote the category of Quasi Coherent,
Coherent, locally free sheaves over
$X$.
$$
\left\{\begin{array}{l}
Denote~d=\dim X, ~k=grade Z\quad 
{\rm (see~Definition~\ref{GradCloQA})}\\ 
Coh^Z(X)=\{M\in Coh(X): \supp{M} \subseteq Z \}\\
\forall M\in Coh(X), ~{\rm denote}~\PDV{{\SF}}:=projective~(locally~free)~dimension\\
{\BM}(X)=\{M\in Coh(X): \PDV(M) < \infty \}\\
{\BM}^Z(X)=\left\{M\in  {\BM}(X): \supp{M} \subseteq Z\right\}\\
C{\BM}^Z(X)=\left\{M\in  {\BM}^Z(X): {\PDV}(M)=k \right\}\quad {\rm see~(\ref{defPerModd})}\\
{\rm Usually,} ~{\SE}(X)\subseteq QCoh(X)~{\rm  would ~denote ~an~exact~subcategory}\\
{\rm Examples:~}{\SE}(X)=Coh(X), Coh^Z(X), {\BM}(X),  C{\BM}^Z(X), ~{and~others} \\
\end{array}\right.
$$
We continue:
$$
\left\{\begin{array}{l}
{\rm For~such~an~exact~category,}\\
Ch^b({\SE}(X))~{\rm denotes~ the ~category ~of~ bounded ~ complexes, ~with~objects~in~}{\SE}(X)\\
Ch^b_Z({\SE}(X))=\left\{{\SF}_{\bullet}\in Ch^b({\SE}(X)): \left({\SF}_{\bullet}\right)_{|X-Z}~is~exact \right\}\quad
{\rm the~full~subcategory.}\\
{\bf dg}{\SE}(X)~{\rm denotes~the~dg~category~corresponding ~to~} Ch^b({\SE}(X))\\
{\bf dg}_Z{\SE}(X)~{\rm denotes~the~dg~category~corresponding ~to~} Ch^b_Z({\SE}(X))\\
{\bfD}^b({\SE}(X))={\rm Bounded ~Derived~Category~of~exact~categories}~{\SE}(X)\\
{\SD}^b_Z({\SE}(X))=\left\{{\SF}_{\bullet}\in {\bfD}^b({\SE}(X)): \left({\SF}_{\bullet}\right)_{|X-Z}~is~exact \right\}\quad
{\rm the~full~subcategory.}\\
\end{array}\right.
$$

For future use, we extend the list (\ref{nota}), as follows
$$
\left\{\begin{array}{l}
{\bf P}erf(X)\subseteq Ch(QCoh(X))~{\rm denote~the~category~of~perfect~complexes}\\
{\bf P}erf_Z(X)=\LRs{
M_{\bullet}\in {\bf P}erf(X): \LRf{
{\CH}_n\LRf{M_{\bullet}}
}_{|U}=0~\forall~n
} ~{\rm ~be~the~full~subcategory.}\\
\end{array}\right.
$$

}
\end{notation}
For a larger part of this article, we work with the following setup.
\begin{notation}\label{setUP}{\rm 
Let $A$ be a noetherian commutative ring. Let $X$ be a quasi projective scheme over $\spec{A}$. 
We will be working under the following setup.
\bE
\item Let $Y=\proj{S}$ for some graded ring $S=A[x_0, x_1, \ldots, x_N]$ with $\deg(x_i)=1$.
\item We assume $X\subseteq Y$ is an open subset. Assume $\dim X=d$. 
\eE 
}
\end{notation}

\printindex

\end{document}